# Optimal reduction

Juan–Pablo Ortega
Institut Nonlinéaire de Nice
Centre National de la Recherche Scientifique
1361, route des Lucioles
F-06560 Valbonne, France
Juan-Pablo.Ortega@inln.cnrs.fr

27 June 2002

**Abstract**

We generalize various symplectic reduction techniques of Marsden, Weinstein, Sjamaar, Bates, Lerman, Marle, Kazhdan, Kostant, and Sternberg to the context of the optimal momentum map. We see that, even though all those reduction procedures had been designed to deal with canonical actions on symplectic manifolds in the presence of a momentum map, our construction allows the construction of symplectic point and orbit reduced spaces purely within the Poisson category under hypotheses that do not necessarily imply the existence of a momentum map.

We construct an orbit reduction procedure for canonical actions on a Poisson manifold that exhibits an interesting interplay with the von Neumann condition previously introduced by the author in his study of the theory of singular dual pairs. More specifically, this condition ensures that the orbits in the momentum space of the optimal momentum map (we call them polar reduced spaces) admit a presymplectic structure that generalizes the Kostant–Kirillov–Souriau symplectic structure of the coadjoint orbits in the dual of a Lie algebra. Using this presymplectic structure, the optimal orbit reduced spaces are symplectic with a form that satisfies a relation identical to the classical one obtained by Marle, Kazhdan, Kostant, and Sternberg for free Hamiltonian actions on a symplectic manifold. In the Poisson case we provide a sufficient condition for the polar reduced spaces to be symplectic. In the symplectic case the polar reduced spaces are symplectic if and only if certain relation between the tangent space to the orbit and its symplectic orthogonal with the tangent space to the isotropy type submanifolds is satisfied. In general, the presymplectic polar reduced spaces are foliated by symplectic submanifolds that are obtained through a generalization to the optimal context of the so called Sjamaar Principle, already existing in the theory of Hamiltonian singular reduction. We call these subspaces the regularized polar reduced spaces.

We use these ideas to shed some light in the problem of orbit reduction of globally Hamiltonian actions when the symmetry group is not compact and in the construction of a family of presymplectic homogeneous manifolds and of its symplectic foliation.

We also show that these reduction techniques can be implemented in stages whenever we are in the presence of certain hypotheses that generalize those already existing for free globally Hamiltonian actions.

# Contents

# 1 Introduction

Let $(M, \{\cdot,\cdot\})$ be a Poisson manifold and $G$ be a Lie group that acts properly on $M$ by Poisson diffeomorphisms via the left action $\Phi : G \times M \to M$. The group of Poisson transformations associated to this action will be denoted by $A_G := \{\Phi_g \mid g \in G\}$ and the canonical projection of $M$ onto the orbit space by $\pi_{A_G} : M \to M/A_G = M/G$. We will denote by $\mathfrak{g}$ the Lie algebra of $G$ and by $\mathfrak{g}\cdot m := T_m(G\cdot m)$ the tangent space at the point $m$ of its $G$–orbit,

The use of the canonical symmetries of $M$ encoded in the action of the Lie group $G$ has been used in [MR86, OR98] to ***reduce*** the Poisson system $(M, \{\cdot,\cdot\})$ into a smaller one where the degeneracies induced by the group invariance have disappeared. When $M$ happens to be a symplectic manifold with symplectic form $\omega$ and the $G$–group action has a momentum map $\mathbf{J} : M \to \mathfrak{g}^*$ associated, the reduction procedure can be adapted to this category using the so called ***symplectic*** or ***Marsden–Weinstein*** ***reduction*** [MW74, SL91, ACG91, BL97, O98, CS01, OR02b]. For the last thirty years, Marsden–Weinstein reduction has been a major tool in the construction of new symplectic manifolds and in the study of mechanical systems with symmetry.

More recently, a new momentum map, we call it ***optimal momentum map***, has been introduced [OR02a]. This object is partially inspired by the distribution theoretical approach to the conservation laws induced by symmetry that one can find in the works of Cartan [C22]. The use of this tool allows the construction of symplectically reduced spaces purely within the Poisson category under hypotheses that do not necessarily imply the existence of a (standard or group valued [AMM98]) momentum map. For a proof of these facts please check with [O02a]. All along this paper we will refer to the construction of Marsden–Weinstein reduced spaces with the help of the optimal momentum map as ***optimal reduction***.

In this paper we will study three main topics:

**(i) Optimal orbit reduction:** in the classical theory of symplectic reduction there are two equivalent approaches to the construction of the symplectic quotients, namely, ***point*** [MW74] and ***orbit reduction*** [Mar76, KKS78]. The analog of point reduction in the optimal context has been carried out in [OR02a, O02a]. In the first part of this paper we will concentrate in the development of an optimal orbit reduction procedure.

**(ii) Symplectic leaves and polar reduced spaces:** when a canonical, proper, and free action of a connected Lie group $G$ on a symplectic manifold $(M, \omega)$ has an equivariant momentum map $\mathbf{J} : M \rightarrow \mathfrak{g}^*$ associated, the diagram

$$\begin{array}{ccc} & M & \\ {\scriptstyle \pi_{A_G}} \swarrow & & \searrow {\scriptstyle \mathbf{J}} \\ M/G & & \mathfrak{g}^*_{\mathbf{J}} := \mathbf{J}(M) \end{array}$$

forms a so called ***dual pair*** [Lie90, W83], that is, the maps $\pi_{A_G}$ and $\mathbf{J}$ have symplectically orthogonal fibers. It has been shown [W83, Bl01] that if $\mathbf{J}$ has connected fibers then there is a bijective correspondence between the symplectic leaves of $M/G$, namely the symplectic orbit reduced spaces, and those of $\mathfrak{g}^*_{\mathbf{J}}$, that is, the coadjoint orbits inside $\mathbf{J}(M)$. The notion of dual pair has been generalized in [O02] to a context that allows the introduction of a similar diagram in the optimal context. In this framework it can also be formulated a symplectic leaf correspondence theorem that provides a bijection between the optimal orbit reduced spaces and the symplectic leaves of the space that in this case plays the role of the dual of the Lie algebra $\mathfrak{g}^*$. We will spell out the smooth and (pre)symplectic structures of these spaces that we will call ***polar reduced spaces***. These constructions generalize the standard results on the Lie–Poisson structures on the dual of Lie algebras and on the symplectic character of their coadjoint orbits, as discovered by Kostant, Kirillov, and Souriau. A condition introduced in [O02] under the name of *von Neumann condition* will play a very important role in this part.

**(iii) Optimal reduction by stages:** Suppose that we are in the same setup as in point **(ii)**. Let $N \subset G$ be a closed normal subgroup of $G$. The Reduction by Stages Theorem [MMPR98, MMOPR02] states that in such a situation we can carry out symplectic reduction in two shots: we first reduce by the $N$ action; the resulting reduced space inherits some symmetry properties from the quotient group $G/N$ that can be used to perform symplectic reduction one more time. The resulting reduced space is naturally symplectomorphic to the one–shot reduced space that one obtains by just using the $G$–action. We will see that this procedure can be reproduced in the optimal context without any hypothesis on the freeness of the action. As a byproduct we will obtain a Singular Reduction by Stages Theorem that will generalize the results in [MMPR98, MMOPR02] to the non free actions case.

# 2 The optimal momentum map

Unless the contrary is explicitly stated, all along this paper we will work on a Poisson manifold $(M, \{\cdot, \cdot\})$ acted canonically and properly upon by the Lie group $G$ via the left action $\Phi : G \times M \to M$. The group of Poisson transformations associated to this action will be denoted by $A_G := \{\Phi_g \mid g \in G\}$ and the projection of $M$ onto the orbit space by $\pi_{A_G} : M \to M/A_G = M/G$. Given a point $m \in M$ with isotropy subgroup $G_m$ we will denote by $\mathfrak{g} \cdot m := T_m(G \cdot m)$ the tangent space at $m$ to the $G$–orbit that goes through $m$. We recall that when the $G$–action is proper, the connected components of the subset $M_{G_m} := \{z \in M \mid G_z = G_m\}$ made out of the points that have the same isotropy as $m$, are smooth submanifolds that we will call ***isotropy type submanifolds***.

## 2.1 The optimal momentum map and the momentum space

The optimal momentum map was introduced in [OR02a] as a general method to find the conservation laws associated to the symmetries of a Poisson system. We briefly recall its definition and elementary properties, as presented in [OR02a]. Let $A'_G$ be the distribution on $M$ defined by the relation:

$$A'_G(m) := \{X_f(m) \mid f \in C^\infty(M)^G\}, \qquad \text{for all} \quad m \in M.$$

Depending on the context, the generalized distribution (in the sequel we will omit the adjective "generalized") $A'_G$ is called the $G$–***characteristic distribution*** or the ***polar distribution*** defined by $A_G$ [O02]. The distribution $A'_G$ is smooth and integrable in the sense of Stefan and Sussman [St74a, St74b, Su73], that is, given any point in $M$ there is an ***accessible set*** or ***maximal integral leaf*** of $A'_G$ going through it. We recall (see [OR02a] for the details) that if $M$ is actually a symplectic manifold with form $\omega$ then

$$A'_G(m) = (\mathfrak{g} \cdot m)^\omega \cap T_m M_{G_m}, \quad \text{for all} \quad m \in M. \tag{2.1}$$

Moreover, if the $G$–action has a standard momentum map $\mathbf{J} : M \to \mathfrak{g}^*$ associated then

$$A'_G(m) = \ker T_m\mathbf{J} \cap T_m M_{G_m}, \quad \text{for all} \quad m \in M. \tag{2.2}$$

The ***optimal momentum map*** $\mathcal{J}$ is defined as the canonical projection onto the leaf space of $A'_G$, that is,

$$\mathcal{J} : M \longrightarrow M/A'_G.$$

By its very definition, the fibers or levels sets of $\mathcal{J}$ are preserved by the Hamiltonian flows associated to $G$–invariant Hamiltonian functions and $\mathcal{J}$ is ***universal*** with respect to this property, that is, any other map whose level sets are preserved by $G$–equivariant Hamiltonian dynamics factors necessarily through $\mathcal{J}$.

The leaf space $M/A'_G$ is called the ***momentum space*** of $\mathcal{J}$. When considered as a topological space endowed with the quotient topology, it is easy to see [O02] that the optimal momentum map is ***continuous and open***.

The pair $(C^\infty(M/A'_G), \{\cdot,\cdot\}_{M/A'_G})$ is a Poisson algebra (the term Poisson variety is also frequently used) when we define

$$C^\infty(M/A'_G) := \{f \in C^0(M/A'_G) \mid f \circ \mathcal{J} \in C^\infty(M)\}, \tag{2.3}$$

and the bracket $\{\cdot,\cdot\}_{M/A'_G}$ defined by

$$\{f,\, g\}_{M/A'_G}(\mathcal{J}(m)) = \{f \circ \mathcal{J},\, g \circ \mathcal{J}\}(m), \tag{2.4}$$

for every $m \in M$ and $f, g \in C^\infty(M/A'_G)$. Note that as $\mathcal{J}$ is open and surjective then, for any real valued function $f$ on $M/A'_G$ such that $f \circ \mathcal{J} \in C^\infty(M)$ we have that $f \in C^0(M/A'_G)$ necessarily. Hence, in this case, the standard definition (2.3) can be rephrased by saying that $f \in C^\infty(M/A'_G)$ iff $f \circ \mathcal{J} \in C^\infty(M)$.

The $G$–action on $M$ naturally induces a smooth action $\Psi : G \times M/A'_G \to M/A'_G$ of $G$ on $M/A'_G$ defined by the expression $\Psi(g, \mathcal{J}(m)) := \mathcal{J}(g \cdot m)$ with respect to which $\mathcal{J}$ is $G$–equivariant. We recall that the term smooth in this context means that $\Psi^* C^\infty(M/A'_G) \subset C^\infty(G \times M/A'_G)$. Notice that since $M/A'_G$ is not Hausdorff in general, there is no guarantee that the isotropy subgroups $G_\rho$ of elements $\rho \in M/A'_G$ are closed, and therefore embedded, subgroups of $G$. Also, even if $G$ is connected $G_\rho$ does not need to be connected (see example in [O02a]). However, there is still something that we can say:

**Proposition 2.1** *Let $G_\rho$ be the isotropy subgroup of the element $\rho \in M/A'_G$ associated to the $G$–action on $M/A'_G$ that we just defined. Then:*

**(i)** *There is a unique smooth structure on $G_\rho$ with respect to which this subgroup is an initial (see below) Lie subgroup of $G$ with Lie algebra $\mathfrak{g}_\rho$ given by*

$$\mathfrak{g}_\rho = \{\xi \in \mathfrak{g} \mid \xi_M(m) \in T_m\mathcal{J}^{-1}(\rho), \text{ for all } m \in \mathcal{J}^{-1}(\rho)\} \tag{2.5}$$

*or, equivalently*

$$\mathfrak{g}_\rho = \{\xi \in \mathfrak{g} \mid \exp t\xi \in G_\rho, \text{ for all } t \in \mathbb{R}\}. \tag{2.6}$$

**(ii)** *With this smooth structure for $G_\rho$, the left action $\Phi^\rho : G_\rho \times \mathcal{J}^{-1}(\rho) \to \mathcal{J}^{-1}(\rho)$ defined by $\Phi^\rho(g, z) := \Phi(g, z)$ is smooth.*

**(iii)** *This action has fixed isotropies, that is, if $z \in \mathcal{J}^{-1}(\rho)$ then $(G_\rho)_z = G_z$, and $G_m = G_z$ for all $m \in \mathcal{J}^{-1}(\rho)$.*

**(iv)** *Let $z \in \mathcal{J}^{-1}(\rho)$ arbitrary. Then,*

$$\mathfrak{g}_\rho \cdot z = A'_G(z) \cap \mathfrak{g} \cdot z = T_z\mathcal{J}^{-1}(\rho) \cap \mathfrak{g} \cdot z. \tag{2.7}$$

**Proof.** For **(i)** through **(iii)** check with [O02a]. We prove **(iv)**: the inclusion $\mathfrak{g}_\rho \cdot z \subset A'_G(z) \cap \mathfrak{g} \cdot z$ is a consequence of (2.5). Conversely, let $X_f(z) = \xi_M(z) \in A'_G(z) \cap \mathfrak{g} \cdot z$, with $f \in C^\infty(M)^G$ and $\xi \in \mathfrak{g}$. The $G$–invariance of the function $f$ implies that $[X_f, \xi_M] = 0$, and hence, if $F_t$ is the flow of $X_f$ and $G_t$ is the flow of $\xi_M$ (more explicitly $G_t(m) = \exp t\xi \cdot m$ for any $m \in M$), then $F_t \circ G_s = G_s \circ F_t$. By one of the Trotter product formulas (see [AMR99, Corollary 4.1.27]), the flow $H_t$ of $X_f - \xi_M$ is given by

$$H_t(m) = \lim_{n\to\infty} \left(F_{t/n} \circ G_{-t/n}\right)^n (m) = \lim_{n\to\infty} \left(F^n_{t/n} \circ G^n_{-t/n}\right)(m) = (F_t \circ G_{-t})(m) = F_t(\exp -t\xi \cdot m),$$

for any $m \in M$. Consequently, as $X_f(z) = \xi_M(z)$, the point $z \in M$ is an equilibrium of $X_f - \xi_M$, hence $F_t(\exp -t\xi \cdot z) = z$ or, analogously $\exp t\xi \cdot z = F_t(z)$. Applying $\mathcal{J}$ on both sides of this equality, and taking into account that $F_t$ is the flow of a $G$–invariant Hamiltonian vector field, it follows that $\exp t\xi \cdot \rho = \rho$, and hence $\xi \in \mathfrak{g}_\rho$ by (2.6) . Thus $\xi_M(z) \in \mathfrak{g}_\rho \cdot z$, as required. ■

Recall that we say that $N$ is an ***initial submanifold*** of the smooth manifold $M$ when the inclusion $i : N \to M$ is a smooth immersion that satisfies that for any manifold $Z$, a mapping $f : Z \to N$ is smooth iff $i \circ f : Z \to M$ is smooth. The initial submanifold structure is unique in the sense that if $N$ admits another smooth structure, call it $N'$, that makes it into an initial submanifold of $M$, then the identity map $id_N : N \to N'$ is a diffeomorphism. Indeed, as the injection $N \hookrightarrow M$ is smooth and $N'$ is by hypothesis initial then, the identity map $id_N : N \to N'$ is smooth. As the same argument can be made for $id_{N'} : N' \to N$, the result follows.

We finish this section by emphasizing that the structure of the momentum space $M/A'_G$ may become very intricate. The following example shows that even when the $G$–action is very simple and the corresponding orbit space $M/G = M/A_G$ is a quotient regular manifold, the associated momentum space $M/A'_G$ does not need to share those properties.

**Example 2.2** Let $M := \mathbb{T}^2 \times \mathbb{T}^2$ be the product of two tori whose elements we will denote by the four–tuples $(e^{i\theta_1}, e^{i\theta_2}, e^{i\psi_1}, e^{i\psi_2})$. We endow $M$ with the symplectic structure $\omega$ defined by $\omega := \mathbf{d}\theta_1 \wedge \mathbf{d}\theta_2 + \sqrt{2}\,\mathbf{d}\psi_1 \wedge \mathbf{d}\psi_2$. We now consider the canonical circle action given by $e^{i\phi} \cdot (e^{i\theta_1}, e^{i\theta_2}, e^{i\psi_1}, e^{i\psi_2}) := (e^{i(\theta_1+\phi)}, e^{i\theta_2}, e^{i(\psi_1+\phi)}, e^{i\psi_2})$. First of all, notice that since the circle is compact and acts freely on $M$, the corresponding orbit space $M/A_{S^1}$ is a smooth manifold such that the projection $\pi_{A_{S^1}} : M \to M/A_{S^1}$ is a surjective submersion. The polar distribution $A'_{S^1}$ does not have that property. Indeed, $C^\infty(M)^{S^1}$ comprises all the functions $f$ of the form $f \equiv f(e^{i\theta_2}, e^{i\psi_2}, e^{i(\theta_1-\psi_1)})$. An inspection of the Hamiltonian flows associated to such functions readily shows that the leaves of $A'_{S^1}$ fill densely the manifold $M$ and that the leaf space $M/A'_{S^1}$ can be identified with the leaf space $\mathbb{T}^2/\mathbb{R}$ of a Kronecker (irrational) foliation of a two–torus $\mathbb{T}^2$. Under these circumstances $M/A'_{S^1}$ cannot possibly be a regular quotient manifold. ♦

## 2.2 The level sets of the momentum map and the associated isotropies

By construction, the fibers of $\mathcal{J}$ are the leaves of an integrable generalized distribution and thereby ***initial immersed submanifolds*** of $M$ [Daz85]. We summarize this and other elementary properties of the fibers of $\mathcal{J}$ in the following proposition.

**Proposition 2.3** *Let $(M, \{\cdot,\cdot\})$ be a Poisson manifold and $G$ be a Lie group that acts properly and canonically on $M$. Let $\mathcal{J} : M \to M/A'_G$ be the associated optimal momentum map. Then for any $\rho \in M/A'_G$ we have that:*

**(i)** *The level set $\mathcal{J}^{-1}(\rho)$ is an immersed initial submanifold of $M$.*

**(ii)** *There is a unique symplectic leaf $\mathcal{L}$ of $(M, \{\cdot,\cdot\})$ such that $\mathcal{J}^{-1}(\rho) \subset \mathcal{L}$.*

**(iii)** *Let $m \in M$ be an arbitrary element of $\mathcal{J}^{-1}(\rho)$. Then, $\mathcal{J}^{-1}(\rho) \subset M_{G_m}$, with $M_{G_m} := \{z \in M \mid G_z = G_m\}$.*

In the sequel we will denote by $\mathcal{L}_\rho$ the unique symplectic leaf of $M$ that contains $\mathcal{J}^{-1}(\rho)$. Notice that as $\mathcal{L}_\rho$ is also an immersed initial submanifold of $M$, the injection $i_{\mathcal{L}_\rho} : \mathcal{J}^{-1}(\rho) \hookrightarrow \mathcal{L}_\rho$ is smooth.

From the point of view of the optimal momentum map the existence of a standard ($\mathfrak{g}^*$ or $G$–valued) momentum map can be seen as an integrability feature of the $G$–characteristic distribution that makes the fibers of $\mathcal{J}$ particularly well–behaved. Indeed, it can be proved that when $M$ is a symplectic manifold and the $G$–action has a standard momentum map associated then, the fibers $\mathcal{J}^{-1}(\rho)$ of the optimal momentum map are closed imbedded submanifolds of $M$. More generally, if $\mathcal{J}^{-1}(\rho)$ is closed as a subset of the isotropy type submanifold $M_H$ in which it is sitting, then (see [OR02a])

- $\mathcal{J}^{-1}(\rho)$ is a closed embedded submanifold of $M_H$ and therefore an embedded submanifold of $M$, and
- the isotropy subgroup $G_\rho$ of $\rho \in A'_G$ is a closed embedded Lie subgroup of $G$.

## 2.3 The dual pair associated to the optimal momentum map

We mentioned in the introduction that the standard momentum map can be used to construct a dual pair in the sense of Lie [Lie90] and Weinstein [W83]. The notions of duality and dual pair have been generalized in [O02] in such way that in many situations the optimal momentum map provides an example of these newly introduced dual pairs. We now briefly recall some of the notions introduced in [O02]. For the details and proofs of the following facts the reader is encouraged to check with the original paper.

**Definition 2.4** *Let $M$ be a smooth manifold. A* ***pseudogroup of transformations*** *or* ***pseudogroup of local diffeomorphisms*** *$A$ of $M$ is a set of local diffeomorphisms of $M$ that satisfy:*

**(i)** *Each $\phi \in A$ is a diffeomorphism of an open set (called the domain of $\phi$) of $M$ onto an open set (called the range of $\phi$) of $M$.*

**(ii)** *Let $U = \cup_{i \in I} U_i$, where each $U_i$ is an open set of $M$. A diffeomorphism $\phi$ of $U$ onto an open set of $M$ belongs to $A$ if and only if the restrictions of $\phi$ to each $U_i$ is in $A$.*

**(iii)** *For every open set $U$ of $M$, the identity transformation of $U$ is in $A$.*

**(iv)** *If $\phi \in A$, then $\phi^{-1} \in A$.*

**(v)** *If $\phi \in A$ is a diffeomorphism of $U$ onto $V$ and $\phi' \in A$ is a diffeomorphism of $U'$ onto $V'$ and $V \cap U'$ is nonempty, then the diffeomorphism $\phi' \circ \phi$ of $\phi^{-1}(V \cap U')$ onto $\phi'(V \cap U')$ is in $A$.*

*Let $A$ be a pseudogroup of transformations on a manifold $M$ and $\sim$ be the relation on $M$ defined by: for any $x, y \in M$, $x \sim y$ if and only if there exists $\phi \in A$ such that $y = \phi(x)$. The relation $\sim$ is an equivalence relation whose space of equivalence classes is denoted by $M/A$.*

*If $M$ is a Poisson manifold with Poisson bracket $\{\cdot,\cdot\}$, we say that a pseudogroup of transformations $A$ of $M$ is a* ***pseudogroup of local Poisson diffeomorphisms*** *when any diffeomorphism $\phi \in A$ of an open set $U$ of $M$ onto an open set $V$ of $M$ is also a Poisson map between $(U, \{\cdot,\cdot\}_U)$ and $(V, \{\cdot,\cdot\}_V)$. The symbols $\{\cdot,\cdot\}_U$ and $\{\cdot,\cdot\}_V$ denote the restrictions of the bracket $\{\cdot,\cdot\}$ to $U$ and $V$, respectively.*

**Definition 2.5** *Let $(M, \{\cdot,\cdot\})$ be a Poisson manifold and $A$ be a pseudogroup of local Poisson diffeomorphisms of $M$. Let $\mathcal{A}'$ be the set of Hamiltonian vector fields associated to all the elements of $C^\infty(U)^A$, for all the open $A$–invariant subsets $U$ of $M$, that is,*

$$\mathcal{A}' = \left\{ X_f \mid f \in C^\infty(U)^A, \text{ with } U \subset M \text{ open and } A\text{–invariant} \right\}. \tag{2.8}$$

*The distribution $A'$ associated to the family $\mathcal{A}'$ will be called the* ***polar distribution*** *defined by $A$ (or equivalently the* ***polar of*** *$A$). Any generating family of vector fields for $A'$ will be called a* ***polar family*** *of $A$. The family $\mathcal{A}'$ will be called the* ***standard polar family*** *of $A$. The pseudogroup of local Poisson diffeomorphisms constructed using finite compositions of flows of the vector fields in any polar family of $A$ will be referred to as a* ***polar pseudogroup*** *induced by $A$. The polar pseudogroup $G_{\mathcal{A}'}$ induced by the standard polar family $\mathcal{A}'$ will be called the* ***standard polar pseudogroup****.*

**Remark 2.6** We say that the pseudogroup $A$ has the ***extension property*** when any $A$–invariant function $f \in C^\infty(U)^A$ defined on any $A$–invariant open subset $U$ satisfies that: for any $z \in U$, there is a $A$–invariant open neighborhood $V \subset U$ of $z$ and a $A$–invariant smooth function $F \in C^\infty(M)^A$ such that $f|_V = F|_V$. If the pseudosubgroup $A$ has the extension property, there is a simpler polar family, we will call it $\mathcal{A}'_{ext}$, that can be used to generate $A'$, namely

$$\mathcal{A}'_{ext} = \left\{ X_f \mid f \in C^\infty(M)^A \right\}.$$

In particular, if $A = A_G$, that is, the Poisson diffeomorphism group associated to a proper canonical $G$–action, the extension property is always satisfied and hence $\mathcal{A}'_{ext} = A'_G$, the $G$–characteristic distribution.

**Definition 2.7** *Let $(M, \{\cdot,\cdot\})$ be a Poisson manifold and $A, B$ be two pseudogroups of local Poisson diffeomorphisms. We say that the diagram*

$$\begin{array}{ccc} & (M, \{\cdot,\cdot\}) & \\ \pi_A \swarrow & & \searrow \pi_B \\ (M/A, \{\cdot,\cdot\}_{M/A}) & & (M/B, \{\cdot,\cdot\}_{M/B}) \end{array}$$

*is a* ***dual pair*** *on $(M, \{\cdot,\cdot\})$ when the polar distributions $A'$ and $B'$ are integrable and they satisfy that*

$$M/A' = M/B \text{ and } M/B' = M/A. \tag{2.9}$$

We now focus on the dual pairs induced by the optimal momentum map. Hence, let $(M, \{\cdot,\cdot\})$ be a Poisson manifold acted canonically and properly upon by a Lie group $G$, $A_G$ be the associated group of canonical transformations and $\mathcal{J} : M \to M/A'_G$ be the optimal momentum map. A natural question to ask is when the diagram

$$\begin{array}{ccc} & (M, \{\cdot,\cdot\}) & \\ \pi_{A_G} \swarrow & & \searrow \mathcal{J} \\ (M/A_G, \{\cdot,\cdot\}_{M/A_G}) & & (M/A'_G, \{\cdot,\cdot\}_{M/A'_G}) \end{array}$$

is a dual pair in the sense of Definition 2.7. Obviously, in this case, condition (2.9) is satisfied iff the double polar $A''_G := (G_{\mathcal{A}'})'$ of $A_G$ is such that

$$M/A''_G = M/A_G.$$

Poisson subgroups satisfying this condition are referred to as ***von Neumann subgroups***. In our discussion on orbit reduction we will use a slightly less demanding condition, namely, we will need group actions such that

$$\mathfrak{g} \cdot z = T_z(A_G(z)) \subset A_G''(z), \quad \text{for all} \quad z \in M. \tag{2.10}$$

A group action that satisfies (2.10) is called ***weakly von Neumann***. Obviously, if $A_G$ is von Neumann it is weakly von Neumann. Given that $A_G''$ is spanned by Hamiltonian vector fields, the weak von Neumann condition (2.10) implies that for any $z \in M$ sitting in the symplectic leaf $\mathcal{L}_z$ we have that $\mathfrak{g} \cdot z \subset T_z\mathcal{L}_z$, in particular, if $G^0$ is the connected component of $G$ containing the identity, the orbit $G^0 \cdot z$ is contained in the symplectic leaf $\mathcal{L}_z$.

We say that the group $A_G$ is ***weakly Hamiltonian*** when for every element $g \in G$ and any $m \in M$ we can write $\Phi_g(m) = F^1_{t_1} \circ F^2_{t_2} \circ \cdots \circ F^k_{t_k}(m)$, with $F^i_{t_i}$ the flow of a Hamiltonian vector field $X_{h_i}$ associated to a function $h_i \in \left(C^\infty(M)^G\right)^c$ that centralizes the $G$–invariant functions on $M$. It is clear that connected Lie group actions that have a standard ($\mathfrak{g}^*$ or $G$–valued) momentum map associated are weakly Hamiltonian. The importance of this condition in relation to our dual pairs is linked to the fact that ***weakly Hamiltonian proper actions induce von Neumann subgroups***. For a proof of this fact and for other situations where the von Neumann condition is satisfied see [O02].

## 2.4 Dual pairs, reduced spaces, and symplectic leaves

Let $(M, \{\cdot,\cdot\})$ be a smooth Poisson manifold, $A$ be a subgroup of its Poisson diffeomorphism group, and $(M/A, \{\cdot,\cdot\}_{M/A})$ be the associated quotient Poisson variety. Let $V \subset M/A$ be an open subset of $M/A$ and $h \in C^\infty(V)$ be a smooth function defined on it. If we call $U := \pi_A^{-1}(V)$ then, the vector field $X_{h\circ\pi_A|_U}$ belongs to the standard polar family $\mathcal{A}'$ and therefore its flow $(F_t, \mathrm{Dom}(F_t))$ uniquely determines a local Poisson diffeomorphism $(\bar{F}_t, \pi_A(\mathrm{Dom}(F_t)))$ of $M/A$. We will say that $(\bar{F}_t, \pi_A(\mathrm{Dom}(F_t)))$ is the ***Hamiltonian flow*** associated to $h$. The symplectic leaves of $M/A$ will be defined as the accessible sets in this quotient by finite compositions of Hamiltonian flows. It is not clear how to define these flows by projection of $A$–equivariant flows when $A$ is a pseudogroup of local transformations of $M$, hence we will restrict in this section to the case in which $A$ is an actual group of Poisson diffeomorphisms.

**Definition 2.8** *Let $(M, \{\cdot,\cdot\})$ be a smooth Poisson manifold, $A$ be a subgroup of its Poisson diffeomorphism group, and $(M/A, \{\cdot,\cdot\}_{M/A})$ be the associated quotient Poisson variety. Given a point $[m]_A \in M/A$, the **symplectic leaf** $\mathcal{L}_{[m]_A}$ going through it is defined as the (path connected) set formed by all the points that can be reached from $[m]_A$ by applying to it a finite number of Hamiltonian flows associated to functions in $C^\infty(V)$, with $V \subset M/A$ any open subset of $M/A$, that is,*

$$\mathcal{L}_{[m]_A} := \{F^1_{t_1} \circ F^2_{t_2} \circ \cdots \circ F^k_{t_k}([m]_A) \mid k \in \mathbb{N},\ F_{t_i} \text{ flow of some } X_{h_i},\ h_i \in C^\infty(V),\ V \subset M/A \text{ open}\}.$$

*The relation* being in the same symplectic leaf *determines an equivalence relation in $M/A$ whose corresponding space of equivalence classes will be denoted by $(M/A)/\{\cdot,\cdot\}_{M/A}$.*

**Theorem 2.9 (Symplectic leaves correspondence)** *Let $(M, \{\cdot,\cdot\})$ be a smooth Poisson manifold, $A, B$ be two groups of Poisson diffeomorphisms of $M$, and $G_{\mathcal{A}'}, G_{\mathcal{B}'}$ be the standard polar pseudogroups. If we denote by $(M/A)/\{\cdot,\cdot\}_{M/A}$ and $(M/B)/\{\cdot,\cdot\}_{M/B}$ the space of symplectic leaves of the Poisson varieties $(M/A, \{\cdot,\cdot\}_{M/A})$ and $(M/B, \{\cdot,\cdot\}_{M/B})$, respectively, we have that:*

**(i)** *The symplectic leaves of $M/A$ and $M/B$ are given by the orbits of the $G_{\mathcal{A}'}$ and $G_{\mathcal{B}'}$ actions on $M/A$ and $M/B$, respectively. As a consequence of this statement, we can write that*

$$(M/A)/\{\cdot,\cdot\}_{M/A} = (M/A)/G_{\mathcal{A}'} \text{ and } (M/B)/\{\cdot,\cdot\}_{M/B} = (M/B)/G_{\mathcal{B}'}. \tag{2.11}$$

**(ii)** *If the diagram $(M/A, \{\cdot,\cdot\}_{M/A}) \overset{\pi_A}{\longleftarrow} (M, \{\cdot,\cdot\}) \overset{\pi_B}{\longrightarrow} (M/B, \{\cdot,\cdot\}_{M/B})$ is a dual pair then the map*

$$\begin{array}{ccc} (M/A)/\{\cdot,\cdot\}_{M/A} & \longrightarrow & (M/B)/\{\cdot,\cdot\}_{M/B} \\ \mathcal{L}_{[m]_A} & \longmapsto & \mathcal{L}_{[m]_B} \end{array} \tag{2.12}$$

*is a bijection. The symbols $\mathcal{L}_{[m]_A}$ and $\mathcal{L}_{[m]_B}$ denote the symplectic leaves in $M/A$ and $M/B$, respectively, going through the point $[m]_A$ and $[m]_B$.*

One of our goals in the following pages will consist of describing the symplectic leaves of the Poisson varieties in the legs of the diagram $(M/G, \{\cdot,\cdot\}_{M/A_G}) \overset{\pi_{A_G}}{\leftarrow} (M, \{\cdot,\cdot\}) \overset{\mathcal{J}}{\rightarrow} (M/A'_G, \{\cdot,\cdot\}_{M/A'_G})$ which, in some situations will coincide with the symplectic reduced spaces that constitute one of the main themes of our work. We emphasize that in order to have well defined symplectic leaves in the Poisson varieties $(M/A_G, \{\cdot,\cdot\}_{M/A_G})$ and $(M/A'_G, \{\cdot,\cdot\}_{M/A'_G})$ it is very important that $A_G$ is an actual group and not just a local group of Poisson transformations and the same with the polar pseudogroup that generates $A'_G$. When the manifold $M$ is symplectic and the $G$–group action is proper it can be proved that there exists a polar family $\mathcal{A}_G^{c'}$ made only of complete vector fields (see [O02]) which shows that $M/A'_G = M/G_{\mathcal{A}_G^{c'}}$ is the quotient space by a Poisson group action and that, therefore, its symplectic leaves are well–defined. In general we say that $A'_G$ is ***completable*** whenever there exists a polar family $\mathcal{A}_G^{c'}$ made only of complete vector fields.

# 3 Optimal reduction

We start by recalling the basics of the classical symplectic or Marsden–Weinstein reduction theory. Let $(M, \omega)$ be a symplectic manifold and $G$ be a compact connected Lie group acting freely on $(M, \omega)$ by symplectomorphisms. Suppose that this action has a standard equivariant momentum map $\mathbf{J} : M \to \mathfrak{g}^*$ associated. There are two equivalent approaches to reduction that can be found in the literature:

- **Point reduction:** it is preferable for applications in dynamics. The point reduction theorem says that for any $\mu \in \mathbf{J}(M) \subset \mathfrak{g}^*$, the quotient $\mathbf{J}^{-1}(\mu)/G_\mu$ is a symplectic manifold with symplectic form $\omega_\mu$ uniquely determined by the equality

$$\pi_\mu^* \omega_\mu = i_\mu^* \omega,$$

  where $G_\mu$ is the isotropy subgroup of the element $\mu \in \mathfrak{g}^*$ with respect to the coadjoint action of $G$ on $\mathfrak{g}^*$, $i_\mu : \mathbf{J}^{-1}(\mu) \hookrightarrow M$ is the canonical injection, and $\pi_\mu : \mathbf{J}^{-1}(\mu) \to \mathbf{J}^{-1}(\mu)/G_\mu$ the projection onto the orbit space.

- **Orbit reduction:** this approach is particularly important in the treatment of quantization questions. Let $\mathcal{O}$ be the coadjoint orbit of some element $\mu \in \mathbf{J}(M)$. The subset $\mathbf{J}^{-1}(\mathcal{O})$ is a smooth submanifold of $M$ and the quotient $\mathbf{J}^{-1}(\mathcal{O})/G$ is a regular symplectic quotient manifold with the symplectic form $\omega_\mathcal{O}$ determined by the equality

$$i_\mathcal{O}^* \omega = \pi_\mathcal{O}^* \omega_\mathcal{O} + \mathbf{J}_\mathcal{O}^* \omega_\mathcal{O}^+, \tag{3.1}$$

where $i_{\mathcal{O}} : \mathbf{J}^{-1}(\mathcal{O}) \hookrightarrow M$ is the inclusion, $\pi_{\mathcal{O}} : \mathbf{J}^{-1}(\mathcal{O}) \rightarrow \mathbf{J}^{-1}(\mathcal{O})/G$ the projection, $\mathbf{J}_{\mathcal{O}} = \mathbf{J}|_{\mathbf{J}^{-1}(\mathcal{O})}$, and $\omega_{\mathcal{O}}^{+}$ the "+" orbit symplectic structure on $\mathcal{O}$ (also called Kostant–Kirillov–Souriau —KKS for short— symplectic structure). The use of the orbit reduction approach is particularly convenient when we are interested in the study of the geometry of the orbit space $M/G$ as a Poisson manifold. Indeed, the connected components of $\mathbf{J}^{-1}(\mathcal{O})/G$ constitute the symplectic leaves of $M/G$ and expression (3.1) appears as a corollary of the theory of dual pairs. Indeed, as we already said in the introduction, the diagram $M/G \overset{\pi_{A_G}}{\leftarrow} (M, \omega) \overset{\mathbf{J}}{\rightarrow} \mathbf{J}(M) \subset \mathfrak{g}^*$ forms a dual pair. It has been shown [W83, Bl01] that whenever we have two Poisson manifolds in the legs of a dual pair $(P_1, \{\cdot,\cdot\}_{P_1}) \overset{\pi_1}{\leftarrow} (M, \omega) \overset{\pi_2}{\rightarrow} (P_2, \{\cdot,\cdot\}_{P_2})$ and $\pi_1$ and $\pi_2$ have connected fibers, its symplectic leaves are in bijection. Moreover, if two symplectic leaves $\mathcal{L}_1 \subset P_1$ and $\mathcal{L}_2 \subset P_2$ are in correspondence, their symplectic structures $\omega_{\mathcal{L}_1}$ and $\omega_{\mathcal{L}_2}$ are linked by the equality

$$i_{\mathcal{K}}^{*}\omega = \pi_1|_{\mathcal{K}}^{*}\omega_{\mathcal{L}_1} + \pi_2|_{\mathcal{K}}^{*}\omega_{\mathcal{L}_2}, \tag{3.2}$$

where $\mathcal{K} \subset M$ is the leaf of the integrable distribution $\ker T\pi_1 + \ker T\pi_2$ that contains both $\pi_1^{-1}(\mathcal{L}_1)$ and $\pi_2^{-1}(\mathcal{L}_2)$. Therefore, if we assume that $\mathbf{J}$ has connected fibers, expression (3.1) appears as a corollary of (3.2), given that $\mathbf{J}(\mathcal{O})/G$ and $\mathcal{O}$ are symplectic leaves in correspondence of $M/G$ and $\mathbf{J}(M) \subset \mathfrak{g}^*$, respectively.

The use of the optimal momentum map allows the extension of these reduction procedures to far more general situations. Indeed, as we will see in the following paragraphs, the optimal approach allows the construction of symplectically reduced spaces purely within the Poisson category under hypothesis that do not necessarily imply the existence of a standard momentum map. Moreover, we will develop an orbit reduction procedure that in the context of the dual pairs reviewed in Section 2.3 reproduces the beautiful interplay between symplectic reduction and Poisson geometry that we just reviewed. We begin our study with point reduction.

## 3.1 Optimal point reduction

The study of this approach has been carried out in [O02a]. We reproduce here the main result in that paper. In the statement we will denote by $\pi_\rho : \mathcal{J}^{-1}(\rho) \rightarrow \mathcal{J}^{-1}(\rho)/G_\rho$ the canonical projection onto the orbit space of the $G_\rho$–action on $\mathcal{J}^{-1}(\rho)$ defined in Proposition 2.1.

**Theorem 3.1 (Optimal point reduction by Poisson actions)** *Let $(M, \{\cdot,\cdot\})$ be a smooth Poisson manifold and $G$ be a Lie group acting canonically and properly on $M$. Let $\mathcal{J} : M \rightarrow M/A'_G$ be the optimal momentum map associated to this action. Then, for any $\rho \in M/A'_G$ whose isotropy subgroup $G_\rho$ acts properly on $\mathcal{J}^{-1}(\rho)$, the orbit space $M_\rho := \mathcal{J}^{-1}(\rho)/G_\rho$ is a smooth symplectic regular quotient manifold with symplectic form $\omega_\rho$ defined by:*

$$\pi_\rho^*\omega_\rho(m)(X_f(m), X_h(m)) = \{f, h\}(m), \quad \text{for any } m \in \mathcal{J}^{-1}(\rho) \text{ and any } f, h \in C^\infty(M)^G. \tag{3.3}$$

*We will refer to the pair $(M_\rho, \omega_\rho)$ as the **(optimal) point reduced space** of $(M, \{\cdot,\cdot\})$ at $\rho$.*

**Remark 3.2** Let $i_{\mathcal{L}_\rho} : \mathcal{J}^{-1}(\rho) \hookrightarrow \mathcal{L}_\rho$ be the natural smooth injection of $\mathcal{J}^{-1}(\rho)$ into the symplectic leaf $(\mathcal{L}_\rho, \omega_{\mathcal{L}_\rho})$ of $(M, \{\cdot,\cdot\})$ in which it is sitting. As $\mathcal{L}_\rho$ is an initial submanifold of $M$, the injection $i_{\mathcal{L}_\rho}$ is a smooth map. The form $\omega_\rho$ can also be written in terms of the symplectic structure of the leaf $\mathcal{L}_\rho$ as

$$\pi_\rho^*\omega_\rho = i_{\mathcal{L}_\rho}^*\omega_{\mathcal{L}_\rho}. \tag{3.4}$$

The reader should be warned that this statement does NOT imply that the previous theorem could be obtained by just performing symplectic optimal reduction [OR02a] in the symplectic leaves of the Poisson manifold, basically because those leaves are not $G$–manifolds. Recall that the fact that the $G$–action is Poisson does not imply that it preserves the symplectic leaves.

In view of this remark we can obtain the standard Symplectic Stratification Theorem of Poisson manifolds as a straightforward corollary of Theorem 3.1 by taking the group $G = \{e\}$. In that case the distribution $A'_G$ coincides with the characteristic distribution of the Poisson manifold and the level sets of the optimal momentum map, and thereby the symplectic quotients $M_\rho$, are exactly the symplectic leaves. We explicitly point this out in our next statement. ♦

**Corollary 3.3 (Symplectic Stratification Theorem)** *Let $(M, \{\cdot, \cdot\})$ be a smooth Poisson manifold. Then, $M$ is the disjoint union of the maximal integral leaves of the integrable generalized distribution $E$ given by*

$$E(m) := \{X_f(m) \mid f \in C^\infty(M)\}, \quad m \in M.$$

*These leaves are symplectic initial submanifolds of $M$.*

**Remark 3.4** The only extra hypothesis in the statement of Theorem 3.1 with respect to the hypotheses used in the classical reduction theorems is the properness of the $G_\rho$–action on $\mathcal{J}^{-1}(\rho)$. This is a real hypothesis in the sense that the properness of the $G_\rho$–action is not automatically inherited from the properness of the $G$–action on $M$, as it used to be the case in the presence of a standard momentum map (see [OR02a]). For an example illustrating that this is really the case the reader may want to check with [O02a]. ♦

The interest of reduction in Poisson dynamics is justified by the following result whose proof is a simple diagram chasing exercise.

**Theorem 3.5 (Optimal point reduction of $G$–equivariant Poisson dynamics)** *Let $(M, \{\cdot, \cdot\})$ be a smooth Poisson manifold and $G$ be a Lie group acting canonically and properly on $M$. Let $\mathcal{J} : M \to M/A'_G$ be the optimal momentum map associated and $\rho \in M/A'_G$ be such that $G_\rho$ acts properly on $\mathcal{J}^{-1}(\rho)$. Let $h \in C^\infty(M)^G$ be a $G$–invariant function on $M$ and $X_h$ be the associated $G$–equivariant Hamiltonian vector field on $M$. Then,*

**(i)** *The flow $F_t$ of $X_h$ leaves $\mathcal{J}^{-1}(\rho)$ invariant, commutes with the $G$–action, and therefore induces a flow $F^\rho_t$ on $M_\rho$ uniquely determined by the relation $\pi_\rho \circ F_t \circ i_\rho = F^\rho_t \circ \pi_\rho$, where $i_\rho : \mathcal{J}^{-1}(\rho) \hookrightarrow M$ is the inclusion.*

**(ii)** *The flow $F^\rho_t$ in $(M_\rho, \omega_\rho)$ is Hamiltonian with the Hamiltonian function $h_\rho \in C^\infty(M_\rho)$ given by the equality $h_\rho \circ \pi_\rho = h \circ i_\rho$.*

**(iii)** *Let $k \in C^\infty(M)^G$ be another $G$–invariant function on $M$ and $\{\cdot, \cdot\}_\rho$ be the Poisson bracket associated to the symplectic form $\omega_\rho$ on $M_\rho$. Then, $\{h, k\}_\rho = \{h_\rho, k_\rho\}_\rho$.*

## 3.2 The symplectic case and Sjamaar's Principle

In the next few paragraphs we will see that when $M$ is a symplectic manifold with form $\omega$, the optimal point reduction by the $G$–action on $M$ produces the same results as the reduction of the isotropy type submanifolds by the relevant remaining group actions on them. In the globally Hamiltonian context, that is, in the presence of a $G$–equivariant momentum map, this idea is usually referred to as ***Sjamaar's principle*** [S90, SL91].

Let $\mathcal{J} : M \to M/A'_G$ be the optimal momentum map corresponding to the proper $G$–action on $(M, \omega)$. Fix $\rho \in M/A'_G$ a momentum value of $\mathcal{J}$ and let $H \subset G$ be the unique $G$–isotropy subgroup such that $\mathcal{J}^{-1}(\rho) \subset M_H$ and $G_\rho \subset H$. Recall that the normalizer $N(H)$ of $H$ in $G$ acts naturally. This action induces a free action of the quotient group $L := N(H)/H$ on $M_H$. Let $M_H^\rho$ be the unique connected component of $M_H$ that contains $\mathcal{J}^{-1}(\rho)$ and $L^\rho$ be the closed subgroup of $L$ that leaves it invariant. Obviously, $L^\rho$ can be written as $L^\rho = N(H)^\rho/H$ for some closed subgroup $N(H)^\rho$ of $N(H)$.

The subset $M_H^\rho$ is a symplectic embedded submanifold of $M$ where the group $L^\rho$ acts freely and canonically. We will denote by $\mathcal{J}_{L^\rho} : M_H^\rho \to M_H^\rho/A'_{L^\rho}$ the associated optimal momentum map. The following proposition explains the interest of this construction. We omit the proof since it is a straightforward consequence of the existence of local $G$–invariant extensions to $M$ for the $L^\rho$–invariant smooth functions defined in $M_H^\rho$ that has been proved in Lemma 4.4 of [OR02a].

**Proposition 3.6 (Optimal Sjamaar's Principle)** *Let $G$ be a Lie group that acts properly and canonically on the symplectic manifold $(M, \omega)$, with associated optimal momentum map $\mathcal{J} : M \to M/A'_G$. Let $\rho \in M/A'_G$ and $H \subset G$ be the unique $G$–isotropy subgroup such that $\mathcal{J}^{-1}(\rho) \subset M_H$ and $G_\rho \subset H$. Then, with the notation introduced in the previous paragraphs we have that:*

**(i)** *Let $i_H^\rho : M_H^\rho \hookrightarrow M$ be the inclusion. For any $z \in M_H^\rho$ we have that $T_z i_H^\rho \cdot A'_{L^\rho}(z) = A'_G(z)$.*

**(ii)** *Let $z \in \mathcal{J}^{-1}(\rho)$ be such that $\mathcal{J}_{L^\rho}(z) = \sigma \in M_H^\rho/A'_{L^\rho}$. Then, $\mathcal{J}^{-1}(\rho) = \mathcal{J}_{L^\rho}^{-1}(\sigma)$.*

**(iii)** $L_\sigma^\rho = G_\rho/H$.

**(iv)** $(M_H^\rho)_\sigma = \mathcal{J}_{L^\rho}^{-1}(\sigma)/L_\sigma^\rho = \mathcal{J}^{-1}(\rho)/(G_\rho/H) = \mathcal{J}^{-1}(\rho)/G_\rho = M_\rho$. *Moreover, if $G_\rho$ acts properly on $\mathcal{J}^{-1}(\rho)$ this equality is true when we consider $M_\rho$ and $(M_H^\rho)_\sigma$ as symplectic spaces, that is,*

$$(M_\rho, \omega_\rho) = ((M_H^\rho)_\sigma, (\omega|_{M_H^\rho})_\sigma).$$

**Definition 3.7** *Suppose that we are under the hypotheses of the previous proposition. We will refer to the symplectic reduced space $((M_H^\rho)_\sigma, (\omega|_{M_H^\rho})_\sigma)$ as the **regularization** of the point reduced space $(M_\rho, \omega_\rho)$.*

## 3.3 The space for optimal orbit reduction

The main difference between the point and orbit reduced spaces is in the invariance properties of the submanifolds out of which they are constructed. More specifically, if we mimic the standard orbit reduction procedure using the momentum map, the optimal orbit reduced space that we should study is $G \cdot \mathcal{J}^{-1}(\rho)/G = \mathcal{J}^{-1}(\mathcal{O}_\rho)/G$, where $\mathcal{O}_\rho := G \cdot \rho \subset M/A'_G$. Hence, the first question that we have to tackle is: is there a canonical smooth structure for $\mathcal{J}^{-1}(\mathcal{O}_\rho)$ and $\mathcal{J}^{-1}(\mathcal{O}_\rho)/G$ that we can use to carry out the orbit reduction scheme in this framework?

We will first show that there is an affirmative answer for the smooth structure of $\mathcal{J}^{-1}(\mathcal{O}_\rho)$. The main idea that we will prove in the following paragraphs is that *$\mathcal{J}^{-1}(\mathcal{O}_\rho)$ can be naturally endowed with* *the unique smooth structure that makes it into an initial submanifold of* $M$. We start with the following proposition whose proof can be found in the appendix.

**Proposition 3.8** *Let* $(M, \{\cdot, \cdot\})$ *be a smooth Poisson manifold and* $G$ *be a Lie group acting canonically and properly on* $M$. *Let* $\mathcal{J} : M \rightarrow M/A'_G$ *be the optimal momentum map associated to this action. Then,*

**(i)** *The generalized distribution* $D$ *on* $M$ *defined by* $D(m) := \mathfrak{g} \cdot m + A'_G(m)$, *for all* $m \in M$, *is integrable.*

**(ii)** *Let* $m \in M$ *be such that* $\mathcal{J}(m) = \rho$, *then* $G^0 \cdot \mathcal{J}^{-1}(\rho)$ *is the maximal integral submanifold of* $D$ *going through the point* $m$. *The symbol* $G^0$ *denotes the connected component of* $G$ *containing the identity.*

As we already said, a general fact about integrable generalized distributions [Daz85] states that the smooth structure on a subset of $M$ that makes it into a maximal integral manifold of a given distribution coincides with the *unique* smooth structure that makes it into an initial submanifold of $M$. Therefore, the previous proposition shows that the sets $G^0 \cdot \mathcal{J}^{-1}(\rho)$ are initial submanifolds of $M$.

For the proof of the following proposition see the appendix.

**Proposition 3.9** *Suppose that we have the same setup as in Proposition 3.8. If either* $G_\rho$ *is closed in* $G$ *or, more generally,* $G_\rho$ *acts properly on* $\mathcal{J}^{-1}(\rho)$, *then:*

**(i)** *The* $G_\rho$ *action on the product* $G \times \mathcal{J}^{-1}(\rho)$ *defined by* $h \cdot (g, z) := (gh, h^{-1} \cdot z)$ *is free and proper and therefore, the corresponding orbit space* $G \times \mathcal{J}^{-1}(\rho)/G_\rho =: G \times_{G_\rho} \mathcal{J}^{-1}(\rho)$ *is a smooth regular quotient manifold. We will denote by* $\pi_{G_\rho} : G \times \mathcal{J}^{-1}(\rho) \rightarrow G \times_{G_\rho} \mathcal{J}^{-1}(\rho)$ *the canonical surjective submersion.*

**(ii)** *The mapping* $i : G \times_{G_\rho} \mathcal{J}^{-1}(\rho) \rightarrow M$ *defined by* $i([g, z]) := g \cdot z$ *is an injective immersion onto* $\mathcal{J}^{-1}(\mathcal{O}_\rho)$ *such that, for any* $[g, z] \in G \times_{G_\rho} \mathcal{J}^{-1}(\rho)$, $T_{[g,z]} i \cdot T_{[g,z]}(G \times_{G_\rho} \mathcal{J}^{-1}(\rho)) = D(g \cdot z)$. *On other words* $i(G \times_{G_\rho} \mathcal{J}^{-1}(\rho)) = \mathcal{J}^{-1}(\mathcal{O}_\rho)$ *is an integral submanifold of* $D$.

By using the previous propositions we will now show that, in the presence of the standard hypotheses for reduction, $\mathcal{J}^{-1}(\mathcal{O}_\rho)$ is an initial submanifold of $M$ whose connected components are the also initial submanifolds $gG^0 \cdot \mathcal{J}^{-1}(\rho)$, $g \in G$. We start with the following definition:

**Definition 3.10** *Let* $(M, \{\cdot, \cdot\})$ *be a smooth Poisson manifold and* $G$ *be a Lie group acting canonically and properly on* $M$. *Let* $\mathcal{J} : M \rightarrow M/A'_G$ *be the optimal momentum map associated to this action and* $\rho \in M/A'_G$. *Suppose that* $G_\rho$ *acts properly on* $\mathcal{J}^{-1}(\rho)$. *In these circumstances, by Proposition 3.9, the twist product* $G \times_{G_\rho} \mathcal{J}^{-1}(\rho)$ *has a canonical smooth structure. Consider in the set* $\mathcal{J}^{-1}(\mathcal{O}_\rho)$ *the smooth structure that makes the bijection* $G \times_{G_\rho} \mathcal{J}^{-1}(\rho) \rightarrow \mathcal{J}^{-1}(\mathcal{O}_\rho)$ *given by* $(g, z) \rightarrow g \cdot z$ *into a diffeomorphism. We will refer to this structure as the* ***initial smooth structure*** *of* $\mathcal{J}^{-1}(\mathcal{O}_\rho)$.

The following theorem justifies the choice of terminology in the previous definition and why we will be able to refer to the smooth structure there introduced as *THE* initial smooth structure of $\mathcal{J}^{-1}(\mathcal{O}_\rho)$.

**Theorem 3.11** *Suppose that we are in the same setup as in Definition 3.10. Then, the set* $\mathcal{J}^{-1}(\mathcal{O}_\rho)$ *endowed with the initial smooth structure is an actual initial submanifold of* $M$ *that can be decomposed as a disjoint union of connected components as*

$$\mathcal{J}^{-1}(\mathcal{O}_\rho) = \dot{\bigcup_{[g] \in G/(G^0 G_\rho)}} gG^0 \cdot \mathcal{J}^{-1}(\rho). \tag{3.5}$$

*Each connected component of $\mathcal{J}^{-1}(\mathcal{O}_\rho)$ is a maximal integral submanifold of the distribution $D$ defined in Proposition 3.8. If, additionally, the subgroup $G_\rho$ is closed in $G$, the topology on $\mathcal{J}^{-1}(\mathcal{O}_\rho)$ induced by its initial smooth structure coincides with the initial topology induced by the map $\mathcal{J}_{\mathcal{J}^{-1}(\mathcal{O}_\rho)} : \mathcal{J}^{-1}(\mathcal{O}_\rho) \rightarrow \mathcal{O}_\rho$ given by $z \longmapsto \mathcal{J}(z)$, where the orbit $\mathcal{O}_\rho$ is endowed with the smooth structure coming from the homogeneous manifold $G/G_\rho$. Finally, notice that (3.5) implies that $\mathcal{J}^{-1}(\mathcal{O}_\rho)$ has as many connected components as the cardinality of the homogeneous manifold $G/(G^0G_\rho)$.*

**Proof.** First of all notice that the sets $gG^0 \cdot \mathcal{J}^{-1}(\rho)$ are clearly maximal integral submanifolds of $D$ by part **(ii)** in Proposition 3.8. As a corollary of this, they are the connected components of $\mathcal{J}^{-1}(\mathcal{O}_\rho)$ endowed with the smooth structure in Definition 3.11. Indeed, let $S$ be the connected component of $\mathcal{J}^{-1}(\mathcal{O}_\rho)$ that contains $gG^0\cdot\mathcal{J}^{-1}(\rho)$, that is, $gG^0\cdot\mathcal{J}^{-1}(\rho) \subset S \subset \mathcal{J}^{-1}(\mathcal{O}_\rho)$. As $\mathcal{J}^{-1}(\mathcal{O}_\rho)$ is a manifold, it is locally connected, and therefore its connected components are open and closed. In particular, since $S$ is an open connected subset of $\mathcal{J}^{-1}(\mathcal{O}_\rho)$, part **(ii)** in Proposition 3.9 shows that $S$ is a connected integral submanifold of $D$. By the maximality of $gG^0\cdot\mathcal{J}^{-1}(\rho)$ as an integral submanifold of $D$, $gG^0\cdot\mathcal{J}^{-1}(\rho) = S$, necessarily. The set $gG^0 \cdot \mathcal{J}^{-1}(\rho)$ is therefore a connected component of $\mathcal{J}^{-1}(\mathcal{O}_\rho)$. As it is a leaf of a smooth integrable distribution on $M$, it is also an initial submanifold of $M$ [Daz85] of dimension $d = \dim \mathcal{J}^{-1}(\mathcal{O}_\rho) = \dim G + \dim \mathcal{J}^{-1}(\rho) - \dim G_\rho$.

We now show that $\mathcal{J}^{-1}(\mathcal{O}_\rho)$ with the smooth structure in Definition 3.11 is an initial submanifold of $M$. First of all part **(ii)** in Proposition 3.9 shows that $\mathcal{J}^{-1}(\mathcal{O}_\rho)$ is an injectively immersed submanifold of $M$. The initial character can be obtained as a consequence of the fact that its connected components are initial together with the following elementary lemma:

**Lemma 3.12** *Let $N$ be an injectively immersed submanifold of the smooth manifold $M$. Suppose that $N$ can be written as the disjoint union of a family $\{S_\alpha\}_{\alpha\in I}$ of open subsets of $N$ such that each $S_\alpha$ is an initial submanifold of $M$. Then, $N$ is initial.*

**Proof of the lemma.** Let $i_N : N \hookrightarrow M$ and $i_\alpha : S_\alpha \hookrightarrow N$ be the injections. Let $Z$ be an arbitrary smooth manifold and $f : Z \rightarrow M$ be a smooth map such that $f(Z) \subset N$. As the sets $S_\alpha$ are open and partition $N$, the manifold $Z$ can be written as a disjoint union of open sets $Z_\alpha := f^{-1}(S_\alpha)$, that is

$$Z = \dot{\bigcup_{\alpha\in I}} f^{-1}(S_\alpha).$$

Given that for each index $\alpha$ the map $f_\alpha : Z_\alpha \rightarrow M$ obtained by restriction of $f$ to $Z_\alpha$ is smooth, the corresponding map $\bar{f}_\alpha : Z_\alpha \rightarrow S_\alpha$ defined by the identity $i_\alpha \circ \bar{f}_\alpha = f_\alpha$ is also smooth by the initial character of $S_\alpha$. Let $\bar{f} : Z \rightarrow N$ be the map obtained by union of the mappings $\bar{f}_\alpha$. This map is smooth and satisfies that $i_N \circ \bar{f} = f$ which proves that $N$ is initial. ▼

We now prove Expression (3.5). First of all notice that as $G^0$ is normal in $G$, the set $G^0G_\rho$ is a (in principle not closed) subgroup of $G$. We obviously have that

$$\mathcal{J}^{-1}(\mathcal{O}_\rho) = \bigcup_{g\in G} gG^0\mathcal{J}^{-1}(\rho). \tag{3.6}$$

Moreover, if $g$ and $g' \in G$ are such that $[g] = [g'] \in G/(G^0G_\rho)$ then we can write that $g' = ghk$ with $h \in G^0$ and $k \in G_\rho$. Consequently, $g'G^0\mathcal{J}^{-1}(\rho) = ghkG^0\mathcal{J}^{-1}(\rho) = gh(G^0k)\mathcal{J}^{-1}(\rho) = g(hG^0)(k\mathcal{J}^{-1}(\rho)) = gG^0\mathcal{J}^{-1}(\rho)$, which implies that (3.6) can be refined to

$$\mathcal{J}^{-1}(\mathcal{O}_\rho) = \bigcup_{[g]\in G/(G^0G_\rho)} gG^0\mathcal{J}^{-1}(\rho). \tag{3.7}$$

It only remains to be shown that this union is disjoint: let $gh \cdot z = lh' \cdot z'$ with $h, h' \in G^0$ and $z, z' \in \mathcal{J}^{-1}(\rho)$. If we apply $\mathcal{J}$ to both sides of this equality we obtain that $gh \cdot \rho = lh' \cdot \rho$. Hence, $(h')^{-1} l^{-1} gh \in G_\rho$ and $l^{-1} g \in h' G_\rho h^{-1} \subset G^0 G_\rho$. This implies that $[l] = [g] \in G/(G^0 G_\rho)$ and $g G^0 \mathcal{J}^{-1}(\rho) = l G^0 \mathcal{J}^{-1}(\rho)$, as required.

We finally show that when $G_\rho$ is closed in $G$, the topology on $\mathcal{J}^{-1}(\mathcal{O}_\rho)$ induced by its initial smooth structure coincides with the initial topology induced by the map $\mathcal{J}_{\mathcal{J}^{-1}(\mathcal{O}_\rho)} : \mathcal{J}^{-1}(\mathcal{O}_\rho) \to \mathcal{O}_\rho$ on $\mathcal{J}^{-1}(\mathcal{O}_\rho)$. Recall first that this topology is characterized by the fact that for any topological space $Z$ and any map $\phi : Z \to \mathcal{J}^{-1}(\mathcal{O}_\rho)$ we have that $\phi : Z \to \mathcal{J}^{-1}(\mathcal{O}_\rho)$ is continuous if and only if $\mathcal{J}_{\mathcal{J}^{-1}(\mathcal{O}_\rho)} \circ \phi$ is continuous. Moreover, as the family $\{\mathcal{J}^{-1}_{\mathcal{J}^{-1}(\mathcal{O}_\rho)}(U) \mid U \text{ open subset of } \mathcal{O}_\rho\}$ is a subbase of this topology, the initial topology on $\mathcal{J}^{-1}(\mathcal{O}_\rho)$ induced by the map $\mathcal{J}_{\mathcal{J}^{-1}(\mathcal{O}_\rho)}$ is first countable. We prove that this topology coincides with the topology induced by the initial smooth structure on $\mathcal{J}^{-1}(\mathcal{O}_\rho)$ by showing that the map $f : G \times_{G_\rho} \mathcal{J}^{-1}(\rho) \to \mathcal{J}^{-1}(\mathcal{O}_\rho)$, $f([g, z]) := g \cdot z$ is a homeomorphism when we consider $\mathcal{J}^{-1}(\mathcal{O}_\rho)$ as a topological space with the initial topology induced by $\mathcal{J}_{\mathcal{J}^{-1}(\mathcal{O}_\rho)}$. Indeed, $f$ is continuous if and only if the map $G \times_{G_\rho} \mathcal{J}^{-1}(\rho) \to \mathcal{O}_\rho$ given by $[g, z] \mapsto g \cdot \rho$ is continuous, which in turn is equivalent to the continuity of the map $G \times \mathcal{J}^{-1}(\rho) \to G/G_\rho$ defined by $(g, z) \longmapsto g G_\rho$, which is true. We now show that the inverse $f^{-1} : \mathcal{J}^{-1}(\mathcal{O}_\rho) \to G \times_{G_\rho} \mathcal{J}^{-1}(\rho)$ of $f$ given by $g \cdot z \mapsto [g, z]$ is continuous. Since the initial topology on $\mathcal{J}^{-1}(\mathcal{O}_\rho)$ induced by $\mathcal{J}_{\mathcal{J}^{-1}(\mathcal{O}_\rho)}$ is first countable it suffices to show that for any convergent sequence $\{z_n\} \subset \mathcal{J}^{-1}(\mathcal{O}_\rho) \to z \in \mathcal{J}^{-1}(\mathcal{O}_\rho)$, we have that $\lim_{n \to \infty} f^{-1}(z_n) = f^{-1}(\lim_{n \to \infty} z_n) = f^{-1}(z)$. Indeed, as $\mathcal{J}_{\mathcal{J}^{-1}(\mathcal{O}_\rho)}$ is continuous, the sequence $\{\mathcal{J}(z_n) = g_n \cdot \rho\} \subset \mathcal{O}_\rho$ converges in $\mathcal{O}_\rho$ to $\mathcal{J}(z) = g \cdot \rho$, for some $g \in G$. Let $j : \mathcal{O}_\rho \to G/G_\rho$ be the standard diffeomorphism and $\sigma : U_{gG_\rho} \subset G/G_\rho \to G$ be a local smooth section of the submersion $G \to G/G_\rho$ in a neighborhood $U_{gG_\rho}$ of $gG_\rho \in G/G_\rho$. Let $V := \mathcal{J}^{-1}_{\mathcal{J}^{-1}(\mathcal{O}_\rho)}(j^{-1}(U_{gG_\rho}))$. $V$ is an open neighborhood of $z$ in $\mathcal{J}^{-1}(\mathcal{O}_\rho)$ because $j \circ \mathcal{J}_{\mathcal{J}^{-1}(\mathcal{O}_\rho)}(z) = j(g \cdot \rho) = gG_\rho \in U_{gG_\rho}$. We now notice that for any $m \in V$ we can write that

$$f^{-1}(m) = [\sigma \circ j \circ \mathcal{J}_{\mathcal{J}^{-1}(\mathcal{O}_\rho)}(m), (\sigma \circ j \circ \mathcal{J}_{\mathcal{J}^{-1}(\mathcal{O}_\rho)}(m))^{-1} \cdot m].$$

Consequently, as

$$\begin{aligned} \lim_{n \to \infty} f^{-1}(z_n) &= \lim_{n \to \infty} [\sigma \circ j \circ \mathcal{J}_{\mathcal{J}^{-1}(\mathcal{O}_\rho)}(z_n), (\sigma \circ j \circ \mathcal{J}_{\mathcal{J}^{-1}(\mathcal{O}_\rho)}(z_n))^{-1} \cdot z_n] \\ &= [\sigma \circ j \circ \mathcal{J}_{\mathcal{J}^{-1}(\mathcal{O}_\rho)}(z), (\sigma \circ j \circ \mathcal{J}_{\mathcal{J}^{-1}(\mathcal{O}_\rho)}(z))^{-1} \cdot z] = f^{-1}(z), \end{aligned}$$

the continuity of $f^{-1}$ is guaranteed. ■

## 3.4 The symplectic orbit reduction quotient

We will know show that the quotient $\mathcal{J}^{-1}(\mathcal{O}_\rho)/G$ can be endowed with a smooth structure that makes it into a regular quotient manifold, that is, the projection $\pi_{\mathcal{O}_\rho} : \mathcal{J}^{-1}(\mathcal{O}_\rho) \to \mathcal{J}^{-1}(\mathcal{O}_\rho)/G$ is a smooth submersion. We will carry this out under the same hypotheses present in Definition 3.10, that is, $G_\rho$ acts properly on $\mathcal{J}^{-1}(\rho)$.

First of all notice that as $\mathcal{J}^{-1}(\mathcal{O}_\rho)$ is an initial $G$–invariant submanifold of $M$, the $G$–action on $\mathcal{J}^{-1}(\mathcal{O}_\rho)$ is smooth. We will prove that $\mathcal{J}^{-1}(\mathcal{O}_\rho)/G$ is a regular quotient manifold by showing that this action is actually proper and satisfies that all the isotropy subgroups are conjugate to a given one. Indeed, recall that the initial manifold structure on $\mathcal{J}^{-1}(\mathcal{O}_\rho)$ is the one that makes it $G$–equivariantly diffeomorphic to the twist product $G \times_{G_\rho} \mathcal{J}^{-1}(\rho)$ when we take in this space the $G$–action given by the expression $g \cdot [h, z] := [gh, z]$, $g \in G$, $[h, z] \in G \times_{G_\rho} \mathcal{J}^{-1}(\rho)$. Therefore, it suffices to show that this $G$–action has the desired properties. First of all this action is proper since a general property about twist products (see [OR02b]) says that the $G$–action on $G \times_{G_\rho} \mathcal{J}^{-1}(\rho)$ is proper iff the $G_\rho$–action on $\mathcal{J}^{-1}(\rho)$ is proper, which we supposed as a hypothesis. We now look at the isotropies of this action: in Proposition 2.3 we saw that all the elements in $\mathcal{J}^{-1}(\rho)$ have the same $G$–isotropy, call it $H$. As $H \subset G_\rho$, this is also their $G_\rho$–isotropy. Now, using a standard property of the isotropies of twist products [OR02b], we have that $G_{[g,z]} = g(G_\rho)_z g^{-1} = gHg^{-1}$, for any $[g,z] \in G \times_{G_\rho} \mathcal{J}^{-1}(\rho)$, as required.

The quotient manifold $\mathcal{J}^{-1}(\mathcal{O}_\rho)/G$ is naturally diffeomorphic to the symplectic point reduced space. Indeed,

$$\mathcal{J}^{-1}(\mathcal{O}_\rho)/G \simeq G \times_{G_\rho} \mathcal{J}^{-1}(\rho)/G \simeq \mathcal{J}^{-1}(\rho)/G_\rho.$$

This diffeomorphism can be explicitly implemented as follows. Let $l_\rho : \mathcal{J}^{-1}(\rho) \to \mathcal{J}^{-1}(\mathcal{O}_\rho)$ be the inclusion. As the inclusion $\mathcal{J}^{-1}(\rho) \hookrightarrow M$ is smooth and $\mathcal{J}^{-1}(\mathcal{O}_\rho)$ is initial $l_\rho$ is smooth. Also, since $l_\rho$ is $(G_\rho, G)$ equivariant it drops to a unique smooth map $L_\rho : \mathcal{J}^{-1}(\rho)/G_\rho \to \mathcal{J}^{-1}(\mathcal{O}_\rho)/G$ that makes the following diagram

$$\begin{array}{ccc} \mathcal{J}^{-1}(\rho) & \xrightarrow{\ l_\rho\ } & \mathcal{J}^{-1}(\mathcal{O}_\rho) \\ {\scriptstyle \pi_\rho}\downarrow & & \downarrow{\scriptstyle \pi_{\mathcal{O}_\rho}} \\ \mathcal{J}^{-1}(\rho)/G_\rho & \xrightarrow{\ L_\rho\ } & \mathcal{J}^{-1}(\mathcal{O}_\rho)/G. \end{array}$$

commutative. $L_\rho$ is a smooth bijection. In order to show that its inverse is also smooth we will think of $\mathcal{J}^{-1}(\mathcal{O}_\rho)$ as $G \times_{G_\rho} \mathcal{J}^{-1}(\rho)$. First of all notice that the projection $G \times \mathcal{J}^{-1}(\rho) \to \mathcal{J}^{-1}(\rho)$ is $G_\rho$–(anti)equivariant and therefore induces a smooth map $G \times_{G_\rho} \mathcal{J}^{-1}(\rho) \to \mathcal{J}^{-1}(\rho)/G_\rho$ given by $[g,z] \mapsto [z]$, $[g,z] \in G \times_{G_\rho} \mathcal{J}^{-1}(\rho)$. This map is $G$–invariant and therefore drops to another smooth mapping $G \times_{G_\rho} \mathcal{J}^{-1}(\rho)/G \to \mathcal{J}^{-1}(\rho)/G_\rho$ that coincides with $L_\rho^{-1}$, the inverse of $L_\rho$, which is consequently a diffeomorphism.

The orbit reduced space $\mathcal{J}^{-1}(\mathcal{O}_\rho)/G$ can be therefore trivially endowed with a symplectic structure $\omega_{\mathcal{O}_\rho}$ by defining $\omega_{\mathcal{O}_\rho} := (L_\rho^{-1})^* \omega_\rho$. We put together all the facts that we just proved in the following theorem–definition:

**Theorem 3.13 (Optimal orbit reduction by Poisson actions)** *Let $(M, \{\cdot,\cdot\})$ be a smooth Poisson manifold and $G$ be a Lie group acting canonically and properly on $M$. Let $\mathcal{J} : M \to M/A'_G$ be the optimal momentum map associated to this action and $\rho \in M/A'_G$. Suppose that $G_\rho$ acts properly on $\mathcal{J}^{-1}(\rho)$. If we denote $\mathcal{O}_\rho := G \cdot \rho$, then:*

**(i)** *There is a unique smooth structure on $\mathcal{J}^{-1}(\mathcal{O}_\rho)$ that makes it into an initial submanifold of $M$.*

**(ii)** *The $G$–action on $\mathcal{J}^{-1}(\mathcal{O}_\rho)$ by restriction of the $G$–action on $M$ is smooth and proper and all its isotropy subgroups are conjugate to a given compact isotropy subgroup of the $G$–action on $M$.*

**(iii)** *The quotient $M_{\mathcal{O}_\rho} := \mathcal{J}^{-1}(\mathcal{O}_\rho)/G$ admits a unique smooth structure that makes the projection $\pi_{\mathcal{O}_\rho} : \mathcal{J}^{-1}(\mathcal{O}_\rho) \to \mathcal{J}^{-1}(\mathcal{O}_\rho)/G$ a surjective submersion.*

**(iv)** *The quotient $M_{\mathcal{O}_\rho} := \mathcal{J}^{-1}(\mathcal{O}_\rho)/G$ admits a unique symplectic structure $\omega_{\mathcal{O}_\rho}$ that makes it symplectomorphic to the point reduced space $M_\rho$. We will refer to the pair $(M_{\mathcal{O}_\rho}, \omega_{\mathcal{O}_\rho})$ as the* ***(optimal) orbit reduced space*** *of $(M, \{\cdot,\cdot\})$ at $\mathcal{O}_\rho$.*

In this setup we can easily formulate an analog of Theorem 3.5.

**Theorem 3.14 (Optimal orbit reduction of $G$–equivariant Poisson dynamics)** *Let $(M, \{\cdot, \cdot\})$ be a smooth Poisson manifold and $G$ be a Lie group acting canonically and properly on $M$. Let $\mathcal{J} : M \to M/A'_G$ be the optimal momentum map associated and $\rho \in M/A'_G$ be such that $G_\rho$ acts properly on $\mathcal{J}^{-1}(\rho)$. Let $h \in C^\infty(M)^G$ be a $G$–invariant function on $M$ and $X_h$ be the associated $G$–equivariant Hamiltonian vector field on $M$. Then,*

**(i)** *The flow $F_t$ of $X_h$ leaves $\mathcal{J}^{-1}(\mathcal{O}_\rho)$ invariant, commutes with the $G$–action, and therefore induces a flow $F_t^{\mathcal{O}_\rho}$ on $M_{\mathcal{O}_\rho}$ uniquely determined by the relation $\pi_{\mathcal{O}_\rho} \circ F_t \circ i_{\mathcal{O}_\rho} = F_t^{\mathcal{O}_\rho} \circ \pi_{\mathcal{O}_\rho}$, where $i_{\mathcal{O}_\rho} : \mathcal{J}^{-1}(\mathcal{O}_\rho) \hookrightarrow M$ is the inclusion.*

**(ii)** *The flow $F_t^{\mathcal{O}_\rho}$ in $(M_{\mathcal{O}_\rho}, \omega_{\mathcal{O}_\rho})$ is Hamiltonian with the Hamiltonian function $h_{\mathcal{O}_\rho} \in C^\infty(M_{\mathcal{O}_\rho})$ given by the equality $h_{\mathcal{O}_\rho} \circ \pi_{\mathcal{O}_\rho} = h \circ i_{\mathcal{O}_\rho}$.*

**(iii)** *Let $k \in C^\infty(M)^G$ be another $G$–invariant function on $M$ and $\{\cdot, \cdot\}_{\mathcal{O}_\rho}$ be the Poisson bracket associated to the symplectic form $\omega_{\mathcal{O}_\rho}$ on $M_{\mathcal{O}_\rho}$. Then, $\{h, k\}_{\mathcal{O}_\rho} = \{h_{\mathcal{O}_\rho}, k_{\mathcal{O}_\rho}\}_{\mathcal{O}_\rho}$.*

We conclude this section with a brief description of the orbit version of the regularized reduced spaces introduced in Definition 3.7 for the symplectic case. If we follow the prescription introduced in Section 3.3 using the $L^\rho$–action on $M_H^\rho$ we are first supposed to study the set $\mathcal{J}_{L^\rho}^{-1}(L^\rho \cdot \sigma)$. The initial smooth structure on this set induced by the twist product $L^\rho \times_{L^\rho_\sigma} \mathcal{J}_{L^\rho}^{-1}(\sigma)$ makes it into an initial submanifold of $M_H^\rho$. Moreover, if we use the statements in Proposition 3.7 it is easy to see that $\mathcal{J}_{L^\rho}^{-1}(L^\rho \cdot \sigma) = L^\rho \cdot \mathcal{J}_{L^\rho}^{-1}(\sigma) = N(H)^\rho \cdot \mathcal{J}^{-1}(\rho) = \mathcal{J}^{-1}(\mathcal{N}_\rho)$, with $\mathcal{N}_\rho := N(H)^\rho \cdot \rho \subset M/A'_G$.

The set $\mathcal{J}_{L^\rho}^{-1}(L^\rho \cdot \sigma) = \mathcal{J}^{-1}(\mathcal{N}_\rho)$ is an embedded submanifold of $\mathcal{J}^{-1}(\mathcal{O}_\rho)$ (since $\mathcal{J}^{-1}(\mathcal{N}_\rho) \simeq N(H)^\rho \times_{G_\rho} \mathcal{J}^{-1}(\rho)$ is embedded in $G \times_{G_\rho} \mathcal{J}^{-1}(\rho) \simeq \mathcal{J}^{-1}(\mathcal{O}_\rho)$). Moreover, a simple diagram chasing shows that the symplectic quotient $(\mathcal{J}_{L^\rho}^{-1}(L^\rho \cdot \sigma)/L^\rho, (\omega|_{M_H^\rho})_{L^\rho \cdot \sigma})$ is naturally symplectomorphic to the orbit reduced space $(\mathcal{J}^{-1}(\mathcal{O}_\rho)/G, \omega_{\mathcal{O}_\rho})$. We will say that $(\mathcal{J}_{L^\rho}^{-1}(L^\rho \cdot \sigma)/L^\rho, (\omega|_{M_H^\rho})_{L^\rho \cdot \sigma})$ is ***an orbit regularization*** of $(\mathcal{J}^{-1}(\mathcal{O}_\rho)/G, \omega_{\mathcal{O}_\rho})$.

We finally show that

$$\mathcal{J}^{-1}(\mathcal{O}_\rho) = \dot{\bigcup}_{[g] \in G/N(H)^\rho} \mathcal{J}^{-1}(\mathcal{N}_{g \cdot \rho}). \tag{3.8}$$

The equality is a straightforward consequence of the fact that for any $g \in G$,

$$M_{gHg^{-1}}^{g\rho} = \Phi_g(M_H^\rho), \quad N(gHg^{-1})^{g\rho} = gN(H)^\rho g^{-1}, \quad \text{and} \quad \mathcal{J}^{-1}(\mathcal{N}_{g \cdot \rho}) = gN(H)^\rho \mathcal{J}^{-1}(\rho). \tag{3.9}$$

The last relation implies that if $g, g' \in G$ are such that $[g] = [g'] \in G/N(H)^\rho$, then $\mathcal{J}^{-1}(\mathcal{N}_{g \cdot \rho}) = \mathcal{J}^{-1}(\mathcal{N}_{g' \cdot \rho})$. We now show that the union in (3.8) is indeed disjoint: let $gn \cdot z \in \mathcal{J}^{-1}(\mathcal{N}_{g \cdot \rho})$ and $g'n' \cdot z' \in \mathcal{J}^{-1}(\mathcal{N}_{g' \cdot \rho})$ be such that $gn \cdot z = g'n' \cdot z'$, with $g, g' \in G$, $n, n' \in N(H)^\rho$, and $z, z' \in \mathcal{J}^{-1}(\rho)$. Since $gn \cdot z = g'n' \cdot z'$, we necessarily have that $G_{gn \cdot z} = G_{g'n' \cdot z'}$ which implies that $gHg^{-1} = g'H(g')^{-1}$, and hence $g^{-1}g' \in N(H)$. We now recall that $M_H^\rho$ is the accessible set going through $z$ or $z'$ of the integrable generalized distribution $B'_G$ defined by

$$B'_G := \operatorname{span}\{X \in \mathfrak{X}(U)^G \mid U \text{ open } G\text{–invariant set in } M\},$$

where the symbol $\mathfrak{X}(U)^G$ denotes the set of $G$–equivariant vector fields defined on $U$. Let $\mathcal{B}'_G$ be the pseudogroup of transformations of $M$ consisting of the $G$–equivariant flows of the vector fields that span $B'_G$. Now, as the points $n \cdot z, n' \cdot z' \in M_H^\rho$, there exists $\mathcal{F}_T \in \mathcal{B}'_G$ such that $n' \cdot z' = \mathcal{F}_T(n \cdot z)$, hence $(g')^{-1}gn\cdot z=\mathcal{F}_T(n\cdot z)$. Moreover, as any element in $M_H^\rho$ can be written as $\mathcal{G}_T(n\cdot z)$ with $\mathcal{G}_T\in\mathcal{B}'_G$, we have that

$$(g')^{-1}g\cdot\mathcal{G}_T(n\cdot z)=\mathcal{G}_T((g')^{-1}gn\cdot z)=\mathcal{G}_T(\mathcal{F}_T(n\cdot z))\in M_H^\rho,$$

which implies that $(g')^{-1}g\in N(H)^\rho$ and therefore $[g]=[g']\in G/N(H)^\rho$, as required.

### 3.5 The polar reduced spaces

As we already recalled in the introduction to this section, the standard theory of orbit reduction provides a characterization of the symplectic form of the orbit reduced spaces in terms of the symplectic structures of the corresponding coadjoint orbits that, from the dual pairs point of view, play the role of the symplectic leaves of the Poisson manifold in duality, namely $\mathbf{J}(M)\subset\mathfrak{g}^*$.

We will now show that when the group of Poisson transformations $A_G$ is von Neumann (actually we just need weakly von Neumann), that is, when the diagram $(M/G,\{\cdot,\cdot\}_{M/A_G})\stackrel{\pi_{A_G}}{\leftarrow}(M,\{\cdot,\cdot\})\stackrel{\mathcal{J}}{\to}(M/A'_G,\{\cdot,\cdot\}_{M/A'_G})$ is a dual pair in the sense of Definition 2.7, the classical picture can be reproduced in this context. More specifically, in this section we will show that:

- The symplectic leaves of $(M/A'_G,\{\cdot,\cdot\}_{M/A'_G})$ admit a smooth presymplectic structure that generalizes the Kostant–Kirillov–Souriau symplectic structure in the coadjoint orbits of the dual of a Lie algebra in the sense that they are homogeneous presymplectic manifolds. We will refer to these "generalized coadjoint orbits" as ***polar reduced spaces***.
- The presymplectic structure of the polar reduced spaces is related to the symplectic form of the orbit reduced spaces introduced in the previous section via an equality that holds strong resemblance with the classical expression (3.1). Also, it is possible to provide a very explicit characterization of the situations in which the polar reduced spaces are actually symplectic.
- When the manifold $M$ is symplectic, the polar reduced space decomposes as a union of embedded symplectic submanifolds that correspond to the polar reduced spaces of the regularizations of the orbit reduced space. Each of these symplectic manifolds is a homogeneous manifold and we will refer to them as the ***regularized polar reduced subspaces***.

We start with a proposition that spells out the smooth structure of the polar reduced spaces. ***In this section*** we use a stronger hypothesis on $G_\rho$ with respect to the one we used in the previous section, namely, ***we will assume that*** $G_\rho$ ***is closed in*** $G$ which, as we point out in the proof of Proposition 3.9, ***implies that the*** $G_\rho$ ***action on*** $\mathcal{J}^{-1}(\rho)$ ***is proper***.

**Proposition 3.15** *Let $(M,\{\cdot,\cdot\})$ be a smooth Poisson manifold and $G$ be a Lie group acting canonically and properly on $M$. Let $\mathcal{J}:M\to M/A'_G$ be the optimal momentum map associated to this action and $\rho\in M/A'_G$. Suppose that $G_\rho$ is closed in $G$. Then, the polar distribution $A'_G$ restricts to a smooth integrable regular distribution on $\mathcal{J}^{-1}(\mathcal{O}_\rho)$, that we will also denote by $A'_G$. The leaf space $M'_{\mathcal{O}_\rho}:=\mathcal{J}^{-1}(\mathcal{O}_\rho)/A'_G$ admits a unique smooth structure that makes it into a regular quotient manifold and diffeomorphic to the homogeneous manifold $G/G_\rho$. With this smooth structure the projection $\mathcal{J}_{\mathcal{O}_\rho}:\mathcal{J}^{-1}(\mathcal{O}_\rho)\to\mathcal{J}^{-1}(\mathcal{O}_\rho)/A'_G$ is a smooth surjective submersion. We will refer to $M'_{\mathcal{O}_\rho}$ as the* ***polar reduced space****.*

**Proof.** Let $m\in\mathcal{J}^{-1}(\mathcal{O}_\rho)$. By Proposition 3.9 we have that $T_m\mathcal{J}^{-1}(\mathcal{O}_\rho)=D(m)=\mathfrak{g}\cdot m+A'_G(m)$, which implies that the restriction of $A'_G$ to $\mathcal{J}^{-1}(\mathcal{O}_\rho)$ is tangent to it. Consequently, as $\mathcal{J}^{-1}(\mathcal{O}_\rho)$ is an immersed submanifold of $M$, there exists for each Hamiltonian vector field $X_f\in\mathfrak{X}(M)$, $f\in C^\infty(M)^G$, a vector field $X'_f \in \mathfrak{X}(\mathcal{J}^{-1}(\mathcal{O}_\rho))$ such that $Ti_{\mathcal{O}_\rho} \circ X'_f = X_f \circ i_{\mathcal{O}_\rho}$, with $i_{\mathcal{O}_\rho} : \mathcal{J}^{-1}(\mathcal{O}_\rho) \hookrightarrow M$ the injection. The restriction $A'_G|_{\mathcal{J}^{-1}(\mathcal{O}_\rho)}$ of $A'_G$ to $\mathcal{J}^{-1}(\mathcal{O}_\rho)$ is generated by the vector fields of the form $X'_f$ and it is therefore smooth. It is also integrable since for any point $m = g \cdot z \in \mathcal{J}^{-1}(\mathcal{O}_\rho)$, $z \in \mathcal{J}^{-1}(\rho)$, the embedded submanifold $\mathcal{J}^{-1}(g \cdot \rho)$ of $\mathcal{J}^{-1}(\mathcal{O}_\rho)$ is the maximal integral submanifold of $A'_G|_{\mathcal{J}^{-1}(\mathcal{O}_\rho)}$. This is so because the flows $F_t$ and $F'_t$ of $X_f$ and $X'_f$, respectively, satisfy that $i_{\mathcal{O}_\rho} \circ F'_t = F_t \circ i_{\mathcal{O}_\rho}$. It is then clear that $A'_G|_{\mathcal{J}^{-1}(\mathcal{O}_\rho)}$ has constant rank since $\dim A'_G|_{\mathcal{J}^{-1}(\mathcal{O}_\rho)} = \dim \mathcal{J}^{-1}(\rho)$. This all shows that the leaf space $\mathcal{J}^{-1}(\mathcal{O}_\rho)/A'_G$ is well defined.

In order to show that the leaf space $\mathcal{J}^{-1}(\mathcal{O}_\rho)/A'_G$ is a regular quotient manifold we first notice that

$$\mathcal{J}^{-1}(\mathcal{O}_\rho)/A'_G \simeq (G \times_{G_\rho} \mathcal{J}^{-1}(\rho))/A'_G$$

is in bijection with the quotient $G/G_\rho$ that, by the hypothesis on the closedness of $G_\rho$ is a smooth homogeneous manifold. Take in $M'_{\mathcal{O}_\rho} := \mathcal{J}^{-1}(\mathcal{O}_\rho)/A'_G$ the smooth structure that makes the bijection with $G/G_\rho$ a diffeomorphism. It turns out that that smooth structure is the unique one that makes $M'_{\mathcal{O}_\rho}$ into a regular quotient manifold since it can be readily verified that the map

$$\begin{array}{cccc} \mathcal{J}_{\mathcal{O}_\rho} : & \mathcal{J}^{-1}(\mathcal{O}_\rho) \simeq G \times_{G_\rho} \mathcal{J}^{-1}(\rho) & \longrightarrow & \mathcal{J}^{-1}(\mathcal{O}_\rho)/A'_G \simeq G/G_\rho \\ & [g, z] & \longmapsto & gG_\rho \end{array}$$

is a surjective submersion. ■

We now introduce the regularized polar reduced subspaces of $M'_{\mathcal{O}_\rho}$, available when $M$ is symplectic. We retake the ideas and notations introduced just above (3.8). Let $(\mathcal{J}_{L^\rho}^{-1}(L^\rho \cdot \sigma)/L^\rho, (\omega|_{M_H^\rho})_{L^\rho \cdot \sigma})$ be an orbit regularization of $(\mathcal{J}^{-1}(\mathcal{O}_\rho)/G, \omega_{\mathcal{O}_\rho})$. A straightforward application of Proposition 3.6 implies that the reduced space polar to $(\mathcal{J}_{L^\rho}^{-1}(L^\rho \cdot \sigma)/L^\rho, (\omega|_{M_H^\rho})_{L^\rho \cdot \sigma})$ equals

$$\mathcal{J}_{L^\rho}^{-1}(L^\rho \cdot \sigma)/A'_{L^\rho} = \mathcal{J}^{-1}(\mathcal{N}_\rho)/A'_G$$

which is naturally diffeomorphic to $N(H)^\rho/G_\rho$. We will say that $\mathcal{J}^{-1}(\mathcal{N}_\rho)/A'_G$ is ***a regularized polar reduced subspace*** of $M'_{\mathcal{O}_\rho}$. We will write $M'_{\mathcal{N}_\rho} := \mathcal{J}^{-1}(\mathcal{N}_\rho)/A'_G$ and denote by $\mathcal{J}_{\mathcal{N}_\rho} : \mathcal{J}^{-1}(\mathcal{N}_\rho) \to \mathcal{J}^{-1}(\mathcal{N}_\rho)/A'_G$ the canonical projection. Notice that the spaces $M'_{\mathcal{N}_\rho}$ are embedded submanifolds of $M'_{\mathcal{O}_\rho}$. Finally, the decomposition (3.8) implies that the polar reduced space can be written as the following disjoint union of regularized polar reduced subspaces:

$$M'_{\mathcal{O}_\rho} = \mathcal{J}^{-1}(\mathcal{O}_\rho)/A'_G = \dot{\bigcup}_{[g] \in G/N(H)^\rho} \mathcal{J}^{-1}(\mathcal{N}_{g \cdot \rho})/A'_G = \dot{\bigcup}_{[g] \in G/N(H)^\rho} M'_{\mathcal{N}_{g \cdot \rho}}. \tag{3.10}$$

Equivalently, we have that

$$G/G_\rho = \dot{\bigcup}_{[g] \in G/N(H)^\rho} gN(H)^\rho/G_\rho, \tag{3.11}$$

where the quotient $gN(H)^\rho/G_\rho$ denotes the orbit space of the free and proper action of $G_\rho$ on $gN(H)^\rho$ by $h \cdot gn := gnh$, $h \in G_\rho$, $n \in N(H)^\rho$.

Before we state our next result we need some terminology. We will denote by $C^\infty\left(\mathcal{J}^{-1}(\mathcal{O}_\rho)/A'_G\right)$ the set of smooth real valued functions on $M'_{\mathcal{O}_\rho}$ with the smooth structure introduced in Proposition 3.15. Recall now that, as we pointed out in (2.3), there is a notion of smooth function on $M/A'_G$, namely $C^\infty(M/A'_G) := \{f \in C^0(M/A'_G) \mid f \circ \mathcal{J} \in C^\infty(M)^{A'_G}\}$. Analogously, for each open $A'_G$–invariant subset $U$ of $M$ we can define $C^\infty(U/A'_G) := \{f \in C^0(U/A'_G) \mid f \circ \mathcal{J}|_U \in C^\infty(U)^{A'_G}\}$.We define the set of ***Whitney smooth functions*** $W^\infty\left(\mathcal{J}^{-1}(\mathcal{O}_\rho)/A'_G\right)$ on $\mathcal{J}^{-1}(\mathcal{O}_\rho)/A'_G$ as

$$W^\infty\left(\mathcal{J}^{-1}(\mathcal{O}_\rho)/A'_G\right) := \{f \text{ real function on } M'_{\mathcal{O}_\rho} \mid f = F|_{M'_{\mathcal{O}_\rho}}, \text{ with } F \in C^\infty(M/A'_G)\}.$$

The definitions and the fact that $\mathcal{J}_{\mathcal{O}_\rho}$ is a submersion imply that

$$W^\infty\left(\mathcal{J}^{-1}(\mathcal{O}_\rho)/A'_G\right) \subset C^\infty\left(\mathcal{J}^{-1}(\mathcal{O}_\rho)/A'_G\right).$$

Indeed, let $f \in W^\infty\left(\mathcal{J}^{-1}(\mathcal{O}_\rho)/A'_G\right)$ arbitrary. By definition, there exist $F \in C^\infty(M/A'_G)$ such that $f = F|_{M'_{\mathcal{O}_\rho}}$. As $F \in C^\infty(M/A'_G)$ we have that $F\circ\mathcal{J} \in C^\infty(M)$. Also, as $\mathcal{J}^{-1}(\mathcal{O}_\rho)$ is an immersed initial submanifold of $M$, the injection $i_{\mathcal{O}_\rho} : \mathcal{J}^{-1}(\mathcal{O}_\rho) \hookrightarrow M$ is smooth, and therefore so is $F\circ\mathcal{J}\circ i_{\mathcal{O}_\rho} = F\circ\mathcal{J}_{\mathcal{O}_\rho}$. Consequently, $f \circ \mathcal{J}_{\mathcal{O}_\rho} = F \circ \mathcal{J}_{\mathcal{O}_\rho}$ is smooth. As $\mathcal{J}_{\mathcal{O}_\rho}$ is a submersion $f$ is necessarily smooth, that is, $f \in C^\infty\left(\mathcal{J}^{-1}(\mathcal{O}_\rho)/A'_G\right)$, as required.

**Definition 3.16** *We say that $M'_{\mathcal{O}_\rho}$ is **Whitney spanned** when the differentials of its Whitney smooth functions span its cotangent bundle, that is,*

$$\operatorname{span}\{\mathbf{d}f(\sigma) \mid f \in W^\infty(M'_{\mathcal{O}_\rho})\} = T^*_\sigma M'_{\mathcal{O}_\rho}, \quad \textit{for all} \quad \sigma \in M'_{\mathcal{O}_\rho}.$$

A sufficient (but not necessary!) condition for $M'_{\mathcal{O}_\rho}$ to be Whitney spanned is that $W^\infty(M'_{\mathcal{O}_\rho}) = C^\infty(M'_{\mathcal{O}_\rho})$.

We are now in the position to state the main results of this section.

**Theorem 3.17 (Polar reduction of a Poisson manifold)** *Let $(M, \{\cdot,\cdot\})$ be a smooth Poisson manifold and $G$ be a Lie group acting canonically and properly on $M$. Let $\mathcal{J} : M \to M/A'_G$ be the optimal momentum map associated to this action and $\rho \in M/A'_G$ be such that $G_\rho$ is closed in $G$. If $A_G$ is weakly von Neumann then, for each point $z \in \mathcal{J}^{-1}(\mathcal{O}_\rho)$ and vectors $v, w \in T_z\mathcal{J}^{-1}(\mathcal{O}_\rho)$, there exists an open $A'_G$–invariant neighborhood $U$ of $z$ and two smooth functions $f, g \in C^\infty(U)$ such that $v = X_f(z)$ and $w = X_g(z)$. Moreover, there is a unique presymplectic form $\omega'_{\mathcal{O}_\rho}$ on the polar reduced space $M'_{\mathcal{O}_\rho}$ that satisfies*

$$\{f, g\}|_U(z) = \pi^*_{\mathcal{O}_\rho}\omega_{\mathcal{O}_\rho}(z)(v, w) + \mathcal{J}^*_{\mathcal{O}_\rho}\omega'_{\mathcal{O}_\rho}(z)(v, w) \tag{3.12}$$

*If $M'_{\mathcal{O}_\rho}$ is Whitney spanned then the form $\omega'_{\mathcal{O}_\rho}$ is symplectic.*

**Remark 3.18** It can be proved that when $A_G$ is von Neumann and $A'_G$ satisfies the extension property (see Remark 2.6) the symplecticity of $\omega'_{\mathcal{O}_\rho}$ is equivalent to $M'_{\mathcal{O}_\rho}$ being Whitney spanned. ⧫

When the Poisson manifold $(M, \{\cdot,\cdot\})$ is actually a symplectic manifold with symplectic form $\omega$ the von Neumann condition in the previous result is no longer needed. Moreover, the conditions under which the form $\omega'_{\mathcal{O}_\rho}$ is symplectic can be completely characterized and the regularized polar subspaces appear as symplectic submanifolds of the polar space that contains them.

**Theorem 3.19 (Polar reduction of a symplectic manifold)** *Let $(M, \omega)$ be a smooth symplectic manifold and $G$ be a Lie group acting canonically and properly on $M$. Let $\mathcal{J} : M \to M/A'_G$ be the optimal momentum map associated to this action and $\rho \in M/A'_G$ be such that $G_\rho$ is closed in $G$.*

**(i)** *There is a unique presymplectic form $\omega'_{\mathcal{O}_\rho}$ on the polar reduced space $M'_{\mathcal{O}_\rho} \simeq G/G_\rho$ that satisfies*

$$i^*_{\mathcal{O}_\rho}\omega = \pi^*_{\mathcal{O}_\rho}\omega_{\mathcal{O}_\rho} + \mathcal{J}^*_{\mathcal{O}_\rho}\omega'_{\mathcal{O}_\rho}. \tag{3.13}$$

*The form $\omega'_{\mathcal{O}_\rho}$ is symplectic if and only if for one point $z \in \mathcal{J}^{-1}(\mathcal{O}_\rho)$ (and hence for all) we have that*

$$\mathfrak{g} \cdot z \cap (\mathfrak{g} \cdot z)^\omega \subset T_z M_{G_z} \tag{3.14}$$

**(ii)** *Let $M'_{\mathcal{N}_\rho} = \mathcal{J}^{-1}(\mathcal{N}_\rho)/A'_G \simeq N(H)^\rho/G_\rho$ be a regularized polar reduced subspace of $M'_{\mathcal{O}_\rho}$. Let $j_{\mathcal{N}_\rho} : \mathcal{J}^{-1}(\mathcal{N}_\rho)/A'_G \hookrightarrow \mathcal{J}^{-1}(\mathcal{O}_\rho)/A'_G$ be the inclusion and $\omega'_{\mathcal{O}_\rho}$ the presymplectic form defined in* **(i)**. *Then, the form*

$$\omega'_{\mathcal{N}_\rho} := j^*_{\mathcal{N}_\rho}\omega'_{\mathcal{O}_\rho} \tag{3.15}$$

*is symplectic, that is, the regularized polar subspaces are symplectic submanifolds of the polar space that contains them.*

**Remark 3.20** The characterization (3.14) of the symplecticity of $\omega'_{\mathcal{O}_\rho}$ admits a particularly convenient reformulation when the $G$–action on the symplectic manifold $(M, \omega)$ admits an equivariant momentum map $\mathbf{J} : M \to \mathfrak{g}^*$. Indeed, let $z \in M$ be such that $\mathbf{J}(z) = \mu \in \mathfrak{g}^*$ and $G_z = H$. Then, if the symbol $G_\mu$ denotes the coadjoint isotropy of $\mu$, (3.14) is equivalent to

$$\mathfrak{g} \cdot z \cap (\mathfrak{g} \cdot z)^\omega = \mathfrak{g}_\mu \cdot z \subset T_z M_H,$$

which in turn amounts to $\mathfrak{g}_\mu \cdot z \subset \mathfrak{g}_\mu \cdot z \cap T_z M_H = \mathrm{Lie}(N(H) \cap G_\mu) \cdot z$. Let $N_{G_\mu}(H) := N(H) \cap G_\mu$. With this notation, the condition can be rewritten as $\mathfrak{g}_\mu + \mathfrak{h} \subset \mathrm{Lie}(N_{G_\mu}(H)) + \mathfrak{h} \subset \mathfrak{g}_\mu$ or, equivalently, as

$$\mathfrak{g}_\mu = \mathrm{Lie}(N_{G_\mu}(H)). \tag{3.16}$$

♦

**Proof of Theorem 3.17.** As $A_G$ is weakly von Neumann we have that for any $z \in M$ $\mathfrak{g} \cdot z \subset A''_G(z)$ or, equivalently, that for any $z \in M$ and any $\xi \in \mathfrak{g}$, there is a $A'_G$–invariant neighborhood $U$ of $z$ and a function $F \in C^\infty(U/A'_G)$ such that $\xi_M(z) = X_{F\circ\mathcal{J}}(z)$. Consequently, for any vector $v \in T_z\mathcal{J}^{-1}(\mathcal{O}_\rho)$ there exists $f \in C^\infty(M)^G$ and $F \in C^\infty(U/A'_G)$ (shrink $U$ if necessary) such that $v = X_f(z) + X_{F\circ\mathcal{J}}(z) = X_{f|_U + F\circ\mathcal{J}}(z)$. Let $w \in T_z\mathcal{J}^{-1}(\mathcal{O}_\rho)$, $l \in C^\infty(M)^G$, and $L \in C^\infty(U/A'_G)$ be such that $w = X_l(z) + X_{L\circ\mathcal{J}}(z) = X_{l|_U + L\circ\mathcal{J}}(z)$. Expression (3.12) can then be rewritten as

$$\begin{aligned} \mathcal{J}^*_{\mathcal{O}_\rho}\omega'_{\mathcal{O}_\rho}(z)(v, w) &= \mathcal{J}^*_{\mathcal{O}_\rho}\omega'_{\mathcal{O}_\rho}(z)(X_{f|_U + F\circ\mathcal{J}}(z), X_{l|_U + L\circ\mathcal{J}}(z)) \\ &= \{f + F \circ \mathcal{J}, l + L \circ \mathcal{J}\}|_U(z) - \pi^*_{\mathcal{O}_\rho}\omega_{\mathcal{O}_\rho}(z)(X_{f|_U + F\circ\mathcal{J}}(z), X_{l|_U + L\circ\mathcal{J}}(z)) \\ &= \{F \circ \mathcal{J}, L \circ \mathcal{J}\}|_U(z) \end{aligned} \tag{3.17}$$

We now show that $\omega'_{\mathcal{O}_\rho}$ is well defined. Indeed, let $z' \in \mathcal{J}^{-1}(\mathcal{O}_\rho)$ and $v', w' \in T_{z'}\mathcal{J}^{-1}(\mathcal{O}_\rho)$ be such that $T_z\mathcal{J}_{\mathcal{O}_\rho} \cdot v = T_{z'}\mathcal{J}_{\mathcal{O}_\rho} \cdot v'$ and $T_z\mathcal{J}_{\mathcal{O}_\rho} \cdot w = T_{z'}\mathcal{J}_{\mathcal{O}_\rho} \cdot w'$. First of all these equalities imply the existence of an element $\mathcal{F}_T$ in the polar pseudogroup of $A_G$ such that $z' = \mathcal{F}_T(z)$. As $\mathcal{F}_T$ is a local diffeomorphism that such that $\mathcal{J}_{\mathcal{O}_\rho} \circ \mathcal{F}_T = \mathcal{J}_{\mathcal{O}_\rho}$, we have that $T_z\mathcal{J}_{\mathcal{O}_\rho} = T_{z'}\mathcal{J}_{\mathcal{O}_\rho} \cdot T_z\mathcal{F}_T$. Now, we can rewrite the conditions $T_z\mathcal{J}_{\mathcal{O}_\rho} \cdot v = T_{z'}\mathcal{J}_{\mathcal{O}_\rho} \cdot v'$ and $T_z\mathcal{J}_{\mathcal{O}_\rho} \cdot w = T_{z'}\mathcal{J}_{\mathcal{O}_\rho} \cdot w'$ as $T_{z'}\mathcal{J}_{\mathcal{O}_\rho} \cdot T_z\mathcal{F}_T \cdot v = T_{z'}\mathcal{J}_{\mathcal{O}_\rho} \cdot v'$ and $T_{z'}\mathcal{J}_{\mathcal{O}_\rho} \cdot T_z\mathcal{F}_T \cdot w = T_{z'}\mathcal{J}_{\mathcal{O}_\rho} \cdot w'$, respectively, which implies the existence of two functions $f', l' \in C^\infty(M)^G$ such that

$$\begin{aligned} v' &= T_z\mathcal{F}_T(X_f(z) + X_{F\circ\mathcal{J}}(z)) + X_{f'}(\mathcal{F}_T(z)) \\ w' &= T_z\mathcal{F}_T(X_l(z) + X_{L\circ\mathcal{J}}(z)) + X_{l'}(\mathcal{F}_T(z)) \end{aligned}$$

or, equivalently:

$$\begin{aligned} v' &= X_{f\circ\mathcal{F}_{-T}}(\mathcal{F}_T(z)) + X_{F\circ\mathcal{J}}(\mathcal{F}_T(z)) + X_{f'}(\mathcal{F}_T(z)) \\ w' &= X_{l\circ\mathcal{F}_{-T}}(\mathcal{F}_T(z)) + X_{L\circ\mathcal{J}}(\mathcal{F}_T(z)) + X_{l'}(\mathcal{F}_T(z)). \end{aligned}$$

Therefore, using (3.12), we have that

$$\begin{aligned} \mathcal{J}_{\mathcal{O}_\rho}^*\omega'_{\mathcal{O}_\rho}(z')(v', w') &= \{f \circ \mathcal{F}_{-T} + F \circ \mathcal{J} + f', l \circ \mathcal{F}_{-T} + L \circ \mathcal{J} + l'\}|_V(\mathcal{F}_T(z)) \\ &\quad -\pi_{\mathcal{O}_\rho}^*\omega_{\mathcal{O}_\rho}(\mathcal{F}_T(z))(X_{f\circ\mathcal{F}_{-T}|_V + F\circ\mathcal{J} + f'}(z), X_{l\circ\mathcal{F}_{-T}|_V + L\circ\mathcal{J} + l'}(z)) \\ &= \{F \circ \mathcal{J}, L \circ \mathcal{J}\}|_V(\mathcal{F}_T(z)) = \{F \circ \mathcal{J}, L \circ \mathcal{J}\}|_U(z) \\ &= \mathcal{J}_{\mathcal{O}_\rho}^*\omega'_{\mathcal{O}_\rho}(z)(v, w), \end{aligned}$$

where $V = U \cap \mathcal{F}_T(\mathrm{Dom}(\mathcal{F}_T)) = \mathcal{F}_T(U \cap \mathrm{Dom}(\mathcal{F}_T))$. Hence, the form $\omega'_{\mathcal{O}_\rho}$ is well defined. The closedness and skew symmetric character of $\omega'_{\mathcal{O}_\rho}$ is obtained as a consequence of $\mathcal{J}_{\mathcal{O}_\rho}$ being a surjective submersion, $\omega_{\mathcal{O}_\rho}$ being closed and skew symmetric, and the $\{\cdot, \cdot\}$ being a Poisson bracket. An equivalent fashion to realize this is by writing $\omega'_{\mathcal{O}_\rho}$ in terms of the symplectic structure of the leaves of $M$. Indeed, as $A_G$ is weakly von Neumann, each connected component of $\mathcal{J}^{-1}(\mathcal{O}_\rho)$ lies in a single symplectic leaf of $(M, \{\cdot, \cdot\})$. In order to simplify the exposition suppose that $\mathcal{J}^{-1}(\mathcal{O}_\rho)$ is connected and let $\mathcal{L}_{\mathcal{O}_\rho}$ be the unique symplectic leaf of $M$ that contains it (otherwise one has just to proceed connected component by connected component). Let $i_{\mathcal{L}_{\mathcal{O}_\rho}} : \mathcal{J}^{-1}(\mathcal{O}_\rho) \to \mathcal{L}_{\mathcal{O}_\rho}$ be the natural injection. Given that $i_{\mathcal{O}_\rho} : \mathcal{J}^{-1}(\mathcal{O}_\rho) \to M$ is smooth and $\mathcal{L}_{\mathcal{O}_\rho}$ is an initial submanifold of $M$, the map $i_{\mathcal{L}_{\mathcal{O}_\rho}}$ is therefore smooth. If we denote by $\omega_{\mathcal{L}_{\mathcal{O}_\rho}}$ the symplectic form of the leaf $\mathcal{L}_{\mathcal{O}_\rho}$, expression (3.12) can be rewritten as:

$$i_{\mathcal{L}_{\mathcal{O}_\rho}}^*\omega_{\mathcal{L}_{\mathcal{O}_\rho}} = \pi_{\mathcal{O}_\rho}^*\omega_{\mathcal{O}_\rho} + \mathcal{J}_{\mathcal{O}_\rho}^*\omega'_{\mathcal{O}_\rho}. \tag{3.18}$$

The antisymmetry and closedness of $\omega'_{\mathcal{O}_\rho}$ appears then as a consequence of the antisymmetry and closedness of $\omega_{\mathcal{O}_\rho}$ and $\omega_{\mathcal{L}_{\mathcal{O}_\rho}}$.

It just remains to be shown that if $M'_{\mathcal{O}_\rho}$ is Whitney spanned then the form $\omega'_{\mathcal{O}_\rho}$ is non degenerate. Let $z \in \mathcal{J}^{-1}(\mathcal{O}_\rho)$ and $v \in T_z\mathcal{J}^{-1}(\mathcal{O}_\rho)$ be such that

$$\omega'_{\mathcal{O}_\rho}(\mathcal{J}_{\mathcal{O}_\rho}(z))(T_z\mathcal{J}_{\mathcal{O}_\rho} \cdot v, T_z\mathcal{J}_{\mathcal{O}_\rho} \cdot w) = 0, \quad \text{for all} \quad w \in T_z\mathcal{J}^{-1}(\mathcal{O}_\rho). \tag{3.19}$$

Take now $f \in C^\infty(M)^G$ and $F \in C^\infty(U/A'_G)$ such that $v = X_f(z) + X_{F\circ\mathcal{J}}(z)$. Condition (3.19) is equivalent to having that

$$\omega'_{\mathcal{O}_\rho}(\mathcal{J}_{\mathcal{O}_\rho}(z))(T_z\mathcal{J}_{\mathcal{O}_\rho} \cdot X_{F\circ\mathcal{J}}(z), T_z\mathcal{J}_{\mathcal{O}_\rho} \cdot X_{L\circ\mathcal{J}}(z)) = 0, \tag{3.20}$$

for all $L \in C^\infty(V/A'_G)$ and all open $A'_G$–invariant neighborhoods $V$ of $z$. By (3.17) we can rewrite (3.23) as

$$\{F \circ \mathcal{J}, L \circ \mathcal{J}\}|_{U\cap V}(z) = 0. \tag{3.21}$$

Now, notice that for any $h \in W^\infty(M'_{\mathcal{O}_\rho})$ there exists a function $H \in C^\infty(M/A'_G)$ such that $H|_{M'_{\mathcal{O}_\rho}} = h$. Moreover, by (3.24) we have that:

$$\mathbf{d}h(\mathcal{J}_{\mathcal{O}_\rho}(z)) \cdot (T_z\mathcal{J}_{\mathcal{O}_\rho} \cdot X_{F\circ\mathcal{J}}(z)) = \mathbf{d}(h \circ \mathcal{J}_{\mathcal{O}_\rho})(z) \cdot X_{F\circ\mathcal{J}}(z) = \mathbf{d}(H \circ \mathcal{J})(z) \cdot X_{F\circ\mathcal{J}}(z) = 0.$$

Given that the previous equality holds for any $h \in W^\infty(M'_{\mathcal{O}_\rho})$ and $M'_{\mathcal{O}_\rho}$ is Whitney spanned we have that $T_z\mathcal{J}_{\mathcal{O}_\rho} \cdot X_{F\circ\mathcal{J}}(z) = T_z\mathcal{J}_{\mathcal{O}_\rho} \cdot v = 0$, as required. ■

**Proof of Theorem 3.19.** **(i)** The well definiteness and presymplectic character of $\omega'_{\mathcal{O}_\rho}$ in this case can be obtained as a consequence of Theorem 3.17. This is particularly evident when we think of $\omega'_{\mathcal{O}_\rho}$ as the form characterized by equality (3.18) and we recall that in the symplectic case $\omega_{\mathcal{L}_{\mathcal{O}_\rho}} = \omega$.

It just remains to be shown that the form $\omega'_{\mathcal{O}_\rho}$ is non degenerate if and only if condition (3.14) holds. We proceed by showing first that if condition (3.14) holds for the point $z \in \mathcal{J}^{-1}(\mathcal{O}_\rho)$ then it holds for all the points in $\mathcal{J}^{-1}(\mathcal{O}_\rho)$. We will then prove that (3.14) at the point $z$ is equivalent to the non degeneracy of $\omega'_{\mathcal{O}_\rho}$ at $\mathcal{J}_{\mathcal{O}_\rho}(z)$.

Suppose first that the point $z \in \mathcal{J}^{-1}(\mathcal{O}_\rho)$ is such that $\mathfrak{g}\cdot z \cap (\mathfrak{g}\cdot z)^\omega \subset T_zM_{G_z}$. Notice now that any element in $\mathcal{J}^{-1}(\mathcal{O}_\rho)$ can be written as $\Phi_g(\mathcal{F}_T(z))$ with $g \in G$ and $\mathcal{F}_T$ in the polar pseudogroup of $A_G$. It is easy to show that the relation

$$\mathfrak{g}\cdot(\Phi_g(\mathcal{F}_T(z))) \cap (\mathfrak{g}\cdot(\Phi_g(\mathcal{F}_T(z))))^\omega \subset T_{\Phi_g(\mathcal{F}_T(z))}M_{G_{\Phi_g(\mathcal{F}_T(z))}}$$

is equivalent to $T_z(\Phi\circ\mathcal{F}_T)(\mathfrak{g}\cdot z \cap (\mathfrak{g}\cdot z)^\omega) \subset T_z(\Phi\circ\mathcal{F}_T)M_{G_z}$ and therefore to $\mathfrak{g}\cdot z \cap (\mathfrak{g}\cdot z)^\omega \subset T_zM_{G_z}$.

Let now $v \in T_z\mathcal{J}^{-1}(\mathcal{O}_\rho)$ be such that

$$\omega'_{\mathcal{O}_\rho}(\mathcal{J}_{\mathcal{O}_\rho}(z))(T_z\mathcal{J}_{\mathcal{O}_\rho}\cdot v, T_z\mathcal{J}_{\mathcal{O}_\rho}\cdot w) = 0, \quad \text{for all} \quad w \in T_z\mathcal{J}^{-1}(\mathcal{O}_\rho). \tag{3.22}$$

Take now $f \in C^\infty(M)^G$ and $\xi \in \mathfrak{g}$ such that $v = X_f(z) + \xi_M(z)$. Condition (3.22) is equivalent to having that

$$\omega'_{\mathcal{O}_\rho}(\mathcal{J}_{\mathcal{O}_\rho}(z))(T_z\mathcal{J}_{\mathcal{O}_\rho}\cdot \xi_M(z), T_z\mathcal{J}_{\mathcal{O}_\rho}\cdot \eta_M(z)) = 0, \qquad \text{for all } \eta \in \mathfrak{g} \tag{3.23}$$

which by (3.13) can be rewritten as

$$\omega(z)(\xi_M(z), \eta_M(z)) = 0, \qquad \text{for all } \eta \in \mathfrak{g}, \tag{3.24}$$

and thereby amounts to having that $\xi_M(z) \in \mathfrak{g}\cdot z \cap (\mathfrak{g}\cdot z)^\omega$. Hence, $\omega'_{\mathcal{O}_\rho}(\mathcal{J}_{\mathcal{O}_\rho}(z))$ is non degenerate if and only if $\xi_M(z) \in \ker T_z\mathcal{J}_{\mathcal{O}_\rho} = A'_G(z)$. Suppose now that condition (3.14) holds; then, as $\xi_M(z) \in \mathfrak{g}\cdot z \cap (\mathfrak{g}\cdot z)^\omega$ we have that $\xi_M(z) \in T_zM_{G_z}$. Using (2.1) we can conclude that $\xi_M(z) \in A'_G(z)$, as required. Conversely, suppose that $\omega'_{\mathcal{O}_\rho}$ is symplectic. The previous equalities immediately imply that $\mathfrak{g}\cdot z \cap (\mathfrak{g}\cdot z)^\omega \subset A'_G(z) \subset T_zM_{G_z}$, as required.

**(ii)** The form $\omega'_{\mathcal{N}_\rho}$ is clearly closed and antisymmetric. We now show that it is non degenerate. Recall firs that the tangent space to $T_z\mathcal{J}^{-1}(\mathcal{N}_\rho)$ at a given point $z \in \mathcal{J}^{-1}(\mathcal{N}_\rho)$ is given by the vectors of the form $v = X_f(z)+\xi_M(z)$, with $f \in C^\infty(M)^G$ and $\xi \in \mathrm{Lie}(N(H)^\rho)$. Let $v = X_f(z)+\xi_M(z) \in T_z\mathcal{J}^{-1}(\mathcal{N}_\rho)$ be such that

$$\mathcal{J}^*_{\mathcal{N}_\rho}(j^*_{\mathcal{N}_\rho}\omega'_{\mathcal{O}_\rho})(z)(X_f(z)+\xi_M(z), X_g(z)+\eta_M(z)) = 0, \quad \text{for all } \eta \in \mathrm{Lie}(N(H)^\rho) \text{ and } g \in C^\infty(M)^G.$$

If we plug into the previous expression the definition of the form $\omega'_{\mathcal{O}_\rho}$ we obtain that

$$\omega(z)(\xi_M(z), \eta_M(z)) = 0, \qquad \text{for all } \eta \in \mathrm{Lie}(N(H)^\rho),$$

that is, $\xi_M(z) \in (\mathrm{Lie}(N(H)^\rho) \cdot z) \cap (\mathrm{Lie}(N(H)^\rho) \cdot z)^\omega = (\mathrm{Lie}(N(H)^\rho) \cdot z) \cap (\mathrm{Lie}(N(H)^\rho) \cdot z)^{\omega|_{M_H^\rho}} = (\mathrm{Lie}(N(H)^\rho/H) \cdot z) \cap A'_{N(H)^\rho/H}(z)$, where the last equality follows from (2.1) and the freeness of the natural $N(H)^\rho/H$–action on $M_H^\rho$. We now recall (see Lemma 4.4 in [OR02a]) that any $N(H)^\rho/H$–invariant function on $M_H^\rho$ admits a local extension to a $G$–invariant function on $M$, hence $\xi_M(z) \in (\mathrm{Lie}(N(H)^\rho/H) \cdot z) \cap A'_G(z)$, and consequently $T_z\mathcal{J}_{\mathcal{O}_\rho} \cdot \xi_M(z) = T_z\mathcal{J}_{\mathcal{O}_\rho} \cdot v = T_z\mathcal{J}_{\mathcal{N}_\rho} \cdot v = 0$, as required. ■

## 3.6 Symplectic leaves and the reduction diagram

Suppose that $A'_G$ is completable so that the symplectic leaves of $M/A'_G$ are well defined. We recall that this is automatically the case when $(M, \omega)$ is symplectic and the $G$–group action is proper (see [O02]). Assume also that $A_G$ is von Neumann so that the diagram $(M/G, \{\cdot,\cdot\}_{M/A_G}) \stackrel{\pi_{A_G}}{\leftarrow} (M, \{\cdot,\cdot\}) \stackrel{\mathcal{J}}{\rightarrow} (M/A'_G, \{\cdot,\cdot\}_{M/A'_G})$ constitutes a dual pair.

Notice that by Definition 2.8, the symplectic leaves of $M/A_G$ and $M/A'_G$ coincide with the connected components of the orbit reduced spaces $M_{\mathcal{O}_\rho}$ and polar reduced spaces $M'_{\mathcal{O}_\rho}$, that we studied in sections 3.4 and 3.5, respectively. We saw that whenever $G_\rho$ is closed in $G$ and the Whitney spanning condition is satisfied these spaces are actual symplectic manifolds. When $M$ is symplectic, the symplecticity of the leaves of $M/A'_G$ is characterized by condition (3.14) or even by (3.16), provided that the $G$–action has a standard equivariant momentum map $\mathbf{J} : M \to \mathfrak{g}^*$ associated. Moreover, when $M_{\mathcal{O}_\rho}$ and $M'_{\mathcal{O}_\rho}$ are corresponding leaves, their symplectic structures are connected to each other by an identity that naturally generalizes the classical relation that we recalled in (3.2).

The following diagram represents all the spaces that we worked with and their relations. The part of the diagram dealing with the regularized spaces refers only to the situation in which $M$ is symplectic.

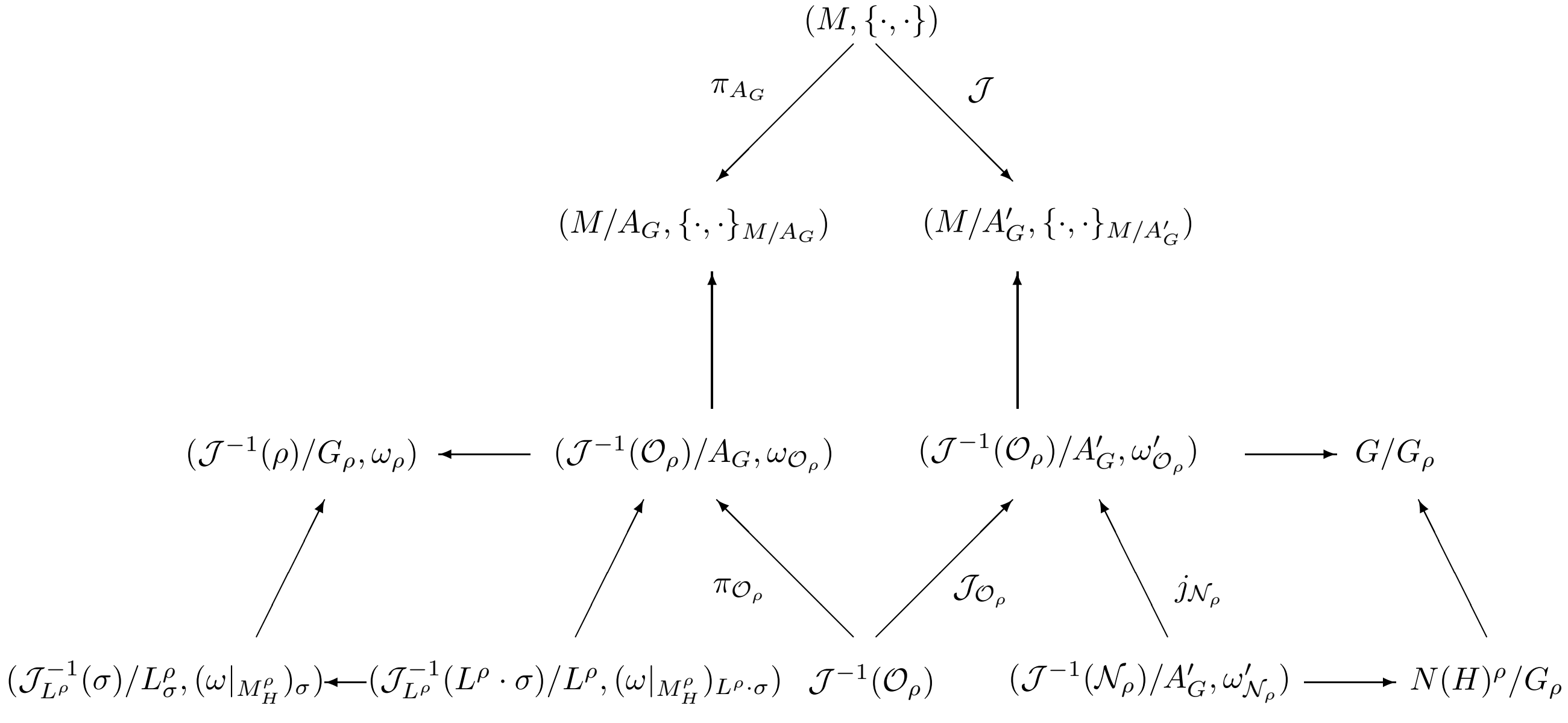

## 3.7 Orbit reduction using the standard momentum map. Beyond compact groups

The approach to optimal orbit reduction developed in the last few sections sheds some light on how to carry out orbit reduction with a standard momentum map when the symmetry group is not compact. This absence of compactness poses some technical problems that have been tackled by various people over the years using different approaches. Since these problems already arise in the free actions case we will restrict ourselves to this situation. More specifically we will assume that we have a Lie group $G$ (not necessarily compact) acting freely and canonically on the symplectic manifold $(M, \omega)$. We will suppose that this action has a coadjoint equivariant momentum map $\mathbf{J} : M \rightarrow \mathfrak{g}^*$ associated. For the sake of simplicity in the exposition and in order to have a better identification with the material presented in the previous sections we will assume that $\mathbf{J}$ has connected fibers. This assumption is not fundamental. The reader interested in the general case with no connectedness hypothesis in the fibers and non free actions may want to check with [OR02b].

In the presence of the hypotheses that we just stated, the momentum map $\mathbf{J}$ is a submersion that maps $M$ onto an open coadjoint equivariant subset $\mathfrak{g}^*_{\mathbf{J}}$ of $\mathfrak{g}^*$. Moreover, any value $\mu \in \mathfrak{g}^*_{\mathbf{J}}$ of $\mathbf{J}$ is regular and has a smooth Marsden–Weinstein symplectic reduced space $\mathbf{J}^{-1}(\mu)/G_\mu$ associated. What about the orbit reduced space $\mathbf{J}^{-1}(\mathcal{O}_\mu)/G$? When the Lie group $G$ is compact there is no problem to canonically endow $\mathbf{J}^{-1}(\mathcal{O}_\mu)/G$ with a smooth structure. Indeed, in this case the coadjoint orbit $\mathcal{O}_\mu$ is an embedded submanifold of $\mathfrak{g}^*$ transverse to the momentum mapping. The Transversal Mapping Theorem ensures that $\mathbf{J}^{-1}(\mathcal{O}_\mu)$ is a $G$–invariant embedded submanifold of $M$ and hence the quotient $\mathbf{J}^{-1}(\mathcal{O}_\mu)/G$ is smooth and symplectic with the form spelled out in (3.1). In the non compact case this argument breaks down due to the non embedded character of $\mathcal{O}_\mu$ in $\mathfrak{g}^*$. In trying to fix this problem this has lead to the assumption of locally closedness on the coadjoint orbits that one can see in a number of papers (see for instance [BL97]). Nevertheless, this hypothesis is not needed to carry out point reduction, and therefore makes the two approaches non equivalent. The first work where this hypothesis has been eliminated is [CS01]. In this paper the authors use a combination of distribution theory with Sikorski differential spaces to show that the orbit reduced space is a symplectic manifold. Nevertheless, the first reference where the standard formula (3.1) appears at this level of generality is [Bl01]. In that paper the author only deals with the free case. Nevertheless the use of a standard technique of reduction to the isotropy type manifolds that the reader can find in [SL91, O98, CS01, OR02b] generalizes the results of [Bl01] to singular situations.

In the next few paragraphs we will illustrate Theorem 3.17 by showing that the results in [CS01, Bl01] can be obtained as a corollary of it.

We start by identifying in this setup all the elements in that result. First of all, we have that the polar distribution satisfies $A'_G = \ker T\mathbf{J}$ (see [OR02a]) and the connectedness hypothesis on the fibers of $\mathbf{J}$ implies that the optimal momentum map $\mathcal{J} : M \rightarrow M/A'_G$ in this case can be identified with $\mathbf{J} : M \rightarrow \mathfrak{g}^*_{\mathbf{J}}$. This immediately implies that for any $\mu \in \mathfrak{g}^*_{\mathbf{J}} \simeq M/A'_G$, the isotropy $G_\mu$ is closed in $G$ and, by Theorem 3.11 there is a unique smooth structure on $\mathbf{J}^{-1}(\mathcal{O}_\mu)$ that makes it into an initial submanifold of $M$ and, at the same time, an integral manifold of the distribution $D = A'_G + \mathfrak{g} \cdot m = \ker T\mathbf{J} + \mathfrak{g} \cdot m$. This structure coincides with the one given in [Bl01]. Also, by Theorem 3.13, the quotient $\mathbf{J}^{-1}(\mathcal{O}_\mu)/G$ admits a unique symplectic structure $\omega_{\mathcal{O}_\mu}$ that makes it symplectomorphic to the Marsden–Weinstein point reduced space $(\mathbf{J}^{-1}(\mu)/G_\mu, \omega_\mu)$. It remains to be shown that we can use (3.13) in this case and that the resulting formula coincides with the standard one (3.1) provided by [Bl01]. An analysis of the polar reduced space in this setup will provide an affirmative answer to this question.

By Proposition 3.15 the polar reduced space $\mathbf{J}(\mathcal{O}_\mu)/A'_G$ is endowed with the only smooth structure that makes it diffeomorphic to the homogeneous space $G/G_\mu \simeq \mathcal{O}_\mu$. Hence, in this case $\mathbf{J}_{\mathcal{O}_\mu} : \mathbf{J}^{-1}(\mathcal{O}_\mu) \rightarrow$ $\mathcal{O}_\mu$ is the map given by $\mathbf{J}_{\mathcal{O}_\mu}(z) := \mathbf{J}(z)$ which is smooth because the coadjoint orbits are always initial submanifolds of $\mathfrak{g}^*$. Therefore we can already compute the polar symplectic form $\omega'_{\mathcal{O}_\mu}$. By (3.13) we have that for any $\xi, \eta \in \mathfrak{g}$ and any $z \in \mathbf{J}^{-1}(\mathcal{O}_\mu)$ (for simplicity in the exposition we take $\mathbf{J}(z) = \mu$):

$$\mathbf{J}^*_{\mathcal{O}_\mu}\omega'_{\mathcal{O}_\mu}(z)(\xi_M(z), \eta_M(z)) = i^*_{\mathcal{O}_\mu}\omega(z)(\xi_M(z), \eta_M(z)) - \pi^*_{\mathcal{O}_\mu}\omega_{\mathcal{O}_\mu}(z)(\xi_M(z), \eta_M(z)),$$

or, equivalently:

$$\omega'_{\mathcal{O}_\mu}(\mu)(\mathrm{ad}^*_\xi\mu, \mathrm{ad}^*_\eta\mu) = \omega(z)(\xi_M(z), \eta_M(z)) = \langle \mathbf{J}(z), [\xi, \eta]\rangle = \langle \mu, [\xi, \eta]\rangle.$$

In conclusion, in this case the polar reduced form $\omega'_{\mathcal{O}_\mu}$ coincides with the "+"–Kostant–Kirillov–Souriau symplectic form on the coadjoint orbit $\mathcal{O}_\mu$. Therefore, the general optimal orbit reduction formula (3.13) coincides with the standard one (3.1).

## 3.8 Examples: the polar reduction of the coadjoint action

We now provide two examples on how we can use the coadjoint action along with Theorems 3.17 and 3.19 to easily produce symplectic manifolds and symplectically decomposed presymplectic manifolds.

### 3.8.1 The coadjoint orbits as polar reduced spaces

Let $G$ be a Lie group, $\mathfrak{g}$ be its Lie algebra, and $\mathfrak{g}^*$ be its dual considered as a Lie–Poisson space. In this elementary example we show how the coadjoint orbits appear as the polar reduced spaces of the coadjoint $G$–action on $\mathfrak{g}^*$.

A straightforward computation shows that the coadjoint action of $G$ on the Lie–Poisson space $\mathfrak{g}^*$ is canonical. Moreover, the polar distribution $A'_G(\mu) = 0$ for all $\mu \in \mathfrak{g}^*$ and therefore the optimal momentum map $\mathcal{J} : \mathfrak{g}^* \to \mathfrak{g}^*$ is the identity map on $\mathfrak{g}^*$. This immediately implies that any open set $U \subset \mathfrak{g}^*$ is $A'_G$–invariant, that $C^\infty(U)^{A'_G} = C^\infty(U)$, and that therefore $\mathfrak{g} \cdot \mu \subset A''_G(\mu)$, for any $\mu \in \mathfrak{g}^*$. The coadjoint action on $\mathfrak{g}^*$ is therefore weakly von Neumann (actually, if $G$ is connected $A_G$ is von Neumann).

We now look at the corresponding reduced spaces. On one hand the orbit reduced spaces $\mathcal{J}^{-1}(\mathcal{O}_\rho)/G$ are the quotients $G\cdot\mu/G$ and therefore amount to points. At the same time, we have that $\mathcal{J}^{-1}(\mathcal{O}_\rho)/A'_G = \mathcal{O}_\mu/A'_G = \mathcal{O}_\mu$, that is, the polar reduced spaces are the coadjoint orbits which, by Theorem 3.17, are symplectic. Indeed, the Whitney spanning condition necessary for the application of this result is satisfied since in this case $\mathrm{span}\{\mathbf{d}f(\mu) \mid f \in W^\infty(M'_{\mathcal{O}_\rho})\} = \mathrm{span}\{\mathbf{d}h|_{\mathcal{O}_\mu}(\mu) \mid h \in C^\infty(\mathfrak{g}^*)\} = T^*_\mu\mathcal{O}_\mu$. Note that the last equality is a consequence of the immersed character of the coadjoint orbits $\mathcal{O}_\mu$ as submanifolds of $\mathfrak{g}^*$ (the equality is easily proved using immersion charts around the point $\mu$).

### 3.8.2 Symplectic decomposition of presymplectic homogeneous manifolds

Let $G$ be a Lie group, $\mathfrak{g}$ be its Lie algebra, and $\mathfrak{g}^*$ be its dual. Let $\mathcal{O}_{\mu_1}$ and $\mathcal{O}_{\mu_2}$ be two coadjoint orbits of $\mathfrak{g}^*$ that we will consider as symplectic manifolds endowed with the KKS–symplectic forms $\omega_{\mathcal{O}_{\mu_1}}$ and $\omega_{\mathcal{O}_{\mu_2}}$, respectively. The cartesian product $\mathcal{O}_{\mu_1} \times \mathcal{O}_{\mu_2}$ is also a symplectic manifold with the sum symplectic form $\omega_{\mathcal{O}_{\mu_1}} + \omega_{\mathcal{O}_{\mu_2}}$. The diagonal action of $G$ on $\mathcal{O}_{\mu_1} \times \mathcal{O}_{\mu_2}$ is canonical with respect to this symplectic structure and, moreover, it has a standard equivariant momentum map $\mathbf{J} : \mathcal{O}_{\mu_1} \times \mathcal{O}_{\mu_2} \to \mathfrak{g}^*$ associated given by $\mathbf{J}(\nu, \eta) = \nu + \eta$. We now suppose that this action is proper and we will study, in this particular case, the orbit and polar reduced spaces introduced in the previous sections.

We start by looking at the level sets of the optimal momentum map $\mathcal{J} : \mathcal{O}_{\mu_1} \times \mathcal{O}_{\mu_2} \to \mathcal{O}_{\mu_1} \times \mathcal{O}_{\mu_2}/A'_G$. A general result (see Theorem 3.6 in [OR02a]) states that in the presence of a standard momentum map the fibers of the optimal momentum map coincide with the connected components of the intersections of the level sets of the momentum map with the isotropy type submanifolds. Hence, in our case, if $\rho = \mathcal{J}(\mu_1, \mu_2)$, we have that

$$\mathcal{J}^{-1}(\rho) = (\mathbf{J}^{-1}(\mu_1 + \mu_2) \cap (\mathcal{O}_{\mu_1} \times \mathcal{O}_{\mu_2})_{G_{(\mu_1,\mu_2)}})_c, \tag{3.25}$$

where the subscript $c$ in the previous expression stands for the connected component of $\mathbf{J}^{-1}(\mu_1+\mu_2) \cap (\mathcal{O}_{\mu_1} \times \mathcal{O}_{\mu_2})_{G_{(\mu_1,\mu_2)}}$ that contains $\mathcal{J}^{-1}(\rho)$. Given that the isotropy $G_{(\mu_1,\mu_2)} = G_{\mu_1} \cap G_{\mu_2}$, with $G_{\mu_1}$ and $G_{\mu_2}$ the coadjoint isotropies of $\mu_1$ and $\mu_2$, respectively, the expression (3.25) can be rewritten as

$$\begin{aligned}\mathcal{J}^{-1}(\rho) = (\{(\mathrm{Ad}^*_{g^{-1}}\mu_1, \mathrm{Ad}^*_{h^{-1}}\mu_2) \mid g, h \in G, \text{ such that} \\ \mathrm{Ad}^*_{g^{-1}}\mu_1 + \mathrm{Ad}^*_{h^{-1}}\mu_2 = \mu_1 + \mu_2,\ gG_{\mu_1}g^{-1} \cap hG_{\mu_2}h^{-1} = G_{\mu_1} \cap G_{\mu_2}\})_c.\end{aligned}$$

It is easy to show that in this case

$$G_\rho = N_{G_{\mu_1+\mu_2}}(G_{\mu_1} \cap G_{\mu_2})^c, \tag{3.26}$$

where the superscript $c$ denotes the closed subgroup of $N_{G_{\mu_1+\mu_2}}(G_{\mu_1} \cap G_{\mu_2}) := N(G_{\mu_1} \cap G_{\mu_2}) \cap G_{\mu_1+\mu_2}$ that leaves $\mathcal{J}^{-1}(\rho)$ invariant. Theorems 3.1 and 3.13 guarantee that the quotients $\mathcal{J}^{-1}(\rho)/G_\rho \simeq \mathcal{J}^{-1}(\mathcal{O}_\rho)/G$ are symplectic. Nevertheless, we will focus our attention in the corresponding polar reduced spaces.

According to Theorem 3.19 and to (3.26), the polar reduced space corresponding to $\mathcal{J}^{-1}(\mathcal{O}_\rho)/G$ is the homogeneous presymplectic manifold

$$G/N_{G_{\mu_1+\mu_2}}(G_{\mu_1} \cap G_{\mu_2})^c. \tag{3.27}$$

Expression (3.16) states that $G/N_{G_{\mu_1+\mu_2}}(G_{\mu_1} \cap G_{\mu_2})^c$ is symplectic if and only if

$$\mathfrak{g}_{\mu_1+\mu_2} = \mathrm{Lie}(N_{G_{\mu_1+\mu_2}}(G_{\mu_1} \cap G_{\mu_2})),$$

which is obviously true when, for instance, $G_{\mu_1} \cap G_{\mu_2}$ is a normal subgroup of $G_{\mu_1+\mu_2}$. In any case, using (3.11) we can write the polar reduced space (3.27) as a disjoint union of its regularized symplectic reduced subspaces that, that in this case are of the form $gN(G_{\mu_1} \cap G_{\mu_2})^\rho/N_{G_{\mu_1+\mu_2}}(G_{\mu_1} \cap G_{\mu_2})^c$ with $g \in G$ and where the superscript $\rho$ denotes the closed subgroup of $N(G_{\mu_1} \cap G_{\mu_2})$ that leaves invariant the connected component of $(\mathcal{O}_{\mu_1} \times \mathcal{O}_{\mu_2})_{G_{\mu_1} \cap G_{\mu_2}}$ that contains $\mathcal{J}^{-1}(\rho)$. More explicitly, we can write the following symplectic decomposition of the polar reduced space:

$$G/N_{G_{\mu_1+\mu_2}}(G_{\mu_1} \cap G_{\mu_2})^c = \dot{\bigcup}_{[g] \in G/N(G_{\mu_1} \cap G_{\mu_2})^\rho} gN(G_{\mu_1} \cap G_{\mu_2})^\rho/N_{G_{\mu_1+\mu_2}}(G_{\mu_1} \cap G_{\mu_2})^c.$$

What we just did in the previous paragraphs for two coadjoint orbits can be inductively generalized to $n$ orbits. We collect the results of that construction under the form of a proposition.

**Proposition 3.21** *Let $G$ be a Lie group, $\mathfrak{g}$ be its Lie algebra, and $\mathfrak{g}^*$ be its dual. Let $\mu_1, \ldots, \mu_n \in \mathfrak{g}^*$. Then, the homogeneous manifold*

$$G/N_{G_{\mu_1+\cdots+\mu_n}}(G_{\mu_1} \cap \ldots \cap G_{\mu_n})^c \tag{3.28}$$

*has a natural presymplectic structure that is nondegenerate if and only if*

$$\mathfrak{g}_{\mu_1+\cdots+\mu_n} = \mathrm{Lie}(N_{G_{\mu_1+\cdots+\mu_n}}(G_{\mu_1} \cap \ldots \cap G_{\mu_n})).$$

*Moreover, (3.28) can be written as a the following disjoint union of symplectic submanifolds*

$$G/N_{G_{\mu_1+\cdots+\mu_n}}(G_{\mu_1}\cap\ldots\cap G_{\mu_n})^c$$
$$=\dot{\bigcup}_{[g]\in G/N(G_{\mu_1}\cap\ldots\cap G_{\mu_n})^\rho} gN(G_{\mu_1}\cap\ldots\cap G_{\mu_n})^\rho/N_{G_{\mu_1+\cdots+\mu_n}}(G_{\mu_1}\cap\ldots\cap G_{\mu_n})^c.$$

# 4 Optimal reduction by stages

As we already described in the introduction, the reduction by stages procedure consists of carrying out reduction in two shots using the normal subgroups of the symmetry group. To be more specific, suppose that we are in the same setup as Theorem 3.1 and that the symmetry group $G$ has a closed normal subgroup $N$. In this section we will spell out the conditions under which reduction by $G$ renders the same result as reduction in the following two stages: we first reduce by $N$; the resulting space inherits symmetry properties coming from the quotient Lie group $G/N$ that can be used to reduce one more time.

In the presence of an equivariant momentum map and freeness in the $G$–action this procedure has been studied in [MMPR98, MMOPR02]. We will extend the results in those papers to the optimal setup and, as a byproduct, we will obtain a generalization to the singular case (non free actions) of the reduction by stages theorem in the presence of an standard equivariant momentum map.

## 4.1 The polar distribution of a normal subgroup

All along this section we will work on a Poisson manifold $(M,\{\cdot,\cdot\})$ acted properly and canonically upon by a Lie group $G$. We will assume that $G$ has a closed normal subgroup that we will denote by $N$. The closedness of $N$ implies that the $N$–action on $M$ by restriction is still proper and that $G/N$ is a Lie group when considered as a homogenous manifold. We will denote by $A'_G$ and $A'_N$ the polar distributions associated to the $G$ and $N$–actions, respectively, and by $\mathcal{J}_G : M\to M/A'_G$ and $\mathcal{J}_N : M\to M/A'_N$ the corresponding optimal momentum maps.

The following proposition provides a characterization of the conditions under which the polar distribution $A'_H$ associated to a closed subgroup $H$ of $G$ is invariant under the lifted action of $G$ to the tangent bundle $TM$.

**Proposition 4.1** *Let $(M,\{\cdot,\cdot\})$ be a Poisson manifold acted properly and canonically upon by a Lie group $G$ via the map $\Phi : G\times M\to M$. Let $H$ be a closed Lie subgroup of $G$. Then:*

**(i)** *The lifted action of $G$ to the tangent bundle $TM$ leaves the $H$–polar distribution $A'_H$ invariant if and only if $f\circ\Phi_{g^{-1}}\in C^\infty(M)^H$, for any $f\in C^\infty(M)^H$ and any $g\in G$. This condition holds if and only if for all $g\in G$, $h\in H$, and $m\in M$, there exists an element $h'\in H$ such that*

$$gh\cdot m = h'g\cdot m. \tag{4.1}$$

**(ii)** *If $G$ acts on $A'_H$ so it does on the corresponding momentum space $M/A'_H$ with a natural action that makes the $H$–optimal momentum map $\mathcal{J}_H : M\to M/A'_H$ $G$–equivariant.*

**(iii)** *For any $m\in M$ we have that $A'_G(m)\subset A'_H(m)$. There is consequently a natural projection $\pi_H : M/A'_G\to M/A'_H$ such that*

$$\mathcal{J}_H = \pi_H\circ\mathcal{J}_G. \tag{4.2}$$

*Moreover, if $G$ acts on $A'_H$ and consequently on $M/A'_H$, the map $\pi_H$ is $G$–equivariant.*

**Proof.** **(i)** Since $H$ is closed in $G$, its action on $M$ by restriction of the $G$–action is still proper. Therefore, $A'_H = \{X_f \mid f \in C^\infty(M)^H\}$. Given that for any $f \in C^\infty(M)^H$ and any $g \in G$ we have that $T\Phi_g \circ X_f = X_{f\circ\Phi_{g^{-1}}} \circ \Phi_g$, we can conclude that the polar distribution $A'_H$ is $G$–invariant iff $f \circ \Phi_{g^{-1}} \in C^\infty(M)^H$, for any $f \in C^\infty(M)^H$ and any $g \in G$. We now check that this condition is equivalent to (4.1).

First of all suppose that $f \circ \Phi_g \in C^\infty(M)^H$, for all $f \in C^\infty(M)^H$ and $g \in G$. Consequently, if we take $m \in M$ and $h \in H$ arbitrary, we have that $f(gh\cdot m) = f(g\cdot m)$. Since the $H$–action on $M$ is proper, the set $C^\infty(M)^H$ of $H$–invariant functions on $M$ separates the $H$–orbits. Therefore, the points $gh \cdot m$ and $g \cdot m$ are in the same $H$–orbit and hence there exists an element $h' \in H$ such that $gh \cdot m = h'g \cdot m$.

Conversely, suppose that for all $g \in G$, $h \in H$, and $m \in M$, there exists an element $h' \in H$ such that $gh \cdot m = h'g \cdot m$. Then, if $f \in C^\infty(M)^H$ we have that

$$f \circ \Phi_g(h \cdot m) = f(gh \cdot m) = f(h'g \cdot m) = f(g \cdot m) = f \circ \Phi_g(m).$$

Consequently, $f \circ \Phi_g \in C^\infty(M)^H$, as required.

**(ii)** Suppose that the lifted action of $G$ to the tangent bundle $TM$ leaves the $H$–polar distribution $A'_H$ invariant. We define the action $G \times M/A'_H \to M/A'_H$ by $g \cdot \mathcal{J}_H(m) := \mathcal{J}_H(g \cdot m)$. It is clearly a left action so all we have to do is showing that it is well defined. Indeed, let $m' \in M$ be such that $m' = \mathcal{F}_T(m)$, with $\mathcal{F}_T \in G_{\mathcal{A}'_H}$. For the sake of simplicity in the exposition suppose that $\mathcal{F}_T = F_T$ with $F_T$ the Hamiltonian flow associated to $f \in C^\infty(M)^H$. Then, for any $g \in G$ we have that

$$g \cdot \mathcal{J}_H(m') = \mathcal{J}_H(g \cdot F_T(m)) = \mathcal{J}_H\left(G_T^{f\circ\Phi_{g^{-1}}}(g \cdot m)\right) = \mathcal{J}_H(g \cdot m) = g \cdot \mathcal{J}_H(m),$$

where $G_T^{f\circ\Phi_{g^{-1}}}$ is the Hamiltonian flow associated to the function $f \circ \Phi_{g^{-1}}$ that, by the hypothesis on the $G$–invariance of $A'_H$, is $H$–invariant and guarantees the equality $\mathcal{J}_H\left(G_T^{f\circ\Phi_{g^{-1}}}(g \cdot m)\right) = \mathcal{J}_H(g \cdot m)$.

**(iii)** The inclusion $A'_G(m) \subset A'_H(m)$ is a direct consequence of the definition of the polar distributions and it implies that each maximal integral leaf of $A'_G$ is included in a single maximal integral leaf of $A'_H$. This feature constitutes the definition of $\pi_H$ that assigns to each leaf in $M/A'_G$ the unique leaf in $M/A'_H$ in which it is sitting. With this definition it is straightforward that $\mathcal{J}_H = \pi_H \circ \mathcal{J}_G$. Now, if $G$ acts on $A'_H$ the map $\mathcal{J}_H$ is $G$–equivariant by part **(ii)**. The $G$–equivariance of $\mathcal{J}_G$ plus the relation $\mathcal{J}_H = \pi_H \circ \mathcal{J}_G$ implies that $\pi_H$ is $G$–invariant. ■

**Remark 4.2** If $H$ is normal in $G$ then, condition (4.1) is trivially satisfied and therefore $G$ acts on $A'_H$. Conversely, if $G$ acts on $A'_H$ and the identity element is an isotropy subgroup of the $G$–action on $M$ then $H$ is necessarily normal in $G$. Indeed, in that case for any $m \in M$, $g \in G$, and $h \in H$, there exists an element $h' \in H$ such that $gh \cdot m = h'g \cdot m$. In particular, if we take an element $m \in M_{\{e\}}$ we have that $gh = h'g$ or, equivalently that $gHg^{-1} \subset H$, for all $g \in G$, which implies that $H$ is normal in $G$. ♦

For future reference we state in the following corollary the claims of Proposition 4.1 in the particular case in which $H$ is a normal subgroup of $G$.

**Corollary 4.3** *Let $(M, \{\cdot,\cdot\})$ be a Poisson manifold acted properly and canonically upon by a Lie group $G$. Let $N$ be a closed normal Lie subgroup of $G$. Then:*

**(i)** *The group $G$ acts on $A'_N$ and on the corresponding momentum space $M/A'_N$ with a natural action that makes the $N$–optimal momentum map $\mathcal{J}_N : M \to M/A'_N$ $G$–equivariant.*

**(ii)** *There is a natural $G$–equivariant projection $\pi_N : M/A'_G \to M/A'_N$ such that $\mathcal{J}_N = \pi_N \circ \mathcal{J}_G$.*

## 4.2 Isotropy subgroups and quotient groups

In this section we introduce the relevant groups and spaces for optimal reduction in two stages.

**Lemma 4.4** *Let $(M, \{\cdot,\cdot\})$ be a Poisson manifold acted properly and canonically upon by a Lie group $G$. Let $N$ be a closed normal Lie subgroup of $G$. Let $\rho \in M/A'_G$ and $\nu := \pi_N(\rho) \in M/A'_N$.*

**(i)** *Let $G_\rho$ and $G_\nu$ be the isotropy subgroups of $\rho \in M/A'_G$ and $\nu := \pi_N(\rho) \in M/A'_N$ with respect to the $G$–actions on $M/A'_G$ and $M/A'_N$, respectively. Then, $G_\rho \subset G_\nu$.*

**(ii)** *Let $N_\nu$ be the $N$–isotropy subgroup of $\nu \in M/A'_N$. Then $N_\nu = N \cap G_\nu$ and $N_\nu$ is normal in $G_\nu$.*

**(iii)** *Endow $N_\nu$ and $G_\nu$ with the unique smooth structures that make them into initial Lie subgroups of $G$. Then, $N_\nu$ is closed in $G_\nu$ and therefore the quotient $H_\nu := G_\nu/N_\nu$ is a Lie group.*

**Proof.** **(i)** It is a consequence of the $G$–equivariance of the projection $\pi_N : M/A'_G \to M/A'_N$. **(ii)** It is straightforward. **(iii)** Let $A$ and $B$ two subsets of a smooth manifold $M$ such that $A \subset B \subset M$. It can be checked by simply using the definition of initial submanifold that if $A$ and $B$ are initial submanifolds of $M$ then $A$ is an initial submanifold of $B$. In our setup, this fact implies that $N_\nu$ is an initial Lie subgroup of $G_\nu$. We actually check that it is a closed Lie subgroup of $G_\nu$. Indeed, let $g \in G_\nu$ be an element in the closure of $N_\nu$ in $G_\nu$. Let $\{g_n\}_{n\in\mathbb{N}} \subset N_\nu$ be a sequence of elements in $N_\nu$ that converges to $g$ in the topology of $G_\nu$. As $G_\nu$ is initial in $G$ we have that $g_n \to g$ also in the topology of $G$. Now, as $\{g_n\}_{n\in\mathbb{N}} \subset N$ and $N$ is closed in $G$, $g \in N$ necessarily. Hence $g \in N \cap G_\nu = N_\nu$, as required. ■

Suppose now that the value $\nu \in M/A'_N$ is such that the action of $N_\nu$ on the level set $\mathcal{J}_N^{-1}(\nu)$ is proper. We emphasize that this property is not automatically inherited from the properness of the $N$–action on $M$. Theorem 3.1 guarantees in that situation that the orbit space $M_\nu := \mathcal{J}_N^{-1}(\nu)/N_\nu$ is a smooth symplectic regular quotient manifold with symplectic form $\omega_\nu$ defined by:

$$\pi_\nu^*\omega_\nu(m)(X_f(m), X_h(m)) = \{f, h\}(m), \text{ for any } m \in \mathcal{J}_N^{-1}(\nu) \text{ and any } f, h \in C^\infty(M)^N.$$

As customary $\pi_\nu : \mathcal{J}_N^{-1}(\nu) \to \mathcal{J}_N^{-1}(\nu)/N_\nu$ denotes the canonical projection and $i_\nu : \mathcal{J}_N^{-1}(\nu) \hookrightarrow M$ the inclusion. We will refer to the pair $(M_\nu, \omega_\nu)$ as the ***first stage reduced space***.

**Proposition 4.5** *Let $(M, \{\cdot,\cdot\})$ be a Poisson manifold acted properly and canonically upon by a Lie group $G$ via the map $\Phi : G \times M \to M$. Let $N$ be a closed normal Lie subgroup of $G$. Let $\rho = \mathcal{J}_G(m) \in M/A'_G$, for some $m \in M$, and $\nu := \pi_N(\rho) = \mathcal{J}_N(m) \in M/A'_N$.*

**(i)** *If the Lie group $N_\nu$ acts properly on the level set $\mathcal{J}_N^{-1}(\nu)$ then the Lie group $H_\nu := G_\nu/N_\nu$ acts smoothly and canonically on the first stage reduced space $(\mathcal{J}_N^{-1}(\nu)/N_\nu, \omega_\nu)$ via the map*

$$gN_\nu \cdot \pi_\nu(m) := \pi_\nu(g \cdot m), \tag{4.3}$$

*for all $gN_\nu \in H_\nu$ and $m \in \mathcal{J}_N^{-1}(\nu)$.*

**(ii)** *Suppose that $N_\nu$ and $H_\nu$ act properly on $\mathcal{J}_N^{-1}(\nu)$ and $M_\nu$, respectively. Let $\mathcal{J}_{H_\nu} : M_\nu \to M_\nu/A'_{H_\nu}$ be the optimal momentum map associated to the $H_\nu$–action on $M_\nu = \mathcal{J}_N^{-1}(\nu)/N_\nu$ and $\sigma = \mathcal{J}_{H_\nu}(\pi_\nu(m))$. Then,*

$$\mathcal{J}_{H_\nu}\left(\pi_\nu(\mathcal{J}_G^{-1}(\rho))\right) = \sigma. \tag{4.4}$$

**Proof.** **(i)** We first show that the action given by expression (4.3) is well defined and is smooth. The action $\varphi^\nu : G_\nu \times \mathcal{J}_N^{-1}(\nu) \to \mathcal{J}_N^{-1}(\nu)$ obtained by restriction of the domain and range of $\Phi$ is smooth since $G_\nu$ and $\mathcal{J}_N^{-1}(\nu)$ are initial submanifolds of $G$ and $M$, respectively. Also, this map is compatible with the action of $N_\nu \times N_\nu$ on $G_\nu \times \mathcal{J}_N^{-1}(\nu)$ via $(n, n') \cdot (g, z) := (gn^{-1}, n' \cdot z)$, and the $N_\nu$–action on $\mathcal{J}_N^{-1}(\nu)$. Indeed, for any $(n, n') \in N_\nu \times N_\nu$ and any $(g, z) \in G_\nu \times \mathcal{J}_N^{-1}(\nu)$, the point $(gn^{-1}, n' \cdot z)$ gets sent by this map to $gn^{-1}n' \cdot z$. As $N_\nu$ is normal in $G_\nu$ there exists some $n'' \in N_\nu$ such that $gn^{-1}n' \cdot z = n''g \cdot z$ which is in the same $N_\nu$–orbit as $g \cdot z$. Consequently, the map $\varphi^\nu : G_\nu \times \mathcal{J}_N^{-1}(\nu) \to \mathcal{J}_N^{-1}(\nu)$ drops to a smooth map $\phi^\nu : G_\nu/N_\nu \times \mathcal{J}_N^{-1}(\nu)/N_\nu \to \mathcal{J}_N^{-1}(\nu)/N_\nu$ that coincides with (4.3) and therefore satisfies that $\phi^\nu_{kN_\nu} \circ \pi_\nu = \pi_\nu \circ \varphi^\nu_k$, for any $kN_\nu \in H_\nu$.

We now show that the action given by the map $\phi^\nu$ is canonical. Let $kN_\nu \in H_\nu$, $m \in \mathcal{J}_N^{-1}(\nu)$, and $f, h \in C^\infty(M)^N$ arbitrary. Then, taking into account that $\phi^\nu_{kN_\nu} \circ \pi_\nu = \pi_\nu \circ \varphi^\nu_k$ and that by part **(i)** in Proposition 4.1 the functions $f \circ \Phi_{k^{-1}}$ and $h \circ \Phi_{k^{-1}}$ are $N$–invariant, we can write:

$$\begin{aligned}
\pi_\nu^*((\phi^\nu_{kN_\nu})^*\omega_\nu)(m)(X_f(m), X_h(m)) &= ((\phi^\nu_{kN_\nu} \circ \pi_\nu)^*\omega_\nu)(m)(X_f(m), X_h(m)) \\
&= ((\pi_\nu \circ \varphi^\nu_k)^*\omega_\nu)(m)(X_f(m), X_h(m)) \\
&= (\varphi^\nu_k)^*(\pi_\nu^*\omega_\nu)(m)(X_f(m), X_h(m)) \\
&= \pi_\nu^*\omega_\nu(k \cdot m)(T_m\Phi_k \cdot X_f(m), T_m\Phi_k \cdot X_h(m)) \\
&= \pi_\nu^*\omega_\nu(k \cdot m)(X_{f\circ\Phi_{k^{-1}}}(k \cdot m), X_{h\circ\Phi_{k^{-1}}}(k \cdot m)) \\
&= \{f \circ \Phi_{k^{-1}}, h \circ \Phi_{k^{-1}}\}(k \cdot m) = \{f, h\}(m) \\
&= \pi_\nu^*\omega_\nu(m)(X_f(m), X_h(m)).
\end{aligned}$$

Since the map $\pi_\nu$ is a surjective submersion, this chain of equalities implies that $(\phi^\nu_{kN_\nu})^*\omega_\nu = \omega_\nu$, as required.

**(ii)** Let $m' \in \mathcal{J}_G^{-1}(\rho)$ be such that $m' \neq m$. Then, there exists $\mathcal{F}_T \in G_{\mathcal{A}'_G}$ such that $m' = \mathcal{F}_T(m)$. For simplicity in the exposition take $\mathcal{F}_T = F_T$, with $F_T$ the Hamiltonian flow associated to the $G$–invariant function $f \in C^\infty(M)^G$. Let now $f_\nu \in C^\infty(M_\nu)^{H_\nu}$ be the $H_\nu$–invariant function on $M_\nu$ uniquely determined by the relation $f_\nu \circ \pi_\nu = f \circ i_\nu$. The Hamiltonian flow $F^\nu_T$ associated to $f_\nu$ is related to $F_T$ by the relation $F^\nu_T \circ \pi_\nu = \pi_\nu \circ F_T \circ i_\nu$. Therefore, by Noether's Theorem applied to $\mathcal{J}_{H_\nu}$ we have that:

$$\mathcal{J}_{H_\nu}(\pi_\nu(m')) = \mathcal{J}_{H_\nu}(\pi_\nu(F_T(m))) = \mathcal{J}_{H_\nu}(F^\nu_T(\pi_\nu(m))) = \mathcal{J}_{H_\nu}(\pi_\nu(m)) = \sigma,$$

as required. ■

## 4.3 The optimal reduction by stages theorem

Let $m \in M$ be such that $\rho = \mathcal{J}_G(m)$. Also, let $\nu = \mathcal{J}_N(m)$ and $\sigma = \mathcal{J}_{H_\nu}(\pi_\nu(m))$. The second part of Proposition 4.5 guarantees that the restriction of $\pi_\nu$ to $\mathcal{J}_G^{-1}(\rho)$ gives us a well defined map

$$\pi_\nu|_{\mathcal{J}_G^{-1}(\rho)} : \mathcal{J}_G^{-1}(\rho) \longrightarrow \mathcal{J}_{H_\nu}^{-1}(\sigma).$$

This map is smooth because $\mathcal{J}_{H_\nu}^{-1}(\sigma)$ is an initial submanifold of $M_\nu$ and also because $\mathcal{J}_G^{-1}(\rho)$ and $\mathcal{J}_N^{-1}(\nu)$ are initial submanifolds of $M$ such that $\mathcal{J}_G^{-1}(\rho) \subset \mathcal{J}_N^{-1}(\nu) \subset M$, which implies that $\mathcal{J}_G^{-1}(\rho)$ is an initial submanifold of $\mathcal{J}_N^{-1}(\nu)$ (this argument is a straightforward consequence of the definition of initial submanifold). Denote by $i_{\rho,\nu} : \mathcal{J}_G^{-1}(\rho) \hookrightarrow \mathcal{J}_N^{-1}(\nu)$ the corresponding smooth injection. Let $(H_\nu)_\sigma$ be the $H_\nu$–isotropy subgroup of the element $\sigma \in M_\nu/A'_{H_\nu}$. Then, the map $\pi_\nu|_{\mathcal{J}_G^{-1}(\rho)} = \pi_\nu \circ i_{\rho,\nu} : \mathcal{J}_G^{-1}(\rho) \to \mathcal{J}_{H_\nu}^{-1}(\sigma)$ is smooth and $(G_\rho, (H_\nu)_\sigma)$–equivariant. Indeed, let $g \in G_\rho$ and $m \in \mathcal{J}_G^{-1}(\rho)$ arbitrary. By Lemma 4.4 we know that as $G_\rho \subset G_\nu$, then $g \in G_\nu$ and $gN_\nu \in G_\nu/N_\nu$. Using Definition 4.3, we have that $\pi_\nu(g \cdot m) = gN_\nu \cdot \pi_\nu(m)$. Additionally, by (4.4) we have that

$$\mathcal{J}_{H_\nu}(gN_\nu \cdot \pi_\nu(m)) = \mathcal{J}_{H_\nu}(\pi_\nu(g \cdot m)) = \sigma,$$

because $g \cdot m \in \mathcal{J}_G^{-1}(\rho)$, which shows that $gN_\nu \in (H_\nu)_\sigma$ and therefore guarantees the $(G_\rho, (H_\nu)_\sigma)$–equivariance of $\pi_\nu|_{\mathcal{J}_G^{-1}(\rho)}$. Consequently, the map $\pi_\nu|_{\mathcal{J}_G^{-1}(\rho)}$ drops to a well defined map $F$ that makes the following diagram

$$\begin{array}{ccc} \mathcal{J}_G^{-1}(\rho) & \xrightarrow{\pi_\nu|_{\mathcal{J}_G^{-1}(\rho)}} & \mathcal{J}_{H_\nu}^{-1}(\sigma) \\ {\scriptstyle \pi_\rho}\downarrow & & \downarrow{\scriptstyle \pi_\sigma} \\ \mathcal{J}_G^{-1}(\rho)/G_\rho & \xrightarrow{F} & \mathcal{J}_{H_\nu}^{-1}(\sigma)/(H_\nu)_\sigma. \end{array}$$

commutative. We remind the reader once more that the $G_\rho$ and $(H_\nu)_\sigma$–actions on $\mathcal{J}_G^{-1}(\rho)$ and $\mathcal{J}_{H_\nu}^{-1}(\sigma)$, respectively are not automatically proper as a consequence of the properness of the $G$–action on $M$. If that happens to be the case, the map $F$ is smooth. Moreover, in that situation Theorem 3.1 guarantees that the quotients $M_\rho := \mathcal{J}_G^{-1}(\rho)/G_\rho$ and $(M_\nu)_\sigma := \mathcal{J}_{H_\nu}^{-1}(\sigma)/(H_\nu)_\sigma$ are symplectic manifolds. We will refer to the symplectic manifold $(\mathcal{J}_{H_\nu}^{-1}(\sigma)/(H_\nu)_\sigma, \omega_\sigma)$ as the ***second stage reduced space***. Recall that the symplectic form $\omega_\sigma$ is uniquely determined by the equality $\pi_\sigma^*\omega_\sigma = i_\sigma^*\omega_\nu$, where $i_\sigma : \mathcal{J}_{H_\nu}^{-1}(\sigma) \hookrightarrow \mathcal{J}_N^{-1}(\nu)/N_\nu$ is the injection and $\pi_\sigma : \mathcal{J}_{H_\nu}^{-1}(\sigma) \to \mathcal{J}_{H_\nu}^{-1}(\sigma)/(H_\nu)_\sigma$ the projection.

Our goal in this section will consist of proving a theorem that under certain hypotheses states that the map $F$ is a symplectomorphism between the one–shot reduced space $(\mathcal{J}_G^{-1}(\rho)/G_\rho, \omega_\rho)$ and the reduced space in two shots $(\mathcal{J}_{H_\nu}^{-1}(\sigma)/(H_\nu)_\sigma, \omega_\sigma)$.

Given that the properness assumptions appear profusely we will simplify the exposition by grouping them all in the following definition.

**Definition 4.6** *Let $(M, \{\cdot,\cdot\})$ be a Poisson manifold acted properly and canonically upon by a Lie group $G$ via the map $\Phi : G \times M \to M$. Let $N$ be a closed normal Lie subgroup of $G$. Let $\rho \in M/A'_G$, $\nu := \pi_N(\rho)$, $H_\nu := G_\nu/N_\nu$, and $\sigma = \mathcal{J}_{H_\nu}\left(\pi_\nu(\mathcal{J}_G^{-1}(\rho))\right) \in M_\nu/A'_{H_\nu}$. We will say that* ***we have proper actions at*** *$\rho$ whenever $G_\rho$ acts properly on $\mathcal{J}_G^{-1}(\rho)$, $N_\nu$ acts properly on $\mathcal{J}_N^{-1}(\nu)$, $H_\nu$ acts properly on $\mathcal{J}_N^{-1}(\nu)/N_\nu$, and $(H_\nu)_\sigma$ acts properly on $\mathcal{J}_{H_\nu}^{-1}(\sigma)$.*

*Let $(G_\nu)_\sigma \subset G_\nu$ be the unique subgroup of $G_\nu$ such that $(H_\nu)_\sigma = (G_\nu)_\sigma/N_\nu$. We say that the element $\rho \in M/A'_G$ satisfies the* ***stages hypothesis*** *when for any other element $\rho' \in M/A'_G$ such that*

$$\pi_N(\rho) = \pi_N(\rho') =: \nu \quad \text{and} \quad \mathcal{J}_{H_\nu}(\pi_\nu(\mathcal{J}_G^{-1}(\rho))) = \mathcal{J}_{H_\nu}(\pi_\nu(\mathcal{J}_G^{-1}(\rho'))) = \sigma$$

*there exists an element $h \in (G_\nu)_\sigma$ such that $\rho' = h \cdot \rho$.*

*We say that the element $\nu \in M/A'_N$ has the* ***local extension property*** *when any function $f \in C^\infty(\mathcal{J}_N^{-1}(\nu))^{G_\nu}$ is such that for any $m \in \mathcal{J}_N^{-1}(\nu)$ there is an open $N$–invariant neighborhood $U$ of $m$ and a function $F \in C^\infty(M)^G$ such that $F|_U = f|_U$.*

**Theorem 4.7 (Optimal Reduction by Stages)** *Let $(M, \{\cdot,\cdot\})$ be a Poisson manifold acted properly and canonically upon by a Lie group $G$ via the map $\Phi : G \times M \to M$. Let $N$ be a closed normal Lie subgroup of $G$. Let $\rho \in M/A'_G$, $\nu := \pi_N(\rho)$, $H_\nu := G_\nu/N_\nu$, and $\sigma = \mathcal{J}_{H_\nu}\left(\pi_\nu(\mathcal{J}_G^{-1}(\rho))\right) \in M_\nu/A'_{H_\nu}$. Then,* *if $\rho$ satisfies the stages hypothesis, we have proper actions at $\rho$, and the quotient manifold $\mathcal{J}_G^{-1}(\rho)/G_\rho$ is either Lindelöf or paracompact, the map*

$$\begin{array}{rccc} F: & \left(\mathcal{J}_G^{-1}(\rho)/G_\rho, \omega_\rho\right) & \longrightarrow & \left(\mathcal{J}_{H_\nu}^{-1}(\sigma)/(H_\nu)_\sigma, \omega_\sigma\right) \\ & \pi_\rho(m) & \longmapsto & \pi_\sigma(\pi_\nu(m)) \end{array}$$

*is a symplectomorphism between the one shot reduced space $(\mathcal{J}_G^{-1}(\rho)/G_\rho, \omega_\rho)$ and $(\mathcal{J}_{H_\nu}^{-1}(\sigma)/(H_\nu)_\sigma, \omega_\sigma)$ that was obtained by reduction in two stages.*

**Proof of the theorem.** $F$ **is injective**: let $\pi_\rho(m)$ and $\pi_\rho(m') \in \mathcal{J}_G^{-1}(\rho)/G_\rho$ be such that $F(\pi_\rho(m)) = F(\pi_\rho(m'))$. By the definition of $F$ this implies that $\pi_\sigma(\pi_\nu(m)) = \pi_\sigma(\pi_\nu(m'))$. Hence, there exists an element $gN_\nu \in (H_\nu)_\sigma$ such that $\pi_\nu(m') = gN_\nu \cdot \pi_\nu(m)$ which, by the definition (4.3), is equivalent to $\pi_\nu(m') = \pi_\nu(g \cdot m)$. Therefore, there exists a $n \in N_\nu$ such that $m' = ng \cdot m$. However, since both $m$ and $m'$ sit in $\mathcal{J}_G^{-1}(\rho)$ we have that $ng \in G_\rho$, necessarily and, consequently $\pi_\rho(m) = \pi_\rho(m')$, as required.

$F$ **is surjective**: let $\pi_\sigma(\bar{z}) \in (M_\nu)_\sigma = \mathcal{J}_{H_\nu}^{-1}(\sigma)/(H_\nu)_\sigma$. Take any $z \in \mathcal{J}_N^{-1}(\nu)$ such that $\pi_\nu(z) = \bar{z}$ and let $\rho' := \mathcal{J}_G(z)$. It is clear that $\pi_N(\rho') = \pi_N(\mathcal{J}_G(z)) = \mathcal{J}_N(z) = \nu = \pi_N(\rho)$ and also, as $\mathcal{J}_{H_\nu}(\pi_\nu(z)) = \sigma$, Lemma 4.5 guarantees that $\mathcal{J}_{H_\nu}(\pi_\nu(\mathcal{J}_G^{-1}(\rho'))) = \sigma = \mathcal{J}_{H_\nu}(\pi_\nu(\mathcal{J}_G^{-1}(\rho)))$. By the stages hypothesis, there exists $h \in (G_\nu)_\sigma$ such that $\rho' = h \cdot \rho$. Now, we have that

$$F(\pi_\rho(h^{-1} \cdot z)) = \pi_\sigma(\pi_\nu(h^{-1} \cdot z)) = \pi_\sigma(h^{-1}N_\nu \cdot \pi_\nu(z)) = \pi_\sigma(\pi_\nu(z)) = \pi_\sigma(\bar{z}),$$

which proves the surjectivity of $F$.

$F$ **is a symplectic map**: we will show that $F^*\omega_\sigma = \omega_\rho$. Let $m \in \mathcal{J}_G^{-1}(\rho)$ and $f, g \in C^\infty(M)^G$ arbitrary. Then,

$$\begin{aligned} \pi_\rho^*(F^*\omega_\sigma)(m)(X_f(m), X_g(m)) &= (F \circ \pi_\rho)^*\omega_\sigma(m)(X_f(m), X_g(m)) \\ &= \left(\pi_\sigma \circ \pi_\nu|_{\mathcal{J}_G^{-1}(\rho)}\right)^* \omega_\sigma(m)(X_f(m), X_g(m)) \\ &= (\pi_\sigma \circ \pi_\nu \circ i_{\rho,\nu})^*\omega_\sigma(m)(X_f(m), X_g(m)) \\ &= ((\pi_\nu \circ i_{\rho,\nu})^*(\pi_\sigma^*\omega_\sigma))(m)(X_f(m), X_g(m)) \\ &= ((\pi_\nu \circ i_{\rho,\nu})^*(i_\sigma^*\omega_\nu))(m)(X_f(m), X_g(m)) \\ &= (i_\sigma \circ \pi_\nu \circ i_{\rho,\nu})^*\omega_\nu(m)(X_f(m), X_g(m)) \\ &= \pi_\nu^*\omega_\nu(m)(X_f(m), X_g(m)) = \{f, g\}(m) \\ &= \pi_\rho^*\omega_\rho(m)(X_f(m), X_g(m)). \end{aligned}$$

This chain of equalities guarantees that $\pi_\rho^*(F^*\omega_\sigma) = \pi_\rho^*\omega_\rho$. Since the map $\pi_\rho$ is a surjective submersion we have that $F^*\omega_\sigma = \omega_\rho$, and consequently $F$ is a symplectic map.

$F$ **is a symplectomorphism**: given that $F$ is a bijective symplectic map, it is necessarily an immersion. Since by hypothesis the space $\mathcal{J}_G^{-1}(\rho)/G_\rho$ is either Lindelöf or paracompact, a standard result in manifolds theory guarantees that $F$ is actually a diffeomorphism. ■

**Proposition 4.8** *Let $(M, \{\cdot, \cdot\})$ be a Poisson manifold acted properly and canonically upon by a Lie group $G$ via the map $\Phi: G \times M \to M$. Let $N$ be a closed normal Lie subgroup of $G$. Let $\rho \in M/A'_G$, $\nu := \pi_N(\rho)$, $H_\nu := G_\nu/N_\nu$, and $\sigma = \mathcal{J}_{H_\nu}\left(\pi_\nu(\mathcal{J}_G^{-1}(\rho))\right) \in M_\nu/A'_{H_\nu}$. If $\nu$ has the local extension property and $N_\nu$ acts properly on $\mathcal{J}_N^{-1}(\nu)$, then $\pi_\nu(\mathcal{J}_G^{-1}(\rho)) = \mathcal{J}_{H_\nu}^{-1}(\sigma)$ and $\rho$ satisfies the stages hypothesis.*

**Proof.** The inclusion $\pi_\nu(\mathcal{J}_G^{-1}(\rho)) \subset \mathcal{J}_{H_\nu}^{-1}(\sigma)$ is guaranteed by (4.4). In order to prove the equality take $\pi_\nu(m) \in \pi_\nu(\mathcal{J}_G^{-1}(\rho))$ and $f \in C^\infty(M_\nu)^{H_\nu}$ arbitrary, such that the Hamiltonian vector field $X_f$ on $M_\nu$ has flow $F_t$. Let $\bar{f} \in C^\infty(\mathcal{J}_N^{-1}(\nu))^{G_\nu}$ be the function defined by $\bar{f} := f \circ \pi_\nu$. The $H_\nu$–invariance of $f$ implies that $\bar{f}$ is $G_\nu$–invariant. In principle, the point $F_T(\pi_\nu(m))$ lies somewhere in $\mathcal{J}_{H_\nu}^{-1}(\sigma)$. However, we will show that it actually stays in $\pi_\nu(\mathcal{J}_G^{-1}(\rho))$, which will prove the desired equality. Indeed, as the curve $\{F_t(\pi_\nu(m))\}_{t\in[0,T]}$ is compact it can be covered by a finite number of open sets $\{U_1, \ldots, U_n\}$. Suppose that we have chosen the neighborhoods $U_i$ such that $\pi_\nu(m) \in U_1$, $F_T(\pi_\nu(m)) \in U_n$, $U_i \cap U_j \neq \emptyset$ iff $|j-i| = 1$, and for each open $N$–invariant set $\pi_\nu^{-1}(U_i)$, there is a $g_i \in C^\infty(M)^G$ such that $\bar{f}|_{\pi_\nu^{-1}(U_i)} = g_i|_{\pi_\nu^{-1}(U_i)}$. where the function $\bar{f}$ admits local extensions to $G$–invariant functions on $M$. We call $G_t^i$ the flow of the Hamiltonian vector field $X_{g_i}$ on $M$ associated to $g_i \in C^\infty(M)^G$. The flows $G_t^i$ and $F_t$ are related by the equality $F_t \circ \pi_\nu|_{\mathcal{J}_N^{-1}(\nu)\cap\pi_\nu^{-1}(U_i)} = \pi_\nu \circ G_t^i \circ i_\nu|_{\mathcal{J}_N^{-1}(\nu)\cap\pi_\nu^{-1}(U_i)}$. Due to the $G$–invariance of $g$ we have that $\mathcal{J}_G \circ G_t^i = \mathcal{J}_G$ and, consequently $\{F_t(\pi_\nu(m))\}_{t\in[0,T]} \subset \pi_\nu(\mathcal{J}_G^{-1}(\rho))$, as required. This proves that $\pi_\nu(\mathcal{J}_G^{-1}(\rho)) = \mathcal{J}_{H_\nu}^{-1}(\sigma)$.

We conclude by showing that this equality implies that $\rho$ satisfies the stages hypothesis. Indeed, if $\rho' \in M/A'_G$ is such that $\mathcal{J}_{H_\nu}(\pi_\nu(\mathcal{J}_G^{-1}(\rho'))) = \sigma$, then $\pi_\nu(\mathcal{J}_G^{-1}(\rho')) \subset \mathcal{J}_{H_\nu}^{-1}(\sigma) = \pi_\nu(\mathcal{J}_G^{-1}(\rho))$. Consequently, for any $\pi_\nu(z') \in \pi_\nu(\mathcal{J}_G^{-1}(\rho'))$, $z' \in \mathcal{J}_G^{-1}(\rho')$, there exists an element $z \in \mathcal{J}_G^{-1}(\rho)$ such that $\pi_\nu(z') = \pi_\nu(z)$. Hence, there is an element $n \in N_\nu \subset (G_\nu)_\sigma$ available such that $z' = n \cdot z$ which, by applying the map $\mathcal{J}_G$ to both sides of this equality implies that $\rho' = n \cdot \rho$. ■

## 4.4 Reduction by stages of globally Hamiltonian actions on symplectic manifolds

In this section we will assume that $M$ is a symplectic manifold and that the $G$–action is proper and canonical, has a standard $\mathfrak{g}^*$–valued equivariant momentum map $\mathbf{J}_G : M \to \mathfrak{g}^*$, and that, as usual, it contains a closed normal subgroup $N \subset G$. Recall that the inclusion $N \subset G$ and the normal character of $N$ in $G$ implies that $\mathfrak{n}$ is an ideal in $\mathfrak{g}$. Let $i : \mathfrak{n} \hookrightarrow \mathfrak{g}$ be the inclusion. As a corollary to these remarks, it is easy to conclude that the $N$–action on $M$ is also globally Hamiltonian with a $G$–equivariant momentum map $\mathbf{J}_N : M \to \mathfrak{n}^*$ given by $\mathbf{J}_N = i^*\mathbf{J}_G$.

When the $G$–action on $M$ is free, symplectic reduction by stages has been studied in [MMPR98, MMOPR02]. In the following pages we will see how our understanding of the optimal reduction by stages procedure allows us to generalize the results in those papers to the non free actions case. More specifically, we will see that the reduced spaces and subgroups involved in the Optimal Reduction by Stages Theorem 4.7 admit in this case a very precise characterization in terms of level sets of the standard momentum maps present in the problem, and of various subgroups of $G$ obtained as a byproduct of isotropy subgroups related to the $G$ and $N$–actions on $M$ and the coadjoint actions on $\mathfrak{g}^*$ and $\mathfrak{n}^*$.

We start our study by looking in this setup at the level sets of the $G$ and $N$–optimal momentum maps. A basic property of the optimal momentum map whose proof can be found in [OR02a], establishes the following characterization: let $m \in M$ be such that $\mathcal{J}_G(m) = \rho$, $\mathbf{J}_G(m) = \mu$, and $G_m =: H$. Then, $\mathcal{J}_G^{-1}(\rho)$ equals the unique connected component of the submanifold $\mathbf{J}_G^{-1}(\mu) \cap M_H$ that contains it. Analogously, if $\mathcal{J}_N(m) = \nu$, $\mathbf{J}_N(m) = \eta$, and $N_m = H \cap N$, then $\mathcal{J}_N^{-1}(\nu)$ equals the unique connected component of the submanifold $\mathbf{J}_N^{-1}(\eta) \cap M_{H\cap N}$ that contains it. Recall that the symbol $M_H$ denotes the ***isotropy type submanifold*** associated to the isotropy subgroup $H$ and that it is defined by $M_H := \{z \in M \mid G_z = H\}$. All along this section we will assume the following

**Connectedness hypothesis:** the submanifolds $\mathbf{J}_G^{-1}(\mu) \cap M_H$ and $\mathbf{J}_N^{-1}(\eta) \cap M_{H\cap N}$ are connected.

This hypothesis is NOT realistic however it will make the presentation that follows much more clear and accessible. The reduction by stages problem does not differ much, qualitatively speaking, no matter if we assume the connectedness hypothesis or not, however the necessary additions in the notation to accommodate the most general case would make the following pages very difficult to read. In order to adapt to the general situation our results, the reader should just take the relevant connected components of $\mathbf{J}_G^{-1}(\mu) \cap M_H$ and $\mathbf{J}_N^{-1}(\eta) \cap M_{H\cap N}$, and each time that we quotient them by a group that leaves them invariant, the reader should take the closed subgroup that leaves invariant the connected component that he has previously chosen. The notation becomes immediately rather convoluted but the ideas involved in the process are the same.

We continue our characterization of the ingredients for reduction by stages in the following proposition.

**Proposition 4.9** *Let $(M, \omega)$ be a symplectic manifold acted properly and canonically upon by a Lie group $G$ and suppose that this action has a standard equivariant momentum map $\mathbf{J}_G : M \to \mathfrak{g}^*$ associated. Let $N \subset G$ be a closed normal subgroup of $G$. Then, if $m \in M$ is such that $\mathcal{J}_G(m) = \rho$, $\mathbf{J}_G(m) = \mu$, and its isotropy subgroup $G_m$ equals $G_m =: H \subset G$ we have that*

**(i)** $\mathcal{J}_G^{-1}(\rho) = \mathbf{J}_G^{-1}(\mu) \cap M_H$.

**(ii)** $\mathcal{J}_N(m) = \pi_N(\rho) =: \nu$, $\mathbf{J}_N(m) = i^*\mu =: \eta$, *and* $\mathcal{J}_N^{-1}(\nu) = \mathbf{J}_N^{-1}(\eta) \cap M_{N_\eta \cap H}$.

**(iii)** $G_\rho = N_{G_\mu}(H)$, $N_\nu = N_{N_\eta}(N_\eta \cap H)$, *and* $G_\nu = N_{G_\eta}(N_\eta \cap H)$. *The symbol* $N_{G_\mu}(H) := N(H) \cap G_\mu$ *where $N(H)$ denotes the normalizer of $H$ in $G$. We will refer to $N_{G_\mu}(H)$ as the normalizer of $H$ in $G_\mu$.*

**Proof.** The proof of the equalities $\mathcal{J}_G^{-1}(\rho) = \mathbf{J}_G^{-1}(\mu) \cap M_H$ and $G_\rho = N_{G_\mu}(H)$ can be found in [OR02a].

We now show that $\mathcal{J}_N^{-1}(\nu) = \mathbf{J}_N^{-1}(\eta) \cap M_{N_\eta \cap H}$. By the results in [OR02a], it suffices to show that $N_m = N_\eta \cap H$. Indeed, as the $G$–equivariance of $\mathbf{J}_N$ implies that $H = G_m \subset G_\eta$, we have that $N_m = H \cap N = H \cap G_\eta \cap N = N_\eta \cap H$. Consequently, the same result in [OR02a] that gave us $G_\rho = N_{G_\mu}(H)$, can be applied to the $N$–action on $M$ to obtain $N_\nu = N_{N_\eta}(N_\eta \cap H)$.

Finally, we prove the identity $G_\nu = N_{G_\eta}(N_\eta \cap H)$ by double inclusion. Let first $g \in G_\nu$. The equality $g \cdot \nu = \nu$ implies that $g \cdot m = \mathcal{F}_T(m)$, with $\mathcal{F}_T \in G_{\mathcal{A}'_N}$. For simplicity suppose that $\mathcal{F}_T = F_t$, with $F_t$ the Hamiltonian flow associated to a $N$–invariant function on $M$. The standard Noether's Theorem implies that $g \cdot m = F_t(m) \in \mathbf{J}_N^{-1}(\eta)$ and therefore $g \in G_\eta$. Also, as the flow $F_t$ is $N$–equivariant we have that

$$N_\eta \cap H = N_m = N_{F_t(m)} = N_{g\cdot m} = gN_m g^{-1} = g(N_\eta \cap H)g^{-1},$$

and consequently $g \in N_{G_\eta}(N_\eta \cap H)$. The reverse inclusion is trivial. ■

**Remark 4.10** A major consequence of the previous proposition is the fact that the subgroups $G_\nu$ and $N_\nu$, and those that will derive from them, are automatically closed subgroups. This circumstance implies that the proper actions hypothesis given in Definition 4.6 and necessary for reduction by stages is automatically satisfied in this setup. ♦

The previous proposition allows us to explicitly write down in our setup the one–shot reduced space:

$$M_\rho := \mathcal{J}_G^{-1}(\rho)/G_\rho = \mathbf{J}_G^{-1}(\mu) \cap M_H / N_{G_\mu}(H), \tag{4.5}$$

as well as the first stage reduced space:

$$M_\nu := \mathcal{J}_N^{-1}(\nu)/N_\nu = \mathbf{J}_N^{-1}(\eta) \cap M_{N_\eta \cap H} / N_{N_\eta}(N_\eta \cap H).$$

We now proceed with the construction of the second stage reduced space. As it was already the case in the general optimal setup, the quotient group

$$\mathcal{H}_\nu := G_\nu/N_\nu = \frac{N_{G_\eta}(N_\eta \cap H)}{N_{N_\eta}(N_\eta \cap H)}$$

acts canonically on the quotient $M_\nu$ with optimal momentum map associated $\mathcal{J}_{\mathcal{H}_\nu} : M_\nu \to M_\nu/A'_{\mathcal{H}_\nu}$. In this setup we can say more. Indeed, in this case the $\mathcal{H}_\nu$–action on $M_\nu$ is automatically proper and has a standard momentum map associated $\mathbf{J}_{\mathcal{H}_\nu} : M_\nu \to \mathrm{Lie}(\mathcal{H}_\nu)^*$, where the symbol $\mathrm{Lie}(\mathcal{H}_\nu)$ denotes the Lie algebra of the group $\mathcal{H}_\nu$. An explicit expression for $\mathbf{J}_{\mathcal{H}_\nu}$ can be obtained by mimicking the computations made in [MMPR98, MMOPR02] for the free case. In order to write it down we introduce the following maps: let $\pi_{G_\nu} : G_\nu \to G_\nu/N_\nu$ be the projection, $r_\nu = T_e\pi_{G_\nu} : \mathfrak{g}_\nu \to \mathrm{Lie}(\mathcal{H}_\nu) \simeq \mathfrak{g}_\nu/\mathfrak{n}_\nu$ be its derivative at the identity, and $r_\nu^* : \mathrm{Lie}(\mathcal{H}_\nu)^* \to \mathfrak{g}_\nu^*$ be the corresponding dual map. Then, for any $\pi_\nu(z) \in M_\nu$ and any $r_\nu(\xi) \in \mathrm{Lie}(\mathcal{H}_\nu)$, the momentum map $\mathbf{J}_{\mathcal{H}_\nu}$ is given by the expression

$$\langle \mathbf{J}_{\mathcal{H}_\nu}(\pi_\nu(z)), r_\nu(\xi)\rangle = \langle \mathbf{J}_G(z), \xi\rangle - \langle \bar{\eta}, \xi\rangle, \tag{4.6}$$

where $\bar{\eta} \in \mathfrak{g}_\nu^*$ is some chosen extension of the restriction $\eta|_{\mathfrak{n}_\nu}$ to a linear functional on $\mathfrak{g}_\nu$. This momentum map is not equivariant. Indeed, its non equivariance cocycle $\bar{\omega}$ is given by the expression

$$r_\nu^*(\bar{\omega}(\pi_{G_\nu}(h))) = \mathrm{Ad}^*_{h^{-1}}\bar{\eta} - \bar{\eta},$$

for any $\pi_{G_\nu}(h) \in G_\nu/N_\nu$. The map $\mathbf{J}_{\mathcal{H}_\nu}$ becomes equivariant if we replace the coadjoint action of $\mathcal{H}_\nu$ on the dual of its Lie algebra by the affine action defined by

$$\pi_{G_\nu}(h) \cdot \lambda := \mathrm{Ad}^*_{(\pi_{G_\nu}(h))^{-1}}\lambda + \bar{\omega}(\pi_{G_\nu}(h)), \tag{4.7}$$

for any $\pi_{G_\nu}(h) \in \mathcal{H}_\nu$ and any $\lambda \in \mathrm{Lie}(\mathcal{H}_\nu)^*$. Let now $\tau \in \mathrm{Lie}(\mathcal{H}_\nu)^*$ be the element defined by

$$\langle \tau, r_\nu(\xi)\rangle = \langle \mu, \xi\rangle - \langle \bar{\nu}, \xi\rangle, \tag{4.8}$$

for any $r_\nu(\xi) \in \mathrm{Lie}(\mathcal{H}_\nu)$. A calculation following the lines of [MMPR98, MMOPR02] shows that the isotropy subgroup $(\mathcal{H}_\nu)_\tau$ of $\tau$ with respect to the affine action (4.7) of $\mathcal{H}_\nu$ on the dual of its Lie algebra, is given by

$$(\mathcal{H}_\nu)_\tau = \pi_{G_\nu}\left((G_\nu)_{\mu|_{\mathfrak{g}_\nu}}\right) = \frac{\left(\left(N_{G_\eta}(N_\eta \cap H)\right)_{\mu|_{\mathrm{Lie}\left(N_{G_\eta}(N_\eta\cap H)\right)}}\right)}{N_{N_\eta}(N_\eta \cap H)} \tag{4.9}$$

Now, for any $m \in \mathcal{J}_G^{-1}(\rho)$, the choice of $\tau \in \mathrm{Lie}(\mathcal{H}_\nu)^*$ in (4.8) guarantees that $\mathbf{J}_{\mathcal{H}_\nu}(\pi_\nu(m)) = \tau$ and, moreover, if $\mathcal{J}_{\mathcal{H}_\nu}(\pi_\nu(m)) = \sigma \in M/A'_{\mathcal{H}_\nu}$ then,

$$\mathcal{J}_{\mathcal{H}_\nu}^{-1}(\sigma) = \mathbf{J}_{\mathcal{H}_\nu}^{-1}(\tau) \cap (M_\nu)_{(\mathcal{H}_\nu)_{\pi_\nu(m)}} \tag{4.10}$$

since, by extension of the connectedness hypothesis we will suppose that $\mathbf{J}_{\mathcal{H}_\nu}^{-1}(\tau) \cap (M_\nu)_{(\mathcal{H}_\nu)_{\pi_\nu(m)}}$ is also connected.

We compute the isotropy subgroup $(\mathcal{H}_\nu)_{\pi_\nu(m)}$ in terms of the groups that already appeared in our study. Indeed, we will now show that

$$(\mathcal{H}_\nu)_{\pi_\nu(m)} = \frac{N_{N_\eta}(H \cap N_\eta)H}{N_{N_\eta}(N_\eta \cap H)}. \tag{4.11}$$

Take first an element $\pi_{G_\nu}(g) \in \mathcal{H}_\nu$ such that $\pi_{G_\nu}(g)\cdot\pi_\nu(m) = \pi_\nu(m)$ or, equivalently, $\pi_\nu(g\cdot m) = \pi_\nu(m)$. Hence, there exists a group element $n \in N_\nu = N_{N_\eta}(N_\eta \cap H)$ such that $g \cdot m = n \cdot m$. Given that $G_m = H$ we have that $n^{-1} \cdot g \in H$, necessarily and hence $g \in N_{N_\eta}(N_\eta \cap H)H$ and $\pi_{G_\nu}(g) \in N_{N_\eta}(N_\eta \cap H)H/N_{N_\eta}(N_\eta \cap H)$. Conversely, if $\pi_{G_\nu}(g) \in N_{N_\eta}(N_\eta \cap H)H/N_{N_\eta}(N_\eta \cap H)$, we can write $g = nh$, with $n \in N_{N_\eta}(N_\eta \cap H)$ and $h \in H$ and therefore $\pi_{G_\nu}(g) \cdot \pi_\nu(m) = \pi_\nu(nh \cdot m) = \pi_\nu(n \cdot m) = \pi_\nu(m)$, as required.

In order to write down the second stage reduced space we have to compute the isotropy subgroup $(\mathcal{H}_\nu)_\sigma$. In view of (4.9) and (4.11), and Proposition 4.9 adapted to the optimal momentum map $\mathcal{J}_{\mathcal{H}_\nu}$ we have that

$$(\mathcal{H}_\nu)_\sigma = N_{(\mathcal{H}_\nu)_\tau}\left(\frac{N_{N_\eta}(H \cap N_\eta)H}{N_{N_\eta}(N_\eta \cap H)}\right), \tag{4.12}$$

where the group $(\mathcal{H}_\nu)_\tau$ is given by Expression (4.9). We now recall a standard result about normalizers that says that if $A \subset B \subset C$ are groups such that $A$ is normal in both $B$ and $C$, then

$$N_{C/A}(B/A) = N_C(B)/A.$$

If we apply this equality to Expression (4.12) we obtain that

$$(\mathcal{H}_\nu)_\sigma = \frac{N_{(G_\nu)_{\mu|_{\mathfrak{g}_\nu}}}(N_{N_\eta}(H \cap N_\eta)H)}{N_{N_\eta}(N_\eta \cap H)} = \frac{N_{\left(\left(N_{G_\eta}(N_\eta\cap H)\right)_{\mu|_{\mathrm{Lie}\left(N_{G_\eta}(N_\eta\cap H)\right)}}\right)}(N_{N_\eta}(H \cap N_\eta)H)}{N_{N_\eta}(N_\eta \cap H)} \tag{4.13}$$

All the computations that we just carried out allow us to explicitly write down the second stage reduced space. Namely, by combination of expressions (4.10), (4.11), and (4.13), we obtain that

$$(M_\nu)_\sigma = \mathcal{J}_{\mathcal{H}_\nu}^{-1}(\sigma)/(\mathcal{H}_\nu)_\sigma = \frac{\mathbf{J}_{\mathcal{H}_\nu}^{-1}(\tau) \cap (M_\nu)_{\frac{N_{N_\eta}(H\cap N_\eta)H}{N_{N_\eta}(N_\eta\cap H)}}}{\frac{N_{(G_\nu)_{\mu|_{\mathfrak{g}_\nu}}}(N_{N_\eta}(H\cap N_\eta)H)}{N_{N_\eta}(N_\eta\cap H)}}, \tag{4.14}$$

where the group $(G_\nu)_{\mu|_{\mathfrak{g}_\nu}} = \left(N_{G_\eta}(N_\eta \cap H)\right)_{\mu|_{\mathrm{Lie}\left(N_{G_\eta}(N_\eta\cap H)\right)}}$.

The Optimal Reduction by Stages Theorem 4.7 guarantees that the second stage reduced space (4.14) is symplectomorphic to the one–shot reduced space (4.5) in the presence of the Stages Hypothesis introduced in Definition 4.6. In this setup, that hypothesis can be completely reformulated in terms of relations between Lie algebraic elements and isotropy subgroups. More specifically, in the globally Hamiltonian framework, the Stages Hypothesis is equivalent to the following condition:

**Hamiltonian Stages Hypothesis:** Let $\mu \in \mathfrak{g}^*$ and $H \subset G$. We say that the pair $(\mu, H)$ satisfies the Hamiltonian Stages Hypothesis whenever for any other similar pair $(\mu', H')$ such that

$$\left\{\begin{array}{l} i^*\mu = i^*\mu' =: \eta \in \mathfrak{n}^* \\ N_\eta \cap H = N_\eta \cap H' =: K \end{array}\right. \quad \text{and} \quad \left\{\begin{array}{l} \mu|_{\mathrm{Lie}\left(N_{G_\eta}(K)\right)} = \mu'|_{\mathrm{Lie}\left(N_{G_\eta}(K)\right)} =: \zeta \in \mathrm{Lie}\left(N_{G_\eta}(K)\right)^* \\ N_{N_\eta}(K)H = N_{N_\eta}(K)H' =: L, \end{array}\right.$$

there exists an element $n \in N_{(N_{G_\eta}(K))_\zeta}(L)$ such that

$$\mu' = \mathrm{Ad}^*_{n^{-1}}\mu \quad \text{and} \quad H' = nHn^{-1}.$$

**Remark 4.11** A quick inspection shows that when the $G$–action is free, that is, when all the isotropy subgroups $H = \{e\}$, the previous condition collapses into the Stages Hypothesis introduced in [MMOPR02]. ♦

We recall that, in the same fashion in which the proper actions hypothesis introduced in Definition 4.6 is automatically satisfied in this setup, so is the Lindelöf hypothesis on the one–shot reduced space $M_\rho$ if we just assume that $M$ is Lindelöf. This is so because closed subsets and continuous images of Lindelöf spaces are always Lindelöf.

The Optimal Reduction by Stages Theorem together with the ideas that we just introduced implies in this setup the following highly non trivial symplectomorphism that we enunciate in the form of a theorem. The following statement is consistent with the previously introduced notations.

**Theorem 4.12 (Hamiltonian Reduction by Stages)** *Let $(M, \omega)$ be a symplectic manifold acted properly and canonically upon by a Lie group $G$ that has a closed normal subgroup $N$. Suppose that this action has an equivariant momentum map $\mathbf{J}_G : M \to \mathfrak{g}^*$ associated. Let $\mu \in \mathfrak{g}^*$ be a value of $\mathbf{J}_G$ and $H$ an isotropy subgroup of the $G$ action on $M$. If the manifold $M$ is Lindelöf and the pair $(\mu, H)$ satisfies the Hamiltonian Stages Hypothesis, then the symplectic reduced spaces*

$$\frac{\mathbf{J}_G^{-1}(\mu) \cap M_H}{N_{G_\mu}(H)} \quad \textit{and} \quad \frac{\mathbf{J}_{\mathcal{H}_\nu}^{-1}(\tau) \cap (M_\nu)_{\frac{N_{N_\eta}(H\cap N_\eta)H}{N_{N_\eta}(N_\eta \cap H)}}}{\frac{N_{(G_\nu)_{\mu|_{\mathfrak{g}_\nu}}}(N_{N_\eta}(H\cap N_\eta)H)}{N_{N_\eta}(N_\eta\cap H)}}$$

*are symplectomorphic. In this expression $\eta = i^*\mu$,*

$$M_\nu = \frac{\mathbf{J}_N^{-1}(\eta) \cap M_{N_\eta \cap H}}{N_{N_\eta}(N_\eta \cap H)}, \quad \mathcal{H}_\nu = \frac{N_{G_\eta}(N_\eta \cap H)}{N_{N_\eta}(N_\eta \cap H)}, \quad (G_\nu)_{\mu|_{\mathfrak{g}_\nu}} = \left(N_{G_\eta}(N_\eta \cap H)\right)_{\mu|_{\mathrm{Lie}\left(N_{G_\eta}(N_\eta\cap H)\right)}},$$

*$\mathbf{J}_{\mathcal{H}_\nu} : M_\nu \to \mathrm{Lie}(\mathcal{H}_\nu)^*$ is the momentum map associated to the $\mathcal{H}_\nu$–action on $M_\nu$ defined in (4.6), and $\tau \in \mathrm{Lie}(\mathcal{H}_\nu)^*$ the element defined in (4.8).*

**Remark 4.13** When the $G$–action is free, the previous theorem coincides with the Reduction by Stages Theorem presented in [MMOPR02]. ♦

A special but very important particular case of Theorem 4.12 takes place when the group $G$ is discrete ($\mathfrak{g} = \{0\}$). In that situation, all the standard momentum maps in the construction vanish and the theorem gives us a highly non trivial relation between quotients of isotropy type submanifolds. We start by reformulating the Hamiltonian Stages Hypothesis in this particular case.

**Discrete Reduction by Stages Hypothesis:** Let $G$ be a discrete group, $N$ a normal subgroup, and $H$ a subgroup. We say that $H$ satisfies the Discrete Reduction by Stages Hypothesis with respect to $N$ if for any other subgroup $H'$ such that

$$N \cap H = N \cap H' =: K \quad \text{and} \quad N_N(K)H = N_N(K)H' =: L,$$

there exists an element $n \in N_{N_G(K)}(L)$ such that $H' = nHn^{-1}$.

**Theorem 4.14 (Discrete Reduction by Stages)** *Let $(M, \omega)$ be a symplectic manifold acted properly and canonically upon by a discrete Lie group $G$ that has a closed normal subgroup $N$. Let $H$ be an* *isotropy subgroup of the $G$ action on $M$ and $K := N \cap H$. If the manifold $M$ is Lindelöf and $H$ satisfies the Discrete Reduction by Stages Hypothesis with respect to $N$, then the symplectic reduced spaces*

$$\frac{M_H}{N_G(H)} \quad \text{and} \quad \frac{\left(\frac{M_K}{N_N(K)}\right)_{\frac{N_N(K)H}{N_N(K)}}}{\frac{N_{N_G(K)}(N_N(K)H)}{N_N(K)}} \tag{4.15}$$

*are symplectomorphic.*

**Remark 4.15** When the $G$–action on $M$ is free, the Discrete Reduction by Stages Hypothesis is trivially satisfied and Theorem 4.14 produces the straightforward symplectomorphism

$$M/G \simeq (M/N)/(G/N).$$

Hence, it is in the presence of singularities that the relation stablished in (4.15) is really visible and non trivial. ♦

# 5 Appendix

## 5.1 Proof of Proposition 3.8

**(i)** The distribution $D$ can be written as the span of globally defined vector fields on $M$, that is,

$$D = \operatorname{span}\{\xi_M,\, X_f \mid \xi \in \mathfrak{g} \text{ and } f \in C^\infty(M)^G\}. \tag{5.1}$$

By the Frobenius–Stefan–Sussman Theorem [St74a, St74b, Su73], the integrability of $D$ can be proved by showing that this distribution is invariant by the flows of the vector fields in (5.1) that we used to generate it. Let $f, l \in C^\infty(M)^G$, $\xi, \eta \in \mathfrak{g}$, $F_t$ be the flow of $X_l$, and $H_t$ be the flow of $\eta_M$. Recall that $\eta_M$ is a complete vector field such that $H_t(m) = \exp t\eta \cdot m$, for all $t \in \mathbb{R}$ and $m \in M$. Now, the integrability of $A'_G$ guarantees that $T_mF_t \cdot X_f(m) \in A'_G(F_t(m)) \subset D(F_t(m))$. Also, the $G$–equivariance of $F_t$ and the invariance of the function $f$ imply that $T_mF_t \cdot \xi_M(m) = \xi_M(F_t(m))$ and $T_mH_t \cdot X_f(m) = X_f(H_t(m))$. Finally, we have that

$$\begin{aligned} T_mH_t \cdot \xi_M(m) &= \left.\frac{d}{ds}\right|_{s=0} \exp t\eta \exp s\xi \cdot m \\ &= \left.\frac{d}{ds}\right|_{s=0} \exp t\eta \exp s\xi \exp -t\eta \exp t\eta \cdot m = (\mathrm{Ad}_{\exp t\eta}\xi)_M(\exp t\eta \cdot m), \end{aligned}$$

which proves that $D$ is integrable.

**(ii)** As $D$ is integrable and is generated by the vector fields (5.1), its maximal integral submanifolds coincide with the orbits of the action of the pseudogroup constructed by finite composition of flows of the vector fields in (5.1), that is, for any $m \in M$, the integral leaf $\mathcal{L}_m$ of $D$ that goes through $m$ is:

$$\mathcal{L}_m = \{F_{t_1} \circ \cdots \circ F_{t_n} \mid \text{with } F_{t_i} \text{ the flow of a vector field in (5.1)}\}.$$

Given that $[X_f, \xi_M] = 0$ for all $f \in C^\infty(M)^G$ and $\xi \in \mathfrak{g}$, the previous expression can be rewritten as

$$\mathcal{L}_m = \{H_{t_1} \circ \cdots \circ H_{t_j} \circ G_{s_1} \circ \cdots \circ G_{s_k} \mid G_{s_i} \text{ flow of } f_i \in C^\infty(M)^G, \text{ and } H_{t_i} \text{ flow of } \xi^i_M,\, \xi^i \in \mathfrak{g}\}.$$

Therefore, $\mathcal{L}_m = G^0 \cdot \mathcal{J}^{-1}(\rho)$, as required. ■

## 5.2 Proof of Proposition 3.9

**(i)** It is easy to check that $G_\rho$ is closed in $G$ iff the action of $G_\rho$ on $G$ by right translations is proper. Additionally, if $G_\rho$ is closed in $G$ then the $G_\rho$–action on $\mathcal{J}^{-1}(\rho)$ is proper. In any case, if the action of $G_\rho$ on either $G$, or on $\mathcal{J}^{-1}(\rho)$, or on both, is proper, so is the action on the product $G \times \mathcal{J}^{-1}(\rho)$ in the statement of the proposition. As to the freeness, it is inherited from the freeness of the $G_\rho$–action on $G$.

**(ii)** First of all, the map $i$ is clearly well defined and smooth since it is the projection onto the orbit space $G \times_{G_\rho} \mathcal{J}^{-1}(\rho)$ of the $G_\rho$–invariant smooth map $G \times \mathcal{J}^{-1}(\rho) \to M$ given by $(g, z) \longmapsto g \cdot z$. It is also injective because if $[g, z], [g', z'] \in G \times_{G_\rho} \mathcal{J}^{-1}(\rho)$ are such that $i([g, z]) = i([g', z'])$, then $g \cdot z = g' \cdot z'$ or, analogously, $g^{-1}g' \cdot z' = z$, which implies that $g^{-1}g' \in G_\rho$. Consequently, $[g, z] = [gg^{-1}g', (g')^{-1}g \cdot z] = [g', z']$, as required.

Finally, we check that $i$ is an immersion. Let $[g, z] \in G \times_{G_\rho} \mathcal{J}^{-1}(\rho)$ arbitrary and $\xi \in \mathfrak{g}$, $f \in C^\infty(M)^G$ be such that $T_{[g,z]}i \cdot T_{(g,z)}\pi_{G_\rho} \cdot (T_eL_g(\xi), X_f(z)) = 0$. If we denote by $F_t$ the flow of $X_f$ we can rewrite this equality as

$$\left.\frac{d}{dt}\right|_{t=0} g \exp t\xi \cdot F_t(z) = 0 \quad \text{or equivalently,} \quad T_z\Phi_g(X_f(z) + \xi_M(z)) = 0.$$

Hence $X_f(z) = -\xi_M(z)$ which by (2.7) implies that $\xi \in \mathfrak{g}_\rho$ and therefore $T_{(g,z)}\pi_{G_\rho} \cdot (T_eL_g(\xi), X_f(z)) = T_{(g,z)}\pi_{G_\rho} \cdot (T_eL_g(\xi), -\xi_M(z)) = 0$, as required.

Given that for any $\xi \in \mathfrak{g}$, $f \in C^\infty(M)^G$, and $[g, z] \in G \times_{G_\rho} \mathcal{J}^{-1}(\rho)$ we have that $T_{[g,z]}i \cdot T_{(g,z)}\pi_{G_\rho} \cdot (T_eL_g(\xi), X_f(z)) = (\mathrm{Ad}_g\xi)_M(g \cdot z) + X_f(g \cdot z)$, it is clear that $T_{[g,z]}i \cdot T_{[g,z]}(G \times_{G_\rho} \mathcal{J}^{-1}(\rho)) = D(g \cdot z)$ and thereby $i(G \times_{G_\rho} \mathcal{J}^{-1}(\rho)) = \mathcal{J}^{-1}(\mathcal{O}_\rho)$ is an integral submanifold of $D$. ■

**Acknowledgments** I thank Richard Cushman, Jerry Marsden, James Montaldi, and Tudor Ratiu for their valuable comments and their encouragement regarding this project. This research was partially supported by the European Commission through funding for the Research Training Network *Mechanics and Symmetry in Europe* (MASIE).

# References

[AM78] Abraham, R., and Marsden, J.E. [1978] *Foundations of Mechanics.* Second edition, Addison–Wesley.

[AMR99] Abraham, R., Marsden, J.E., and Ratiu, T.S. [1988] *Manifolds, Tensor Analysis, and Applications.* Volume 75 of *Applied Mathematical Sciences*, Springer-Verlag.

[AMM98] Alekseev, A., Malkin, A., and Meinrenken, E. [1998] Lie group valued momentum maps. *J. Differential Geom.*, 48:445–495.

[ACG91] Arms, J.M., Cushman, R., and Gotay, M.J. [1991] A universal reduction procedure for Hamiltonian group actions. In *The Geometry of Hamiltonian Systems*, (T. S. Ratiu ed.). Pages 33–51. Springer Verlag.

[BL97] Bates, L. and Lerman, E. [1997] Proper group actions and symplectic stratified spaces. *Pacific J. Math.*, 181(2):201–229.

[Bl01] Blaom, A. [2001] *A Geometric Setting for Hamiltonian Perturbation Theory. Memoirs of the American Mathematical Society*, vol. 153, number 727.

[Bre72] Bredon, G.E. [1972] *Introduction to Compact Transformation Groups.* Academic Press.

[CaWe99] Cannas da Silva, A. and Weinstein, A. [1999] *Geometric Models for Noncommutative Algebras.* Berkeley Math. Lecture Notes. Amer. Math. Soc.

[C22] Cartan, É. [1922] Leçons sur les Invariants Intégraux. Hermann.

[CS01] Cushman, R. and Sniatycki, J. [2001] Differential structure of orbit spaces. *Canad. J. Math.*, 53(4):715–755.

[Daz85] Dazord, P. [1985] Feuilletages à singularités. *Nederl. Akad. Wetensch. Indag. Math.*, 47:21–39.

[DuKo99] Duistermaat, J. J. and Kolk, J. A. C. [1999] *Lie Groups*. Springer Verlag.

[KKS78] Kazhdan, D., Kostant, B., and Sternberg, S. [1978] Hamiltonian group actions dynamical systems of Calogero type. *Comm. Pure Appl. Math*, 31:481–508.

[K76] Kirillov, A. A. [1976] Elements of the Theory of Representations. *Grundlehren der mathematischen Wissenschaften*, volume 220. Springer–Verlag.

[Lie90] Lie, S. [1890] *Theorie der Transformationsgruppen. Zweiter Abschnitt.* Teubner.

[LM87] Libermann, P., and Marle, C.–M. [1987] *Symplectic Geometry and Analytical Mechanics.* Reidel.

[Lie90] Lie, S. [1890] *Theorie der Transformationsgruppen. Zweiter Abschnitt.* Teubner.

[Mar76] Marle, C.–M. [1976] Symplectic manifolds, dynamical groups and Hamiltonian mechanics. In *Differential Geometry and Relativity.* M. Cahen, and M. Flato (eds.). Reidel.

[MMPR98] Marsden, J. E., Misiolek, G., Perlmutter, M., and Ratiu, T. S. [1998] Symplectic reduction for semidirect products and central extensions. *Diff. Geom. and Appl.*, 9:173–212.

[MMOPR02] Marsden, J. E., Misiolek, G., Ortega, J.-P., Perlmutter, M., and Ratiu, T. S. [2002] Symplectic reduction by stages. *In preparation.*

[MR86] Marsden, J.E., and Ratiu, T.S. [1986] Reduction of Poisson manifolds. *Letters in Mathematical Physics*, 11:161–169.

[MRW84] Marsden, J.E., Ratiu, T.S.,and Weinstein, A. [1984] Semidirect products and reduction in mechanics. *Trans. A.M.S.,* 281:147–177.

[MW74] Marsden, J.E., and Weinstein, A. [1974] Reduction of symplectic manifolds with symmetry. *Rep. Math. Phys.*, 5(1):121–130.

[McD88] McDuff, D. [1988] The moment map for circle actions on symplectic manifolds. *J. Geom. Phys.*, 5:149–160.

[Mey73] Meyer, K. R. [1973] Symmetries and integrals in mechanics. In *Dynamical Systems*, pp. 259–273. M.M. Peixoto, ed. Academic Press.

[Mi00] Michor, P. W. [2000] *Topics in Differential Geometry. Preprint.*

[O98] Ortega, J.–P. [1998] *Symmetry, Reduction, and Stability in Hamiltonian Systems.* Ph.D. Thesis. University of California, Santa Cruz. June, 1998.

[O02] Ortega, J.–P. [2002] Singular dual pairs. To appear in *Differential Geometry and its Applications. Preprint available at http://arXiv.org/abs/math.SG/0201192.*

[O02a] Ortega, J.–P. [2002] The symplectic reduced spaces of a Poisson action. To appear in *C. R. Acad. Sci. Paris Sér. I Math.*. Available at http://arXiv.org/abs/math/0204154.

[OR98] Ortega, J.–P. and Ratiu, T.S. [1998] Singular reduction of Poisson manifolds. *Letters in Mathematical Physics*, 46:359–372.

[OR02] Ortega, J.–P. and Ratiu, T. S. [2001] A symplectic slice theorem. *Lett. Math. Phys.*, 59:81–93.

[OR02a] Ortega, J.–P. and Ratiu, T. S. [2002] The optimal momentum map. To appear in *Geometry, Dynamics, and Mechanics: 60th Birthday Volume for J.E. Marsden.* P. Holmes, P. Newton, and A. Weinstein, eds., Springer-Verlag. Available at http://arXiv.org/abs/math.SG/0203040.

[OR02b] Ortega, J.–P. and Ratiu, T. S. [2002] *Hamiltonian Singular Reduction. To appear in* Birkhäuser, Progress in Mathematics.

[S90] Sjamaar, R. [1990] *Singular Orbit Spaces in Riemannian and Symplectic Geometry.* Ph. D. thesis, Rijksuniversiteit te Utrecht.

[SL91] Sjamaar, R. and Lerman, E. [1991] Stratified symplectic spaces and reduction. *Ann. of Math.*, 134:375–422.

[St74a] Stefan, P. [1974] Accessibility and foliations with singularities. *Bull. Amer. Math. Soc.*, 80:1142–1145.

[St74b] Stefan, P. [1974] Accessible sets, orbits and foliations with singularities. *Proc. Lond. Math. Soc.*, 29:699–713.

[Su73] Sussman, H. [1973] Orbits of families of vector fields and integrability of distributions. *Trans. Amer. Math. Soc.*, 180:171–188.

[W83] Weinstein, A. [1983] The local structure of Poisson manifolds. *J. Differential Geometry*, 18:523–557.